# Normed Ω-Group

Aleks Kleyn

Aleks_Kleyn@MailAPS.org
http://AleksKleyn.dyndns-home.com:4080/
http://sites.google.com/site/AleksKleyn/
http://arxiv.org/a/kleyn_a_1
http://AleksKleyn.blogspot.com/

ABSTRACT. Since sum which is not necessarily commutative is defined in $\Omega$-algebra $A$, then $\Omega$-algebra $A$ is called $\Omega$-group. I also considered representation of $\Omega$-group. Norm defined in $\Omega$-group allows us to consider continuity of operations and continuity of representation.

# Contents

CHAPTER 1

# Preface

## 1.1. Preface

It was first time that I drew my attention to the universal algebra which is Abelian Group relative to sum in the chapter [8]-2. Later I found that such universal algebra is called $\Omega$-group.

$\Omega$-group is interesting for me, since I think that we may consider calculus in $\Omega$-group. So there was a decision to study normed $\Omega$-group.

I have devoted much attention to the topology of normed $\Omega$-group.

## 1.2. Conventions

CONVENTION 1.2.1. *Let $A$ be $\Omega_1$-algebra. Let $B$ be $\Omega_2$-algebra. Notation*

$$A \xrightarrow{*} B$$

*means that there is representation of $\Omega_1$-algebra $A$ in $\Omega_2$-algebra $B$.* □

Without a doubt, the reader may have questions, comments, objections. I will appreciate any response.

CHAPTER 2

# Normed $\Omega$-Group

## 2.1. Normed $\Omega$-Group

Let sum which is not necessarily commutative be defined in $\Omega$-algebra[2.1] $A$.

DEFINITION 2.1.1. A map
$$f : A \to A$$
is called **additive map** if
$$f(a + b) = f(a) + f(b)$$
$\square$

DEFINITION 2.1.2. A map
$$f : A^n \to A$$
is called **polyadditive map** if for any $i$, $i = 1, ..., n$,
$$f(a_1, ..., a_i + b_i, ..., a_n) = f(a_1, ..., a_i, ..., a_n) + f(a_1, ..., b_i, ..., a_n)$$
$\square$

DEFINITION 2.1.3. If $\Omega$-algebra $A$ is group relative to sum and any operation $\omega \in \Omega$ is polyadditive map, then $\Omega$-algebra $A$ is called $\Omega$**-group**.[2.2] If $\Omega$-group $A$ is associative group relative to sum, then $\Omega$-algebra $A$ is called **associative $\Omega$-group**. If $\Omega$-group $A$ is Abelian group relative to sum, then $\Omega$-algebra $A$ is called **Abelian $\Omega$-group**. $\square$

EXAMPLE 2.1.4. The most evident example of $\Omega$-group is group itself. $\square$

EXAMPLE 2.1.5. A ring is $\Omega$-group. $\square$

EXAMPLE 2.1.6. Biring of matrices over division ring ([6]) is $\Omega$-group. $\square$

CONVENTION 2.1.7. *We assume that considered $\Omega$-group is Abelian.* $\square$

Element of $\Omega$-group $A$ is called $A$**-number**.

DEFINITION 2.1.8. Let $A$ be $\Omega$-group. If a subgroup $B$ of additive group $A$ is closed relative operations from $\Omega$, then the subgroup $B$ is called **subgroup of $\Omega$-group** $A$. $\square$

---

[2.1]See definition of universal algebra in [1, 11].
[2.2]You can find the definition of $\Omega$-group on the website

http://ncatlab.org/nlab/show/Omega-group

DEFINITION 2.1.9. **Norm on Ω-group** $A$[2.3] is a map

$$d \in A \to \|d\| \in R$$

which satisfies the following axioms

2.1.9.1: $\|a\| \geq 0$
2.1.9.2: $\|a\| = 0$ if, and only if, $a = 0$
2.1.9.3: $\|a + b\| \leq \|a\| + \|b\|$
2.1.9.4: $\|-a\| = \|a\|$

The Ω-group $A$, endowed with the structure defined by a given norm on $A$, is called **normed Ω-group**. □

REMARK 2.1.10. If Ω-group $B$ with norm $\|b\|_B$ is subgroup of Ω-group $A$ with norm $\|b\|_A$, then we require $\|b\|_B = \|b\|_A$ □

THEOREM 2.1.11. *Let $A$ be normed Ω-group. Then*

$$\|a - b\| \geq |\|a\| - \|b\|| \tag{2.1.1}$$

PROOF. Since $a = a - b + b$, then

$$\|a\| \leq \|a - b\| + \|b\| \tag{2.1.2}$$

follows from the statement 2.1.9.3. inequality

$$\|a - b\| \geq \|a\| - \|b\| \tag{2.1.3}$$

follows from the inequality (2.1.2). From the statement 2.1.9.4 and the inequality (2.1.3), it follows that

$$\|a - b\| = \|b - a\| \geq \|b\| - \|a\| \tag{2.1.4}$$

The inequality (2.1.1) follows from inequalities (2.1.3), (2.1.4). □

DEFINITION 2.1.12. Let $A$ be normed Ω-group. For $n$-ary operation $\omega$, the value

$$\|\omega\| = \sup \frac{\|a_1 ... a_n \omega\|}{\|a_1\| ... \|a_n\|} \tag{2.1.5}$$

is called **norm of operation** $\omega$. □

THEOREM 2.1.13. *Let $A$ be normed Ω-group. For $n$-ary operation $\omega$,*

$$\|a_1 ... a_n \omega\| \leq \|\omega\| \|a_1\| ... \|a_n\| \tag{2.1.6}$$

PROOF. From the equation (2.1.5), it follows that

$$\frac{\|a_1 ... a_n \omega\|}{\|a_1\| ... \|a_n\|} \leq \sup \frac{\|a_1 ... a_n \omega\|}{\|a_1\| ... \|a_n\|} = \|\omega\| \tag{2.1.7}$$

The inequality (2.1.6) follows from the inequality (2.1.7). □

DEFINITION 2.1.14. Let $A$ be normed Ω-group. Let $a \in A$. The set

$$B_o(a, R) = \{b \in A : \|b - a\| < R\}$$

is called **open ball** with center at $a$. □

---

[2.3]I made definition according to definition from [14], IX, §3.2 and definition [15]-1.1.12, p. 23.

DEFINITION 2.1.15. Let $A$ be normed $\Omega$-group. Let $a \in A$. The set
$$B_c(a, R) = \{b \in A : \|b - a\| \leq R\}$$
is called **closed ball** with center at $a$. □

THEOREM 2.1.16. *The statement $a \in B_o(b, R)$ follows from the statement[2.4] $b \in B_o(a, R)$.*

PROOF. The theorem follows from the definition 2.1.14. □

DEFINITION 2.1.17. Let $A$ be normed $\Omega$-group. Element $a \in A$ is called **limit of a sequence** $a_n$
$$a = \lim_{n \to \infty} a_n$$
if for every $\epsilon \in R$, $\epsilon > 0$, there exists positive integer $n_0$ depending on $\epsilon$ and such, that $\|a_n - a\| < \epsilon$ for every $n > n_0$. We also say that **sequence** $a_n$ **converges** to $a$. □

THEOREM 2.1.18. *Let $A$ be normed $\Omega$-group. Element $a \in A$ is limit of a sequence $a_n$*
$$a = \lim_{n \to \infty} a_n$$
*if for every $\epsilon \in R$, $\epsilon > 0$, there exists positive integer $n_0$ depending on $\epsilon$ and such, that*
$$a_n \in B_o(a, \epsilon)$$
*for every $n > n_0$.*

PROOF. The theorem follows from definitions 2.1.14, 2.1.17. □

DEFINITION 2.1.19. Let $A$ be normed $\Omega$-group. The sequence $a_n$, $a_n \in A$ is called **fundamental** or **Cauchy sequence**, if for every $\epsilon \in R$, $\epsilon > 0$, there exists positive integer $n_0$ depending on $\epsilon$ and such, that $\|a_p - a_q\| < \epsilon$ for every $p$, $q > n_0$. □

THEOREM 2.1.20. *Let $A$ be normed $\Omega$-group. The sequence $a_n$, $a_n \in A$ is fundamental sequence, if for every $\epsilon \in R$, $\epsilon > 0$, there exists positive integer $n_0$ depending on $\epsilon$ and such, that*
$$a_q \in B_o(a_p, \epsilon)$$
*for every $p$, $q > n_0$.*

PROOF. The theorem follows from definitions 2.1.14, 2.1.19. □

THEOREM 2.1.21. *Let $A$ be normed $\Omega$-group. Let $a_n$, $n = 1, ...,$ be fundamental sequence. Let $b_n$, $n = 1, ...,$ be sequence. Let*

(2.1.8) $$\lim_{n \to \infty} (a_n - b_n) = 0$$

*Then $b_n$ is fundamental sequence.*

---

[2.4]Similar theorem is true for closed balls.

PROOF. From the equation (2.1.8) and the definition 2.1.17, it follows that for given $\epsilon \in R$, $\epsilon > 0$, there exists positive integer $N_1$ depending on $\epsilon$ and such, that

$$\|a_n - b_n\| < \frac{\epsilon}{3} \qquad (2.1.9)$$

for every $n > N_1$. According to the definition 2.1.19, for given $\epsilon \in R$, $\epsilon > 0$, there exists positive integer $N_2$ depending on $\epsilon$ and such, that

$$\|a_p - a_q\| < \frac{\epsilon}{3} \qquad (2.1.10)$$

for every $p, q > N_2$. Let

$$N = \max(N_1, N_2)$$

From inequalities (2.1.9), (2.1.10) it follows that for given $\epsilon \in R$, $\epsilon > 0$, there exists positive integer $N$ depending on $\epsilon$ and such, that

$$\|b_p - b_q\| = \|b_p - a_p + a_p - a_q + a_q - b_q\| \le \|b_p - a_p\| + \|a_p - a_q\| + \|a_q - b_q\| < \epsilon$$

for every $p, q > N$. According to the definition 2.1.19, the sequence $b_n$ is fundamental sequence. □

THEOREM 2.1.22. *Let $A$ be normed Ω-group. Let $a_n$, $b_n$, $n = 1, ...$, be fundamental sequences. Let*

$$\lim_{n \to \infty} (a_n - b_n) = 0 \qquad (2.1.11)$$

*If the sequence $a_n$ converges, then the sequence $b_n$ converges and*

$$\lim_{n \to \infty} a_n = \lim_{n \to \infty} b_n \qquad (2.1.12)$$

PROOF. From the equation (2.1.11) and the definition 2.1.17, it follows that for given $\epsilon \in R$, $\epsilon > 0$, there exists positive integer $N_1$ depending on $\epsilon$ and such, that

$$\|a_n - b_n\| < \frac{\epsilon}{2} \qquad (2.1.13)$$

for every $n > N_1$. According to the definition 2.1.23 and the theorem 2.1.24, there exists limit $a$ of the sequence $a_n$. According to the definition 2.1.17, for given $\epsilon \in R$, $\epsilon > 0$, there exists positive integer $N_2$ depending on $\epsilon$ and such, that

$$\|a_n - a\| < \frac{\epsilon}{2} \qquad (2.1.14)$$

for every $n > N_2$. Let

$$N = \max(N_1, N_2)$$

From inequalities (2.1.13), (2.1.14) and from the statement 2.1.9.3, it follows that for given $\epsilon \in R$, $\epsilon > 0$, there exists positive integer $N$ depending on $\epsilon$ and such, that

$$\|a - b_n\| = \|a - a_n + a_n - b_n\| \le \|a - a_n\| + \|a_n - b_n\| < \epsilon$$

for every $n > N$. According to the definition 2.1.17, the sequence $b_n$ converges to $a$. □

DEFINITION 2.1.23. Normed Ω-group $A$ is called **complete** if any fundamental sequence of elements of Ω-group $A$ converges, i.e. has limit in Ω-group $A$. □

THEOREM 2.1.24. *Let $A$ be complete Ω-group. Any fundamental sequence has one and only one limit.*

PROOF. Let $a_n$, $n = 1, ...$, be fundamental sequence. According to the definition 2.1.19, for given $\epsilon \in R$, $\epsilon > 0$, there exists positive integer $n_1$ depending on $\epsilon$ and such, that

$$\|a_p - a_q\| < \frac{\epsilon}{3} \tag{2.1.15}$$

for every $p$, $q > n_1$. Let

$$a = \lim_{n \to \infty} a_n \tag{2.1.16}$$

and

$$b = \lim_{n \to \infty} a_n \tag{2.1.17}$$

From the equation (2.1.16) and the definition 2.1.17, it follows that for given $\epsilon \in R$, $\epsilon > 0$, there exists positive integer $n_2$ depending on $\epsilon$ and such, that

$$\|a_p - a\| < \frac{\epsilon}{3} \tag{2.1.18}$$

for every $p > n_2$. From the equation (2.1.17) and the definition 2.1.17, it follows that for given $\epsilon \in R$, $\epsilon > 0$, there exists positive integer $n_3$ depending on $\epsilon$ and such, that

$$\|a_q - b\| < \frac{\epsilon}{3} \tag{2.1.19}$$

for every $q > n_3$. Let

$$n_0 = \max(n_1, n_2, n_3)$$

From inequalities (2.1.15), (2.1.18), (2.1.19), it follows that for given $\epsilon \in R$, $\epsilon > 0$, there exists positive integer $n_0$ depending on $\epsilon$ and such, that

$$\|a - b\| = \|a - a_p + a_p - a_q + a_q - b\| \leq \|a - a_p\| + \|a_p - a_q\| + \|a_q - b\| < \epsilon$$

for every $p$, $q > n_0$. Therefore, $\|a - b\| = 0$. According to the statement 2.1.9.2, $a = b$.  □

## 2.2. Topology of Ω-Group

REMARK 2.2.1. Invariant distance on additive group of Ω-group $A$

$$d(a, b) = \|a - b\| \tag{2.2.1}$$

defines topology of metric space, compatible with Ω-group structure of $A$. The set of open balls of normed Ω-group $A$ is base of topology,[2.5] compatible with distance (2.2.1).  □

DEFINITION 2.2.2. Let $A$ be normed Ω-group. A set $U \subset A$ is called **open**,[2.6] if for any $A$-number $a \in U$ there exists $\epsilon \in R$, $\epsilon > 0$, such that $B_o(a, \epsilon) \subset U$.  □

---

[2.5] See definition of base of topology in [3], page 81, definition 3.

[2.6] In topology, we usualy define an open set before we define base of topology. In the case of a metric or normed space, it is more convenient to define an open set, based on the definition of base of topology. In such case, the definition is based on one of the properties of base of topology. An immediate proof allows us to see that defined such an open set satisfies the basic properties.

DEFINITION 2.2.3. Let $B$ be subset of topological space $A$. A point $x \in A$ is called a **contact point of set** $B$,[2.7] if for any $\epsilon \in R$, $\epsilon > 0$, open ball $B_o(x, \epsilon)$. contains at least one point of $B$

$$B_o(x, \epsilon) \cap B \neq \emptyset$$

The set of all contact points of the set $B$ is called **closure of set** $B$. Closure of set $B$ is denoted $[B]$. □

REMARK 2.2.4. A set $B$ is closed iff $B = [B]$.

This statement can be either the definition of a closed set, or theorem. For instance, consider the definition of a closed set as a complement to an open set. Since $x \in A \setminus [B]$. then there exists neighborhood of the point $x$, which does not contain points of the set $B$. Therefore, $A \setminus [B]$ is open set. □

DEFINITION 2.2.5. Let $A$ be normed Ω-group. A set $B \subset A$ is called **dense in set**[2.8] $C \subset A$, if $C \subset [B]$. A set $B \subset A$ is called **everywhere dense**,[2.9] if $[B] = A$. □

DEFINITION 2.2.6 (**second axiom of countability**). In topological space with countable base, there is at least one base containing no more than countably many sets.[2.10] □

DEFINITION 2.2.7 (First axiom of separation). A topological space $A$ is called $T_1$-**space**,[2.11] if for points $x, y \in A$, $x \neq y$, there is neighborhood $U_x$, $x \in U_x$, and neighborhood $U_y$, $y \in U_y$, such that $y \notin U_x$, $x \notin U_y$. □

THEOREM 2.2.8. *Any point of $T_1$-space is closed set.*[2.12]

PROOF. Since $x \neq y$, then according to the definition 2.2.7, there exists a neighborhood $O_y$ of the point $y$ such that $x \notin O_y$. According to the definition 2.2.3, $y \notin [\{x\}]$. Therefore, $\{x\} = [\{x\}]$. According to the remark 2.2.4, $\{x\}$ is closed set. □

THEOREM 2.2.9. *Normed Ω-group $A$ is $T_1$-space iff for points $x, y \in A$, $x \neq y$, there exists open balls $B_o(x, r_x)$, $B_o(y, r_y)$ such that $y \notin B_o(x, r_x)$, $x \notin B_o(y, r_y)$.*

PROOF. According to the definition 2.2.7, there is neighborhood $U_x$, $x \in U_x$, and neighborhood $U_y$, $y \in U_y$, such that $y \notin U_x$, $x \notin U_y$. According to the remark 2.2.1 and the definition [3]-3 on page 81, there exists open balls $B_o(x, r_x)$, $B_o(y, r_y)$ such that $B_o(x, r_x) \subset U_x$, $B_o(y, r_y) \subset U_y$. Therefore, $y \notin B_o(x, r_x)$, $x \notin B_o(y, r_y)$. □

CONVENTION 2.2.10. *In order for topology of normed Ω-group $A$ to be nontrivial, we require following.*

2.2.10.1: *Ω-group $A$ is $T_1$-space.*
2.2.10.2: *Normed Ω-group $A$ satisfies the second axiom of countability.*

---

[2.7]See definition of contact point and closure of set in [3], page 79.

[2.8]See also definition in [3], page 48.

[2.9]We also say that the set $B$ is everywhere dense in Ω-group $A$.

[2.10]See also the definition in [3], page 82. According to the theorem [3]-4 on page 82, in topological space satisfying the second axiom of countability, there exists a countable everywhere dense subset.

[2.11]See similar definition in [3], p. 85, definition 4.

[2.12]See also the theorem [3]-8 on page 85.

2.2.10.3: *The set $\{x\}$ is not an open set.*

□

THEOREM 2.2.11. *Let $A$ be normed $\Omega$-group. For any $x \in A$, $r \in R$, there exists $y \in B_o(x, r)$, $y \neq x$.*

PROOF. If we assume that there is no $y \in B_o(x, r)$, $y \neq x$, then in such case $B_o(x, r) = \{x\}$ and the set $\{x\}$ is open set. This statement contradicts to the convention 2.2.10.3. Therefore, assumption is not correct. □

THEOREM 2.2.12. *Let $A$ be normed $\Omega$-group. For any $a, b \in A$, $b \neq a$, there exists $c \in A$, $c \neq a$, such that*

$$\|c - a\| < \|b - a\| \tag{2.2.2}$$

PROOF. According to the theorem 2.2.9, there exists an open ball $B_o(a, r)$ such that $b \notin B_o(a, r)$. According to the definition 2.1.14

$$r \leq \|b - a\| \tag{2.2.3}$$

According to the theorem 2.2.11, there exists $c \in B_o(a, r)$, $c \neq a$. According to the definition 2.1.14

$$\|c - a\| < r \tag{2.2.4}$$

The inequality (2.2.2) follows from inequalities (2.2.3), (2.2.4). □

THEOREM 2.2.13. *Let $A$ be normed $\Omega$-group. For any $a \in A$, there exists a sequence of $A$-numbers $a_n$, $a_n \neq a$, such that*

$$\lim_{n \to \infty} a_n = a \tag{2.2.5}$$

PROOF. There are different ways to construct a sequence $a_n$.

According to the theorem 2.2.11, for any $n$ there exists $a_n \in B(a, \frac{1}{n})$, $a_n \neq a$. According to the theorem 2.1.18, the sequence $a_n$ converges to $a$.

It is possible that $a_n = a_{n+1}$ in considered sequence. However, it is evident that if $m > \frac{1}{\|a_n\|}$, then $a_m \neq a_n$. We can construct a sequence in which all the elements are different.

Let $a_1 \in A$. Suppose we have chosen $a_n \in A$, $n \geq 1$. Let

$$r = \frac{\|a_n - a\|}{2} \tag{2.2.6}$$

According to the theorem 2.2.11, there exists $a_{n+1} \in B(a, r)$, $a_{n+1} \neq a$. According to the equation (2.2.6) and the definition 2.1.14, $a_{n+1} \neq a_n$. According to the theorem 2.1.18, the sequence $a_n$ converges to $a$. □

THEOREM 2.2.14. *Let $A$ be normed $\Omega$-group. Let a set $B \subset A$ is dense in set $C \subset A$, . Then, for any $A$-number $b \in C$, there exists a sequence of $A$-numbers $b_n$, $b_n \in B$, convergent to $b$*

$$b = \lim_{n \to \infty} b_n \tag{2.2.7}$$

PROOF. According to definitions 2.2.3, 2.2.5, for any $n > 0$ there exists $b_n$ such that

$$b_n \in B \cap B_o(b, 1/n)$$

According to the theorem 2.1.18, the sequence $b_n$ converges to $b$, since $b_n \in B_o(x, \epsilon)$ for any $n > 1/\epsilon$. □

DEFINITION 2.2.15. A set $T$ of topological space is called **compact**, if every open cover of $T$ has finite subcover.[2.13] □

DEFINITION 2.2.16. A set $T$ of topological space is called **connected**,[2.14] if there exist open sets $B$ and $C$ such that conditions

2.2.16.1: $T \subset B \cup C$
2.2.16.2: $T \cap B \neq \emptyset$
2.2.16.3: $T \cap C \neq \emptyset$
imply $T \cap B \cap C \neq \emptyset$. □

EXAMPLE 2.2.17. Open connected set in real field is an open interval.[2.15] □

## 2.3. Continuous Map of Ω-Group

DEFINITION 2.3.1. A map
$$f : A_1 \to A_2$$
of normed $\Omega_1$-group $A_1$ with norm $\|x\|_1$ into normed $\Omega_2$-group $A_2$ with norm $\|y\|_2$ is called **continuous**, if for every as small as we please $\epsilon > 0$ there exist such $\delta > 0$, that
$$\|x' - x\|_1 < \delta$$
implies
$$\|f(x') - f(x)\|_2 < \epsilon$$
□

THEOREM 2.3.2. *A map*
$$f : A_1 \to A_2$$
*of normed $\Omega_1$-group $A_1$ with norm $\|x\|_1$ into normed $\Omega_2$-group $A_2$ with norm $\|y\|_2$ is continuous, iff for every as small as we please $\epsilon > 0$ there exist such $\delta > 0$, that*
$$b \in B_o(a, \delta)$$
*implies*
$$f(b) \in B_o(f(a), \epsilon)$$

PROOF. The theorem follows from definitions 2.1.14, 2.3.1. □

THEOREM 2.3.3. *A map*
$$f : A_1 \to A_2$$
*of normed $\Omega_1$-group $A_1$ with norm $\|x\|_1$ into normed $\Omega_2$-group $A_2$ with norm $\|y\|_2$ is continuous, iff preimage of an open set is the open set.*

PROOF. Let $U \subset f(A_1) \subset A_2$ be open set and
(2.3.1) $$W = f^{-1}(U)$$

---

[2.13]See also definition [3]-1, page 92.
[2.14]See also definitions 1 and 2 in [13], pages 107, 108.
[2.15]See comment in [3], page 55.

- Let the map $f$ be continuous. Let $a \in W$. Then $f(a) \in U$. According to the definition 2.2.2, there exists $\epsilon \in R$, $\epsilon > 0$, such that

(2.3.2) $$B_o(f(a), \epsilon) \subset U$$

  According to the theorem 2.3.2, there exist such $\delta > 0$, that

(2.3.3) $$b \in B_o(a, \delta)$$

  implies

(2.3.4) $$f(b) \in B_o(f(a), \epsilon)$$

  From equations (2.3.2), (2.3.4) it follows that

(2.3.5) $$f(b) \in U$$

  From equations (2.3.1), (2.3.5) it follows that

(2.3.6) $$b \in W$$

  Since from the statement (2.3.3) implies the statement (2.3.6), then $B_o(a, \delta) \subset W$. According to the definition 2.2.2, the set $W$ is open.

- Let preimage of an open set be the open set. Let $U = B_o(f(a), \epsilon)$. Since $W$ is open set, then according to the definition 2.2.2 there exists $\delta > 0$ such that $B_o(a, \delta) \subset W$. Therefore, the statement $b \in B_o(a, \delta)$ implies $f(b) \in B_o(f(a), \epsilon)$. According to the theorem 2.3.2, the map $f$ is continuous.

□

THEOREM 2.3.4. *Let*
$$f : X \to Y$$
*be continues map of topological space $X$ into topological space $Y$. Let $T$ be connected set in topological space $X$. Then $f(T)$ is connected set in topological space $Y$.*[2.16]

PROOF. Let $B$, $C$ be open sets in topological space $Y$. Since the map $f$ is continuous, then according to the theorem 2.3.3, sets $f^{-1}(B)$, $f^{-1}(C)$ are open sets. Since statements 2.2.16.1, 2.2.16.2, 2.2.16.3 are true for sets $B$, $C$, $f(T)$, the statements 2.2.16.1, 2.2.16.2, 2.2.16.3 are true for sets $f^{-1}(B)$, $f^{-1}(C)$, $T$. If we assume that the set $f(T)$ is not connected, then according to the definition 2.2.16

(2.3.7) $$B \cap C \cap f(T) = \emptyset$$

From the equation (2.3.7), it follows that

(2.3.8) $$f^{-1}(B) \cap f^{-1}(C) \cap T = \emptyset$$

According to the definition 2.2.16, the set $T$ is not connected. This statement contradicts the assumption of the theorem. Therefore, the set $f(T)$ is connected.

□

THEOREM 2.3.5. *Let*
$$f : R \to R$$
*be continuous map of real field. Then image of interval is interval.*

PROOF. The theorem follows from the theorem 2.3.4 and the example 2.2.17.

□

---

[2.16] See also proposition 4 in [13], page 109.

THEOREM 2.3.6. *A map*
$$f : A_1 \to A_2$$
*of normed $\Omega_1$-group $A_1$ with norm $\|x\|_1$ into normed $\Omega_2$-group $A_2$ with norm $\|y\|_2$ is continuous, iff the condition that a sequence of $A_1$-numbers $a_n$ converges implies that a sequence of $A_2$-numbers $f(a_n)$ converges and following equation is true*

(2.3.9) $$f(\lim_{n \to \infty} a_n) = \lim_{n \to \infty} f(a_n)$$

PROOF. Let $A_1$-number $a$ be limit of the sequence $a_n$

(2.3.10) $$a = \lim_{n \to \infty} a_n$$

- Let the equation (2.3.9) is true for $A$-number $a$. The equation

(2.3.11) $$f(a) = \lim_{n \to \infty} f(a_n)$$

follows from equations (2.3.9), (2.3.10). From the equation (2.3.11) and the theorem 2.1.18, it follows that for given $\epsilon \in R$, $\epsilon > 0$, there exists positive integer $n_2$ depending on $\epsilon$ and such, that

(2.3.12) $$f(a_n) \in B_o(f(a), \epsilon)$$

for every $n > n_2$. From the equation (2.3.10) and the theorem 2.1.18, it follows that for given $\delta \in R$, $\delta > 0$, there exists positive integer $n_1$ depending on $\delta$ and such, that

(2.3.13) $$a_n \in B_o(a, \delta)$$

for every $n > n_1$. If $n_1 \le n_2$, then estimate of $\delta$ is too large. We assume
$$\delta = \frac{\|a - a_{n_1+1}\|}{2}$$
and repeat estimate of value of $n_1$. It is evident that new value of $n_1$ will be greater than the previous. Therefore, after a finite number of iterations, we find $\delta$ such that $n_1 > n_2$ and the statement (2.3.13) implies the statement (2.3.12). According to the theorem 2.3.2, the map $f$ is continuous.

- Let the map $f$ be continuous. According to the theorem 2.3.2, for every as small as we please $\epsilon > 0$ there exist such $\delta > 0$, that

(2.3.14) $$b \in B_o(a, \delta)$$

implies

(2.3.15) $$f(b) \in B_o(f(a), \epsilon)$$

From the equation (2.3.10) and the theorem 2.1.18, it follows that for given $\delta \in R$, $\delta > 0$, there exists positive integer $n_1$ depending on $\delta$ and such, that

(2.3.16) $$a_n \in B_o(a, \delta)$$

for every $n > n_1$. From equations (2.3.14), (2.3.15), (2.3.16) it follows that for given $\epsilon \in R$, $\epsilon > 0$, there exists positive integer $n_1$ depending on $\epsilon$ and such, that

(2.3.17) $$f(a_n) \in B_o(f(a), \epsilon)$$

for every $n > n_1$. The equation

(2.3.18) $$f(a) = \lim_{n \to \infty} f(a_n)$$

follows from the equation (2.3.17) and the theorem 2.1.18. The equation (2.3.9) follows from equations (2.3.18), (2.3.10).

□

DEFINITION 2.3.7. Let
$$f : A_1 \to A_2$$
be map of normed $\Omega_1$-group $A_1$ with norm $\|x\|_1$ into normed $\Omega_2$-group $A_2$ with norm $\|y\|_2$. Value

(2.3.19) $$\|f\| = \sup \frac{\|f(x)\|_2}{\|x\|_1}$$

is called **norm of map** $f$. □

THEOREM 2.3.8. *Let*
$$f : A_1 \to A_2$$
*be additive map of normed $\Omega_1$-group $A_1$ with norm $\|x\|_1$ into normed $\Omega_2$-group $A_2$ with norm $\|y\|_2$. Since $\|f\| < \infty$, then map $f$ is continuous.*

PROOF. Since map $f$ is additive, then according to definition 2.3.7
$$\|f(x) - f(y)\|_2 = \|f(x - y)\|_2 \le \|f\| \, \|x - y\|_1$$
Let us assume arbitrary $\epsilon > 0$. Assume $\delta = \dfrac{\epsilon}{\|f\|}$. Then
$$\|f(x) - f(y)\|_2 \le \|f\| \, \delta = \epsilon$$
follows from inequality
$$\|x - y\|_1 < \delta$$
According to definition 2.3.1 map $f$ is continuous. □

THEOREM 2.3.9. *Let $A$ be normed $\Omega$-group. Norm defined in $\Omega$-group $A$ is continues map of $\Omega$-group $A$ into real field.*

PROOF. Let $\epsilon \in R$, $\epsilon > 0$, and $a \in A$. According to the theorem 2.2.11, there exists $b \in B_o(a, \epsilon)$, $b \ne a$. According to the definition 2.1.14

(2.3.20) $$\|a - b\| < \epsilon$$

From inequalities (2.3.20), (2.1.1), it follows that

(2.3.21) $$|\, \|a\| - \|b\| \,| \le \|a - b\| < \epsilon$$

Since the inequality (2.3.21) follows from the inequality (2.3.20), then the norm is continuous according to the definition 2.3.1. □

THEOREM 2.3.10. *Let*
$$f : A \to B$$
*be continues map of normed $\Omega$-algebra $A$ into normed $\Omega$-algebra $B$. If $C$ is compact set of $\Omega$-algebra $A$, then the image $f(C)$ is compact set of $\Omega$-algebra $B$.*

PROOF. Let $D$ be open cover of the set $f(C)$. Any set $E \in D$ is open set. According to the theorem 2.3.3, preimage $f^{-1}(E)$ is open set. Let

$$D' = \{E' : E' = f^{-1}(E), E \in D\}$$

The statement

$$C = \bigcup_{E' \in D'} E'$$

follows from the statement

$$f(C) = \bigcup_{E \in D} E$$

Therefore, $D'$ is open cover of the set $C$. According to the definition 2.2.15, there exists finite subcover $E'_1, ..., E'_n$. Therefore, the cover $D$ has finite subcover $E_1 = f(E'_1), ..., E_n = f(E'_n)$. According to the definition 2.2.15, the set $f(C)$ is compact set. □

THEOREM 2.3.11. *Let $C$ be compact set of normed $\Omega$-group $A$. Then the norm $\|x\|$, $x \in C$, is bounded from both sides.*

PROOF. According to the theorem 2.3.9, the norm is continues map of $\Omega$-group $A$ into real field. According to the theorem 2.3.10, the set

$$C' = \{\|x\| : x \in C\}$$

is compact set in real field. According to the statement [2]-3.96.e, page 118, the set $C'$ is bounded and contains its exact boundaries. □

## 2.4. Continuity of Operations in Ω-Group

From the statement 2.1.9.3 and the theorem 2.3.8, it follows that sum in normed $\Omega$-group is continuous. From the definition 2.1.12, it follows that operations in normed $\Omega$-group are continuous.[2.17] In this section we consider statements related with continuity of operations.

THEOREM 2.4.1. *In normed $\Omega$-group $A$, the following inequality is true*

(2.4.1) $$\|(c_1 + c_2) - (a_1 + a_2)\| \le \|(c_1 - a_1)\| + \|(c_2 - a_2)\|$$

PROOF. From the statement 2.1.9.3, it follows that

(2.4.2) $$\begin{aligned}\|(c_1 + c_2) - (a_1 + a_2)\| &= \|(c_1 - a_1) + (c_2 - a_2)\| \\ &\le \|(c_1 - a_1)\| + \|(c_2 - a_2)\|\end{aligned}$$

The inequality (2.4.1) follows from the inequality (2.4.2). □

THEOREM 2.4.2. *Let $A$ be normed $\Omega$-group. For $c_1$, $c_2 \in A$, let $c_1 \in B_c(a_1, R_1)$, $c_2 \in B_c(a_2, R_2)$. Then*

(2.4.3) $$c_1 + c_2 \in B_c(a_1 + a_2, R_1 + R_2)$$

---

[2.17]Compare the definition 2.1.12 of norm of operation and the definition [9]-3.4.8 of norm of polylinear map.

PROOF. According to the definition 2.1.14

$$\|c_1 - a_1\| \leq R_1$$
$$\|c_2 - a_2\| \leq R_2 \tag{2.4.4}$$

From the inequalities (2.4.1), (2.4.4), it follows that

$$\|(c_1 + c_2) - (a_1 + a_2)\| \leq \|(c_1 - a_1)\| + \|(c_2 - a_2)\| \leq R_1 + R_2 \tag{2.4.5}$$

The statement (2.4.3) follows from the inequality (2.4.5) and the definition 2.1.14. □

THEOREM 2.4.3. *Let $A$ be normed $\Omega$-group. Let the sequence of $A$-numbers $a_n$ converge and*

$$\lim_{n \to \infty} a_n = a \tag{2.4.6}$$

*Let the sequence of $A$-numbers $b_n$ converge and*

$$\lim_{n \to \infty} b_n = b \tag{2.4.7}$$

*Then the sequence of $A$-numbers $a_n + b_n$ converges and*

$$\lim_{n \to \infty} a_n + b_n = a + b \tag{2.4.8}$$

PROOF. From the equation (2.4.6) and the theorem 2.1.18, it follows that for given $\epsilon \in R$, $\epsilon > 0$, there exists $N_a$ such that the condition $n > N_a$ implies that

$$a_n \in B_o(a, \epsilon/2) \tag{2.4.9}$$

From the equation (2.4.7) and the theorem 2.1.18, it follows that for given $\epsilon \in R$, $\epsilon > 0$, there exists $N_b$ such that the condition $n > N_b$ implies that

$$b_n \in B_o(b, \epsilon/2) \tag{2.4.10}$$

Let

$$N = \max(N_a, N_b)$$

From equations (2.4.9), (2.4.10), the theorem 2.4.2 and condition $n > N$, it follows that

$$a_n + b_n \in B_o(a + b, \epsilon/2 + \epsilon/2) = B_o(a + b, \epsilon) \tag{2.4.11}$$

The equation (2.4.8) follows from the equation (2.4.11) and the theorem 2.1.18. □

THEOREM 2.4.4. *Let $A$ be normed $\Omega$-group. Let $\omega \in \Omega$ be $n$-ary operation. The following inequality is true*

$$\|c_1...c_n\omega - a_1...a_n\omega\| \leq \|\omega\|(\|c_1 - a_1\|C_2...C_n + ... + C_1...C_{n-1}\|c_n - a_n\|) \tag{2.4.12}$$

*where*

$$C_i = \max(\|a_i\|, \|c_i\|) \quad i = 1, ..., n \tag{2.4.13}$$

PROOF. According to the definitions 2.1.2, 2.1.3,

$$c_1...c_n\omega - a_1...a_n\omega = c_1c_2...c_n\omega - a_1c_2...c_n\omega$$
$$+ a_1c_2...c_n\omega - a_1a_2...c_n\omega$$
$$......$$
$$+ a_1...a_{n-1}c_n\omega - a_1...a_{n-1}a_n\omega \tag{2.4.14}$$
$$= (c_1 - a_1)c_2...c_n\omega + a_1(c_2 - a_2)...c_n\omega$$
$$+ ... + a_1...a_{n-1}(c_n - a_n)\omega$$

From the equation (2.4.14) and the statement 2.1.9.3, it follows that

$$\begin{aligned}\|c_1...c_n\omega - a_1...a_n\omega\| &\leq \|(c_1-a_1)c_2...c_n\omega\| + \|a_1(c_2-a_2)...c_n\omega\| \\ &\quad + ... + \|a_1...a_{n-1}(c_n-a_n)\omega\|\end{aligned}$$
(2.4.15)

From the equation (2.4.13) and the definition 2.1.12, it follows that

$$\|(c_1-a_1)c_2...c_n\omega\| \leq \|c_1-a_1\|C_2...C_n$$
(2.4.16)
$$......$$
$$\|a_1...a_{n-1}(c_n-a_n)\omega\| \leq C_1...C_{n-1}\|c_n-a_n\|$$

The inequality (2.4.12) follows from the inequalities (2.4.15), (2.4.16). □

THEOREM 2.4.5. *Let $A$ be normed $\Omega$-group. Let $\omega \in \Omega$ be n-ary operation. For $c_1, ..., c_n \in A$, let $c_1 \in B_c(a_1, R_1), ..., c_n \in B_c(a_n, R_n)$. Then*

$$c_1...c_n\omega \in B_c(a_1...a_n\omega, R)$$
(2.4.17)

*where*

$$R = \|\omega\|(R_1 C_2...C_n + ... + C_1...C_{n-1}R_n)$$
(2.4.18)

$$C_i = \|a_i\| + R_i \quad i = 1, ..., n$$
(2.4.19)

PROOF. According to the definition 2.1.14

$$\|c_1 - a_1\| \leq R_1$$
(2.4.20) $$...$$
$$\|c_n - a_n\| \leq R_n$$

The inequality

$$\|c_i\| = \|a_i + c_i - a_i\| \leq \|a_i\| + \|c_i - a_i\| \leq \|a_i\| + R_i$$
(2.4.21)

follows from the inequality (2.4.20) and from the statement 2.1.9.3. The inequality

$$C_i \leq \max(\|a_i\|, \|a_i\| + R_i) = \|a_i\| + R_i \quad i = 1, ..., n$$
(2.4.22)

follows from the inequality (2.4.21) and the equation (2.4.13). The equation (2.4.19) follows from the inequality (2.4.22). The statement (2.4.17) follows from inequalities (2.4.12), (2.4.20) and the equation (2.4.19). □

THEOREM 2.4.6. *Let $A$ be normed $\Omega$-group. Let $\omega \in \Omega$ be n-ari operation. Let the sequence of A-numbers $a_{i\cdot m}$ converge and*

$$\lim_{m \to \infty} a_{i\cdot m} = a_i$$
(2.4.23)

*Then the sequence of A-numbers $a_{1\cdot m}...a_{n\cdot m}\omega$ converges and*

$$\lim_{m \to \infty} a_{1\cdot m}...a_{n\cdot m}\omega = a_1...a_n\omega$$
(2.4.24)

PROOF. From the equation (2.4.23) and the theorem 2.1.18, it follows that for given

$$\delta_1 \in R, \ \delta_1 > 0$$
(2.4.25)

there exists $M_i$ such that the condition $m > M_i$ implies that

$$a_{i\cdot m} \in B_o(a_i, \delta_1)$$
(2.4.26)

Let

$$M = \max(M_1, ..., M_n)$$

From equations (2.4.26), the theorem 2.4.5 and condition $m > M$, it follows that

(2.4.27) $$a_{1 \cdot m}...a_{n \cdot m}\omega \in B_o(a_1...a_n\omega, \epsilon_1)$$

(2.4.28) $$\epsilon_1 = \delta_1 C_2...C_n + ... + \delta_1 C_1...C_{n-1}$$

where

(2.4.29) $$C_i = \|a_i\| + \delta_1$$

From the equation (2.4.29) and statements 2.1.9.1, (2.4.25) it follows that

(2.4.30) $$C_i > 0 \quad \frac{dC_i}{d\delta_1} > 0$$

From equations (2.4.28), (2.4.29) and the statement (2.4.30), it follows that $\epsilon_1$ is polynomial strictly monotone increasing function of $\delta_1$ such that

$$\delta_1 = 0 \Rightarrow \epsilon_1 = 0$$

According to the theorem 2.3.5, the map (2.4.28) maps the interval $[0, \delta_1)$ into the interval $[0, \epsilon_1)$. According to the theorem 2.3.3, for given $\epsilon > 0$ there exist $\delta > 0$ such that

$$\epsilon_1(\delta) < \epsilon$$

According to construction, a value of $M$ depends on a value of $\delta_1$. We choose the value of $M$ corresponding to $\delta_1 = \delta$. Therefore, for given $\epsilon \in R, \epsilon > 0$, there exists $M$ such that the condition $m > M$ implies that

(2.4.31) $$a_{1 \cdot m}...a_{n \cdot m}\omega \in B_o(a_1...a_n\omega, \epsilon)$$

The equation (2.4.24) follows from the equation (2.4.31) and the theorem 2.1.18. $\square$

## 2.5. Completion of Normed $\Omega$-Group

DEFINITION 2.5.1. Let $A$ be normed $\Omega$-group. Complete $\Omega$-group $B$ is called **completion of normed $\Omega$-group $A$**,[2.18] if

2.5.1.1: $\Omega$-group $A$ is subgroup of $\Omega$-group $B$.[2.19]
2.5.1.2: $\Omega$-group $A$ is everywhere dense in $\Omega$-group $B$.

$\square$

THEOREM 2.5.2. Let $A$ be normed $\Omega$-group. Let $B_1$, $B_2$ be completion of $\Omega$-group $A$. There exists isomorphism of $\Omega$-group[2.20]

(2.5.1) $$f : B_1 \to B_2$$

such that

(2.5.2) $$f(a) = a \quad a \in A$$

(2.5.3) $$\|f(b)\| = \|b\|$$

---

[2.18]An existence of completion of topological space is very important, because it allows us to use continuity as tool to study topological space. Since a topological space has additional structure, then we expect that the completion has the same structure. See the definition of completion of metric space in the definition [3]-4, page 62. See the definition of completion of valued field in the proposition [4]-6 on page 270 and in the following definition on page 271. See the definition of completion of normed vector space in the theorem [5]-1.11.1 on page 55.

[2.19]More precisely, we use isomorphism to identify $\Omega$-group $A$ with subgroup of $\Omega$-group $B$.

[2.20]We also say that completion of $\Omega$-group $A$ is unique up to considered isomorphism.

PROOF. Let $b_1 \in B_1$. According to the statement 2.5.1.2 and the theorem 2.2.14, there exists a sequence of $A$-numbers $a_n$ such that

$$b_1 = \lim_{n \to \infty} a_n \tag{2.5.4}$$

From the equation (2.5.4), it follows that the sequence of $A$-numbers $a_n$ is fundamental sequence in Ω-group $A$. Therefore, the sequence $a_n$ converges to $b_2$ in complete Ω-group $B_2$. According to the theorem 2.1.22, $B_2$-number $b_2$ does not depend from choice of sequence of $a_n$, convergent to $b_1$. Therefore, $f(b_1) = b_2$.

The equation (2.5.2), follows from the theorem 2.2.13.

The map $f$ is bijective, since the construction is reversible and we can find $B_1$-number $b_1$ corresponding to $B_2$-number $b_2$.

According to the theorem 2.3.9, in Ω-group $B_1$ the following equation is true

$$\|b_1\| = \lim_{n \to \infty} \|a_n\| \tag{2.5.5}$$

and in Ω-group $B_2$ the following equation is true

$$\|b_2\| = \lim_{n \to \infty} \|a_n\| \tag{2.5.6}$$

The equation (2.5.3), follows from the equations (2.5.5), (2.5.6).

From theorems 2.4.3, 2.4.6, it follows that the map $f$ is isomorphism of Ω-group. □

The proof of existence of completion of normed Ω-group $A$ (the theorem 2.5.17) based on the proof of the theorem 2.5.2. Since $B$ is complete Ω-group which is the completion of normed Ω-group $A$, then each $B$-number is a limit of a fundamental sequence of $A$-numbers. So, to construct the Ω-group $B$, we consider the set $B'$ of fundamental sequences of $A$-numbers.

LEMMA 2.5.3. *The relation $\sim$ on the set $B'$*

$$a_n \sim b_n \Leftrightarrow \lim_{n \to \infty} (a_n - b_n) = 0 \tag{2.5.7}$$

*is equivalence relation.*[2.21]

PROOF.

2.5.3.1: Since $a_n = a_n$, then the relation $\sim$ is reflexive relation

$$a_n \sim a_n$$

2.5.3.2: From the statement 2.1.9.4 it follows that

$$\|a_n - b_n\| = \|b_n - a_n\|$$

Therefore, the relation $\sim$ is symmetric relation

$$a_n \sim b_n \Leftrightarrow b_n \sim a_n$$

2.5.3.3: Let $a_n$, $n = 1, ...$, be fundamental sequence of $A$-numbers. Let $b_n$, $c_n$, $n = 1, ...$, be sequences of $A$-numbers. The equation

$$\lim_{n \to \infty} (a_n - b_n) = 0 \tag{2.5.8}$$

---

[2.21] See the definition of equivalence relation in [12], page 5.

follows from the statement $a_n \sim b_n$. The equation

(2.5.9) $$\lim_{n \to \infty} (a_n - c_n) = 0$$

follows from the statement $a_n \sim c_n$.

From the theorem 2.1.21 and the equation (2.5.8), it follows that the sequence $b_n$ is fundamental sequence. According to the definition 2.1.17, from the equation (2.5.8), it follows that for every $\epsilon \in R, \epsilon > 0$, there exists positive integer $N_1$ depending on $\epsilon$ and such, that

(2.5.10) $$\|a_n - b_n\| < \frac{\epsilon}{2}$$

for every $n > N_1$.

From the theorem 2.1.21 and the equation (2.5.9) it follows that the sequence $c_n$ is fundamental sequence. According to the definition 2.1.17, from the equation (2.5.9), it follows that for every $\epsilon \in R, \epsilon > 0$, there exists positive integer $N_2$ depending on $\epsilon$ and such, that

(2.5.11) $$\|a_n - c_n\| < \frac{\epsilon}{2}$$

for every $n > N_2$.

Let

$$N = \max(N_1, N_2)$$

From inequalities (2.5.10), (2.5.11) and from the statement 2.1.9.3, it follows that for given $\epsilon \in R, \epsilon > 0$, there exists positive integer $N$ depending on $\epsilon$ and such, that

$$\|b_n - c_n\| = \|b_n - a_n + a_n - c_n\| \le \|b_n - a_n\| + \|a_n - c_n\| < \epsilon$$

for every $n > N$. The equation

(2.5.12) $$\lim_{n \to \infty} (b_n - c_n) = 0$$

follows from the definition 2.1.17. The statement $b_n \sim c_n$ follows from the equation (2.5.12).

2.5.3.4: From arguments 2.5.3.3 it follows that the relation $\sim$ is transitive relation

$$a_n \sim b_n, a_n \sim c_n \Rightarrow b_n \sim c_n$$

The lemma follows from statements 2.5.3.1, 2.5.3.1, 2.5.3.4. □

LEMMA 2.5.4. $a_n \sim b_n$ iff for every $\epsilon \in R, \epsilon > 0$, there exists positive integer $n_0$ depending on $\epsilon$ and such, that

(2.5.13) $$b_n \in B_o(a_n, \epsilon)$$

for every $n > n_0$.

PROOF. According to the definition 2.1.14, the statement (2.5.13) is true iff

$$\|a_n - b_n\| < \epsilon$$

The theorem is true according to the definition 2.1.17 □

Let $B$ be the set of classes of equivalent fundamental sequences of $A$-numbers. We use notation

(2.5.14) $$[a_n] = \{b_n : a_n \sim b_n\}$$

for class of sequences which are equivalent to the sequence $a_n$.

LEMMA 2.5.5. *The set $B$ is Abelian group relative to operation*

(2.5.15) $$[a_n] + [b_n] = [a_n + b_n]$$

PROOF. Let $A$ be normed Ω-group. Let sequences of $A$-numbers $a_n$, $b_n$ be fundamental sequences.

LEMMA 2.5.6. *The sequence of $A$-numbers $a_n + b_n$ is fundamental sequence.*

PROOF. Since the sequence of $A$-numbers $a_n$ is fundamental sequence, then, according to the definition 2.1.19, for every $\epsilon \in R$, $\epsilon > 0$, there exists positive integer $N_a$ depending on $\epsilon$ and such, that

(2.5.16) $$\|a_p - a_q\| < \frac{\epsilon}{2}$$

for every $p$, $q > N_a$. Since the sequence of $A$-numbers $b_n$ is fundamental sequence, then, according to the definition 2.1.19, for every $\epsilon \in R$, $\epsilon > 0$, there exists positive integer $N_b$ depending on $\epsilon$ and such, that

(2.5.17) $$\|b_p - b_q\| < \frac{\epsilon}{2}$$

for every $p$, $q > N_b$. Let

$$N = \max(N_a, N_b)$$

From inequalities (2.4.1), (2.5.16), (2.5.17) and condition $n > N$, it follows that

$$\|(a_p + b_p) - (a_q + b_q)\| \le \|(a_p - a_q)\| + \|(b_p - b_q)\| < \epsilon$$

Therefore, acording to the definition 2.1.19, the sequence of $A$-numbers $a_n + b_n$ is fundamental sequence.

LEMMA 2.5.7. *Let $a_n \sim c_n$, $b_n \sim d_n$. Then*

(2.5.18) $$a_n + b_n \sim c_n + d_n$$

PROOF. The equation

(2.5.19) $$\lim_{n \to \infty}(a_n - c_n) = 0$$

follows from the statement $a_n \sim c_n$. The equation

(2.5.20) $$\lim_{n \to \infty}(b_n - d_n) = 0$$

follows from the statement $b_n \sim d_n$.

From the theorem 2.1.21 and the equation (2.5.19), it follows that the sequence $c_n$ is fundamental sequence. According to the definition 2.1.17, from the equation (2.5.19), it follows that for every $\epsilon \in R$, $\epsilon > 0$, there exists positive integer $N_1$ depending on $\epsilon$ and such, that

(2.5.21) $$\|a_n - c_n\| < \frac{\epsilon}{2}$$

for every $n > N_1$.

From the theorem 2.1.21 and the equation (2.5.20), it follows that the sequence $d_n$ is fundamental sequence. According to the definition 2.1.17, from the equation (2.5.20), it follows that for every $\epsilon \in R$, $\epsilon > 0$, there exists positive integer $N_2$ depending on $\epsilon$ and such, that

(2.5.22) $$\|b_n - d_n\| < \frac{\epsilon}{2}$$

for every $n > N_2$.

Let
$$N = \max(N_1, N_2)$$

From inequalities (2.4.1), (2.5.21), (2.5.22), it follows that for given $\epsilon \in R$, $\epsilon > 0$, there exists positive integer $N$ depending on $\epsilon$ and such, that

(2.5.23) $$\|(a_n + b_n) - (c_n + d_n)\| \leq \|a_n - c_n\| + \|b_n - d_n\| < \epsilon$$

for every $n > N$. The equation

(2.5.24) $$\lim_{n \to \infty}((a_n + b_n) - (c_n + d_n)) = 0$$

follows from the definition 2.1.17. Therefore, the statement (2.5.18) follows from the equation (2.5.24).

From the lemma 2.5.6, it follows that the sequence $a_n + b_n$ is fundamental sequence. From the lemma 2.5.7, it follows that the operation (2.5.15) is defined properly.

2.5.7.1: From the equation
$$a_n + (b_n + c_n) = (a_n + b_n) + c_n$$
it follows that the operation (2.5.15) is associative
$$[a_n] + ([b_n] + [c_n]) = ([a_n] + [b_n]) + [c_n]$$

2.5.7.2: From the equation
$$a_n + b_n = b_n + a_n$$
it follows that the operation (2.5.15) is commutative
$$[a_n] + [b_n] = [b_n] + [a_n]$$

2.5.7.3: If we assume $b_n = 0$, then the equation
$$a_n + 0 = a_n$$
implies
$$[a_n] + [0] = [a_n]$$

2.5.7.4: The equation
$$a_n + (-a_n) = 0$$
implies
$$[a_n] + [-a_n] = [0]$$
Therefore, we can assume
$$-[a_n] = [-a_n]$$

The lemma follows from statements 2.5.7.1, 2.5.7.2, 2.5.7.3, 2.5.7.4. □

LEMMA 2.5.8. *We define operation $\omega \in \Omega$ on the set $B$ using equation*

$$(2.5.25) \qquad [a_{1 \cdot m}]...[a_{n \cdot m}]\omega = [a_{1 \cdot m}...a_{n \cdot m}\omega]$$

PROOF. Let $A$ be normed $\Omega$-group.

LEMMA 2.5.9. *Let sequences of $A$-numbers $a_{1 \cdot m}$, ..., $a_{n \cdot m}$, $m = 1$, ..., be fundamental sequences. The sequence of $A$-numbers $a_{1 \cdot m}...a_{n \cdot m}\omega$ is fundamental sequence.*

PROOF. Since the sequence $a_{i \cdot m}$ is a fundamental sequence of $A$-numbers, then according to the theorem 2.1.20, it follows that for given

$$(2.5.26) \qquad \delta_1 \in R, \ \delta_1 > 0$$

there exists positive integer $M_i$ depending on $\delta_1$ and such, that

$$(2.5.27) \qquad a_{i \cdot q} \in B_o(a_{i \cdot p}, \delta_1)$$

for every $p$, $q > M_i$. Let

$$M = \max(M_1, ..., M_n)$$

From statements (2.5.27) and the theorem 2.4.5, it follows that

$$a_{1 \cdot q}...a_{n \cdot q}\omega \in B_o(a_{1 \cdot p}...a_{n \cdot p}\omega, \epsilon_1)$$

$$(2.5.28) \qquad \epsilon_1 = \delta_1 C_2...C_n + ... + \delta_1 C_1...C_{n-1}$$

where

$$(2.5.29) \qquad C_i = \|a_{i \cdot p}\|_i + \delta_1$$

From the equation (2.5.29) and statements 2.1.9.1, (2.5.26) it follows that

$$(2.5.30) \qquad C_i > 0 \quad \frac{dC_i}{d\delta_1} > 0$$

From equations (2.5.28), (2.5.29) and the statement (2.5.30), it follows that $\epsilon_1$ is polynomial strictly monotone increasing function of $\delta_1$ such that

$$\delta_1 = 0 \Rightarrow \epsilon_1 = 0$$

According to the theorem 2.3.5, the map (2.5.28) maps the interval $[0, \delta_1)$ into the interval $[0, \epsilon_1)$. According to the theorem 2.3.3, for given $\epsilon > 0$ there exist $\delta > 0$ such that

$$\epsilon_1(\delta) < \epsilon$$

According to construction, a value of $M$ depends on a value of $\delta_1$. We choose the value of $M$ corresponding to $\delta_1 = \delta$. Therefore, for given $\epsilon \in R$, $\epsilon > 0$, there exists $M$ such that the condition $p$, $q > M$ implies that

$$(2.5.31) \qquad a_{1 \cdot q}...a_{n \cdot q}\omega \in B_o(a_{1 \cdot p}...a_{n \cdot p}\omega, \epsilon)$$

From the statement (2.5.31) and the theorem 2.1.20, it follows that the sequence of $A$-numbers $a_{1 \cdot m}...a_{n \cdot m}\omega$ is a fundamental sequence.

LEMMA 2.5.10. Let sequences of $A$-numbers $a_{1 \cdot m}, ..., a_{n \cdot m}$, $m = 1, ...,$ $c_{1 \cdot m}, ..., c_{n \cdot m}$, $m = 1, ...,$ be fundamental sequences. Let

(2.5.32)
$$a_{i \cdot m} \sim c_{i \cdot m}, \quad i = 1, ..., n$$

Then

(2.5.33)
$$a_{1 \cdot m}...a_{n \cdot m}\omega \sim c_{1 \cdot m}...c_{n \cdot m}\omega$$

PROOF. From the statement (2.5.32) and the lemma 2.5.4, it follows that for given

(2.5.34)
$$\delta_1 \in R, \ \delta_1 > 0$$

there exists positive integer $M_1$ depending on $\delta_1$ and such, that

(2.5.35)
$$c_{i \cdot m} \in B_o(a_{i \cdot m}, \delta_1)$$

for every $m > M_1$. Let

$$M = \max(M_1, ..., M_n)$$

From statements (2.5.35) and the theorem 2.4.5, it follows that

$$c_{1 \cdot m}...c_{n \cdot m}\omega \in B_o(a_{1 \cdot m}...a_{n \cdot m}\omega, \epsilon_1)$$

(2.5.36)
$$\epsilon_1 = \delta_1 C_2...C_n + ... + \delta_1 C_1...C_{n-1}$$

where

(2.5.37)
$$C_i = \|a_{i \cdot n}\|_i + \delta_1$$

From the equation (2.5.37) and statements 2.1.9.1, (3.3.23) it follows that

(2.5.38)
$$C_i > 0 \quad \frac{dC_i}{d\delta_1} > 0$$

From equations (2.5.36), (2.5.37) and the statement (2.5.38), it follows that $\epsilon_1$ is strictly monotone increasing function of $\delta_1$ such that

$$\delta_1 = 0 \Rightarrow \epsilon_1 = 0$$

According to the theorem 2.3.5, the map (2.5.36) maps the interval $[0, \delta_1)$ into the interval $[0, \epsilon_1)$. According to the theorem 2.3.3, for given $\epsilon > 0$ there exist $\delta > 0$ such that

$$\epsilon_1(\delta) < \epsilon$$

According to construction, a value of $M$ depends on a value of $\delta_1$. We choose the value of $M$ corresponding to $\delta_1 = \delta$. Therefore, for given $\epsilon \in R, \epsilon > 0$, there exists $M$ such that the condition $m > M$ implies that

(2.5.39)
$$c_{1 \cdot m}...c_{n \cdot m}\omega \in B_o(a_{1 \cdot m}...a_{n \cdot m}\omega, \epsilon)$$

The statement (2.5.33) follows from the statement (2.5.39) and the lemma 2.5.4.

From the lemma 2.5.9, it follows that the sequence $a_{1 \cdot m}...a_{n \cdot m}\omega$ is fundamental sequence. From the lemma 2.5.10, it follows that the operation (2.5.25) is defined properly. □

LEMMA 2.5.11. Abelian group $B$ is Abelian Ω-group.

PROOF. According to the definition 2.1.3, the operation $\omega \in \Omega$ is polyadditive. From the equation

$$a_{1 \cdot m}...(a_{i \cdot m} + b_{i \cdot m})...a_{n \cdot m}\omega = a_{1 \cdot m}...a_{i \cdot m}...a_{n \cdot m}\omega + a_{1 \cdot m}...b_{i \cdot m}...a_{n \cdot m}\omega$$

it follows that the operation (2.5.25) is polyadditive

$$[a_{1 \cdot m}]...[a_{i \cdot m} + b_{i \cdot m}]...[a_{n \cdot m}]\omega = [a_{1 \cdot m}]...[a_{i \cdot m}]...[a_{n \cdot m}]\omega + [a_{1 \cdot m}]...[b_{i \cdot m}]...[a_{n \cdot m}]\omega$$

The lemma follows from the definition 2.1.3. □

LEMMA 2.5.12. *Let the sequence $a_n$ be fundamental sequence. Then the sequence $\|a_n\|$ is fundamental sequence.*

PROOF. According to the definition 2.1.19, for given $\epsilon \in R, \epsilon > 0,$ there exists positive integer $N$ depending on $\epsilon$ and such, that

(2.5.40) $$\|a_p - a_q\| < \epsilon$$

for every $p, q > N$. From inequalities (2.1.1), (2.5.40) it follows that for every $\epsilon \in R, \epsilon > 0,$ there exists positive integer $N$ depending on $\epsilon$ and such, that

(2.5.41) $$|\, \|a_p\| - \|a_q\|\, | \leq \|a_p - a_q\| < \epsilon$$

for every $n > N$. According to the definition 2.1.19, the sequence $\|a_n\|$ is fundamental sequence. □

LEMMA 2.5.13. *Let $a_n \sim b_n$. Then*

(2.5.42) $$\lim_{n \to \infty} \|a_n\| = \lim_{n \to \infty} \|b_n\|$$

PROOF. The equation

(2.5.43) $$\lim_{n \to \infty} (a_n - b_n) = 0$$

follows from the statement $a_n \sim b_n$. According to the definition 2.1.17, from the equation (2.5.43), it follows that for every $\epsilon \in R, \epsilon > 0,$ there exists positive integer $N$ depending on $\epsilon$ and such, that

$$\|a_n - b_n\| < \epsilon$$

for every $n > N$. From the inequality (2.1.1), it follows that for every $\epsilon \in R, \epsilon > 0,$ there exists positive integer $N$ depending on $\epsilon$ and such, that

(2.5.44) $$|\, \|a_n\| - \|b_n\|\, | \leq \|a_n - b_n\| < \epsilon$$

for every $n > N$. According to the definition 2.1.17, it follows that

(2.5.45) $$\lim_{n \to \infty} (\|a_n\| - \|b_n\|) = 0$$

From the lemma 2.5.12, it follows that sequences $\|a_n\|$ and $\|b_n\|$ are fundamental sequences. The equation (2.5.42) follows from the theorem 2.1.22. □

LEMMA 2.5.14. *We define the norm on Ω-group $B$ by equation*

(2.5.46) $$\| [a_n] \| = \lim_{n \to \infty} \|a_n\|$$

PROOF. From lemmas 2.5.12, 2.5.13, it follows that the equation (2.5.46) is a proper definition of the map

$$B \to R$$

2.5.14.1: According to the statement 2.1.9.1, $\|a_n\| \geq 0$. Therefore,[2.22]

$$\lim_{n \to \infty} \|a_n\| \geq 0 \qquad (2.5.47)$$

From the equation (2.5.46) and the inequality (2.5.47), it follows that $\|[a_n]\| \geq 0$. Therefore, the statement 2.1.9.1 is true for the map (2.5.46).

2.5.14.2: Let $\|[a]\| = 0$. From the equation (2.5.46) it is follows that

$$\lim_{n \to \infty} \|a_n\| = 0 \qquad (2.5.48)$$

According to the definition 2.1.17, from the equation (2.5.48), it follows that for every $\epsilon \in R, \epsilon > 0$, there exists positive integer $N$ depending on $\epsilon$ and such, that $\|a_n\| < \epsilon$ for every $n > N$. According to the definition 2.1.17,

$$\lim_{n \to \infty} a_n = \lim_{n \to \infty} (a_n - 0) = 0 \qquad (2.5.49)$$

From the equations (2.5.49), (2.5.7), it follows that $a_n \sim 0$. According to the definition (2.5.14), $[a_n] = [0]$. Therefore, the statement 2.1.9.2 is true for the map (2.5.46).

2.5.14.3: According to the statement 2.1.9.3,

$$\|a_n + b_n\| \leq \|a_n\| + \|b_n\|$$

Therefore,[2.23]

$$\lim_{n \to \infty} \|a_n + b_n\| \leq \lim_{n \to \infty} \|a_n\| + \lim_{n \to \infty} \|b_n\| \qquad (2.5.50)$$

From equations (2.5.46) and the inequality (2.5.50), it follows that

$$\|[a_n + b_n]\| \leq \|[a_n]\| + \|[b_n]\| \qquad (2.5.51)$$

From the inequality (2.5.51) and the equation (2.5.15) it follows that

$$\|[a_n] + [b_n]\| \leq \|[a_n]\| + \|[b_n]\|$$

Therefore, the statement 2.1.9.3 is true for the map (2.5.46).

2.5.14.4: According to the statement 2.1.9.4, $\|-a_n\| = \|a_n\|$. According to the theorem 2.1.22,

$$\lim_{n \to \infty} \|-a_n\| = \lim_{n \to \infty} \|a_n\| \qquad (2.5.52)$$

From equations (2.5.46), (2.5.52), it follows that

$$\|[-a_n]\| = \|[a_n]\| \qquad (2.5.53)$$

From equation (2.5.53) and statement 2.5.7.4, it follows that

$$\|-[a_n]\| = \|[a_n]\|$$

Therefore, the statement 2.1.9.4 is true for the map (2.5.46).

According to statements 2.5.14.1, 2.5.14.2, 2.5.14.3, 2.5.7.4 and the definition 2.1.9, the map (2.5.46) is a norm on Ω-group $B$.  □

---

[2.22] See, for instance, [10], page 50.
[2.23] See, for instance, [10], page 50.

LEMMA 2.5.15. Ω-group $A$ is subgroup of Ω-group $B$. Conditions $a \in A$ and
$$a = \lim_{n \to \infty} a_n$$
imply

(2.5.54) $$\|a\| = \|[a_n]\|$$

PROOF. Let $a \in A$. According to definitions 2.1.17, 2.1.19, the sequence $a_n = a$, $n = 1, ...,$ is fundamental sequence and converges to $a$. According to the definition (2.5.7) and the theorem 2.1.22, condition $a_n \sim b_n$ implies
$$a = \lim_{n \to \infty} b_n$$
So we can identify $a$ and $[a_n]$. According to theorems 2.4.3, 2.4.6, this identification is homomorphism of Ω-group $A$ into Ω-group $B$. From the equation (2.5.46) and the theorem 2.3.9, it follows that

(2.5.55) $$\|[a_n]\| = \lim_{n \to \infty} \|a_n\| = \|a\|$$

The equation (2.5.54) follows from the equation (2.5.55). □

LEMMA 2.5.16. For any $B$-number $a$, there exist a sequence of $A$-numbers such that

(2.5.56) $$a = \lim_{n \to \infty} a_n$$

PROOF. There exists a fundamental sequence of $A$-numbers $a_n$ such that $a = [a_n]$. According to the definition 2.1.19, for given $\epsilon \in R$, $\epsilon > 0$, there exists positive integer $N_1$ depending on $\epsilon$ and such, that

(2.5.57) $$\|a_p - a_q\| < \frac{\epsilon}{2}$$

for every $p, q > N_1$. According to lemmas 2.5.5, 2.5.15, for any $q > N_1$

(2.5.58) $$a - a_q = [a_p] - [a_q] = [a_p - a_q]$$

According to the lemma 2.5.14

(2.5.59) $$\|a - a_q\| = \lim_{p \to \infty} \|a_p - a_q\|$$

According to the definition 2.1.17, from the equation (2.5.59), it follows that for every $\epsilon \in R, \epsilon > 0,$ there exists positive integer $N_2$ depending on $\epsilon$ and such, that

(2.5.60) $$|\|a - a_q\| - \|a_p - a_q\|| < \frac{\epsilon}{2}$$

for every $p > N_2$. The inequality[2.24]

(2.5.63) $$\|a - a_q\| - \|a_p - a_q\| < \frac{\epsilon}{2}$$

---

[2.24]The inequality

(2.5.61) $$\|a - a_q\| - \|a_p - a_q\| > -\frac{\epsilon}{2}$$

also follows from the inequality (2.5.60). The inequality

(2.5.62) $$\|a - a_q\| > \|a_p - a_q\| - \frac{\epsilon}{2}$$

follows from the inequality (2.5.61). From inequalities (2.5.57), (2.5.62) and the statement 2.1.9.1 it follows that
$$\|a - a_q\| \geq 0$$
However, this statement is evident and is not interesting for us.

follows from the inequality (2.5.60). Let

$$N = \max(N_1, N_2)$$

From inequalities (2.5.57), (2.5.63), it follows that for given $\epsilon \in R$, $\epsilon > 0$, there exists positive integer $N$ depending on $\epsilon$ and such, that

(2.5.64) $$\|a - a_q\| < \|a_p - a_q\| + \frac{\epsilon}{2} < \epsilon$$

for every $q > N$. The equation (2.5.56) follows from the definition 2.1.17. □

THEOREM 2.5.17. *There exists completion of normed Ω-group A.*

PROOF. According to the lemma 2.5.15, Ω-group $A$ is subgroup of Ω-group $B$. From the lemma 2.5.16 and the theorem 2.2.14 it follows that Ω-group $A$ is everywhere dense in Ω-group $B$. According to the definition 2.5.1, to prove the theorem, we need to prove that Ω-group $B$ is complete.

Let $b_n$ be fundamental sequence of $B$-numbers. According to definitions 2.1.14, 2.2.3, 2.2.5, for any $n > 0$ there exists $a_n \in A$ such that

(2.5.65) $$\|a_n - b_n\| < \frac{1}{n}$$

The equation

(2.5.66) $$\lim_{n \to \infty} (a_n - b_n) = 0$$

follows from the definition 2.1.17. From the equation (2.5.66) and the theorem 2.1.21, it follows that the sequence of $A$-numbers $a_n$ is fundamental sequence. According to the definition (2.5.14), there exists $B$-number $[a_n]$. According to the theorem 2.1.22 and the lemma 2.5.16

(2.5.67) $$\lim_{n \to \infty} b_n = \lim_{n \to \infty} a_n = [a_n]$$

That Ω-group $B$ is complete follows from equation (2.5.67). □

## 2.6. Ω-Group of Maps

THEOREM 2.6.1. *Let $M(X, A)$ be set of maps of a set $X$ to Ω-group $A$. We can define the structure of Ω-group on the set $M(X, A)$.*

PROOF. Let $f, g \in M(X, A)$. Then we assume

$$(f + g)(x) = f(x) + g(x)$$

Let $\omega \in \Omega$ be $n$-ari operation. For maps $f_i \in M(X, A)$, $i = 1, ..., n$, we assume

(2.6.1) $$(f_1...f_n\omega)(x) = f_1(x)...f_n(x)\omega$$

□

Since $X$ is an arbitrary set, we cannot define norm in Ω-group $M(X, A)$. However we can define the convergence of the sequence in $M(X, A)$; therefore, we can define a topology in $M(X, A)$.

DEFINITION 2.6.2. *Let $f_n \in M(X, A)$, $n = 1, ...$, be sequence of maps into normed Ω-group $A$. The map $f \in M(X, A)$ is called* **limit of sequence** $f_n$*, if for any $x \in X$*

$$f(x) = \lim_{n \to \infty} f_n(x)$$

*We also say that* **sequence $f_n$ converges** *to the map $f$.* □

From definitions 2.1.17, 2.6.2, it follows that since sequence $f_n$ converges to $f$, then for any $\epsilon \in R$, $\epsilon > 0$, there exists $N(x)$ such that
$$\|f_n(x) - f(x)\| < \epsilon$$
for any $n > N(x)$.

DEFINITION 2.6.3. Let $f_n \in M(X, A)$, $n = 1, ...,$ be sequence of maps into normed Ω-group $A$. **Sequence $f_n$ converges uniformly** to the map $f$, if for any $\epsilon \in R$, $\epsilon > 0$, there exists $N$ such that
$$\|f_n(x) - f(x)\| < \epsilon$$
for any $n > N$.  □

THEOREM 2.6.4. *Sequence of maps $f_n \in M(X, A)$, $n = 1, ...,$ into normed Ω-group $A$ converges uniformly to the map $f$, if for any $\epsilon \in R$, $\epsilon > 0$, there exists $N$ such that $f_n(x) \in B_o(f(x), \epsilon)$ for any $n > N$.*

PROOF. The theorem follows from definitions 2.1.14, 2.6.3.  □

THEOREM 2.6.5. *Sequence of maps $f_n \in M(X, A)$, $n = 1, ...,$ into normed Ω-group $A$ converges uniformly to the map $f$, if for any $\epsilon \in R$, $\epsilon > 0$, there exists $N$ such that*

(2.6.2) $$\|f_n(x) - f_m(x)\| < \epsilon$$

*for any $n, m > N$.*

PROOF. According to the definition 2.6.3, for any $\epsilon \in R$, $\epsilon > 0$, there exists $N$ such that

(2.6.3) $$\|f_n(x) - f(x)\| < \frac{\epsilon}{2}$$

for any $n > N$. For any $n, m > N$,

(2.6.4) $$\begin{aligned}\|f_n(x) - f_m(x)\| &= \|f_n(x) - f(x) + f(x) - f_m(x)\| \\ &\leq \|f_n(x) - f(x)\| + \|f_m(x) - f(x)\| < \epsilon\end{aligned}$$

follows from the statement 2.1.9.3 and the inequality (2.6.3). The inequality (2.6.2) follows from the inequality (2.6.4).  □

THEOREM 2.6.6. *Sequence of maps $f_n \in M(X, A)$, $n = 1, ...,$ into normed Ω-group $A$ converges uniformly to the map $f$, if for any $\epsilon \in R$, $\epsilon > 0$, there exists $N$ such that $f_n(x) \in B_o(f_m(x), \epsilon)$ for any $n, m > N$.*

PROOF. The theorem follows from the definition 2.1.14 and from the theorem 2.6.5.  □

THEOREM 2.6.7. *Let sequence of maps $f_n \in M(X, A)$, $n = 1, ...,$ into complete Ω-group $A$ converge uniformly to the map $f$. Let sequence of maps $g_n \in M(X, A)$, $n = 1, ...,$ into complete Ω-group $A$ converge uniformly to the map $g$. Then sequence of maps*

$$h_n = f_n + g_n$$

*into complete Ω-group $A$ converges uniformly to the map*

(2.6.5) $$h = f + g$$

PROOF. From the equation
$$f(x) = \lim_{n \to \infty} f_n(x)$$
and the theorem 2.6.4, it follows that for given $\epsilon \in R$, $\epsilon > 0$, there exists $N_a$ such that the condition $n > N_a$ implies that
$$(2.6.6) \qquad f_n(x) \in B_o(f(x), \epsilon/2)$$
From the equation
$$g(x) = \lim_{n \to \infty} g_n(x)$$
and the theorem 2.6.4, it follows that for given $\epsilon \in R$, $\epsilon > 0$, there exists $N_b$ such that the condition $n > N_b$ implies that
$$(2.6.7) \qquad g_n(x) \in B_o(g(x), \epsilon/2)$$
Let
$$N = \max(N_a, N_b)$$
From equations (2.6.6), (2.6.7), the theorem 2.4.2 and condition $n > N$, it follows that
$$(2.6.8) \qquad h_n(x) = f_n(x) + g_n(x) \in B_o(f(x) + g(x), \epsilon/2 + \epsilon/2) = B_o(h(x), \epsilon)$$
The equation (2.6.5) follows from the equation (2.6.8) and the theorem 2.6.4. □

THEOREM 2.6.8. *Let $A$ be complete $\Omega$-group. Let $\omega \in \Omega$ be n-ari operation. Let sequence of maps $f_{i \cdot m} \in M(X, A)$, $i = 1, ..., n$, $m = 1, ...,$ into complete $\Omega$-group $A$ converge uniformly to the map $f_i$. Let range of the map $f_i$ be compact set. Then sequence of maps*
$$h_m = f_{1 \cdot m} ... f_{n \cdot m} \omega$$
*into complete $\Omega$-group $A$ converges uniformly to the map*
$$(2.6.9) \qquad h = f_1 ... f_n \omega$$

PROOF. From the equation
$$f_i(x) = \lim_{m \to \infty} f_{i \cdot m}(x)$$
and the theorem 2.6.4, it follows that for given
$$(2.6.10) \qquad \delta_1 \in R, \ \delta_1 > 0$$
there exists $M_i$ such that the condition $m > M_i$ implies that
$$(2.6.11) \qquad f_{i \cdot m}(x) \in B_o(f_i(x), \delta_1)$$
Let
$$M = \max(M_1, ..., M_n)$$
From equations (2.6.11), the theorem 2.4.5 and condition $m > M$, it follows that
$$(2.6.12) \qquad h_m(x) = f_{1 \cdot m}(x) ... f_{n \cdot m}(x) \omega \in B_o(f_1(x) ... f_n(x) \omega, \epsilon_1) = B_o(h(x), \epsilon_1)$$
$$(2.6.13) \qquad \epsilon_1 = \delta_1 C_2 ... C_n + ... + \delta_1 C_1 ... C_{n-1}$$
where
$$(2.6.14) \qquad C_i = \|f_i(x)\| + \delta_1$$
The value of $C_i$ depends on $x$. So $\epsilon_1$ also depends on $x$. To estimate uniform convergence, we need to choose a maximum value of $\epsilon_1$. For given $\delta_1$, value of $\epsilon_1$ is polylinear map of values of $C_i$. Therefore, $\epsilon_1$ has maximal value when each $C_i$ has maximal value. Since we consider the estimation of the value of $\epsilon_1$ from the right, whether $C_i$ acquire maximal values at the same time is of no concern to us. Since the range of the map $f_i$ is compact set, then, according to the theorem 2.3.11, the following value is defined

$$F_i = \sup \|f_i(x)\| \tag{2.6.15}$$

From the equation (2.6.15) and the statement 2.1.9.1, it follows that

$$F_i \geq 0 \tag{2.6.16}$$

From equations (2.6.14), (2.6.15), it follows that we can put

$$C_i = F_i + \delta_1 \tag{2.6.17}$$

From the equation (2.6.17) and statements (2.6.16), (2.6.10) it follows that

$$C_i > 0 \quad \frac{dC_i}{d\delta_1} > 0 \tag{2.6.18}$$

From equations (2.6.13), (2.6.17) and the statement (2.6.18), it follows that $\epsilon_1$ is polynomial strictly monotone increasing function of $\delta_1$ such that

$$\delta_1 = 0 \Rightarrow \epsilon_1 = 0$$

According to the theorem 2.3.5, the map (2.6.13) maps the interval $[0, \delta_1)$ into the interval $[0, \epsilon_1)$. According to the theorem 2.3.3, for given $\epsilon > 0$ there exist $\delta > 0$ such that

$$\epsilon_1(\delta) < \epsilon$$

According to construction, a value of $M$ depends on a value of $\delta_1$. We choose the value of $M$ corresponding to $\delta_1 = \delta$. Therefore, for given $\epsilon \in R, \epsilon > 0$, there exists $M$ such that the condition $m > M$ implies that

$$h_m(x) = f_{1 \cdot m}(x)...f_{n \cdot m}(x)\omega \in B_o(f_1(x)...f_n(x)\omega, \epsilon) = B_o(h(x), \epsilon) \tag{2.6.19}$$

The equation (2.6.9) follows from the equation (2.6.19) and the theorem 2.6.4. $\square$

CHAPTER 3

# Representation of $\Omega$-Group

### 3.1. Representation of $\Omega$-Group

DEFINITION 3.1.1. Let
$$f : A_1 \dashrightarrow A_2$$
be representation[3.1] of $\Omega_1$-group $A_1$ with norm $\|x\|_1$ in $\Omega_2$-group $A_2$ with norm $\|x\|_2$. The value

(3.1.1) $$\|f\| = \sup \frac{\|f(a_1)(a_2)\|_2}{\|a_1\|_1 \|a_2\|_2}$$

is called **norm of representation** $f$. □

THEOREM 3.1.2. *Let*
$$f : A_1 \dashrightarrow A_2$$
*be representation of $\Omega_1$-group $A_1$ with norm $\|x\|_1$ in $\Omega_2$-group $A_2$ with norm $\|x\|_2$. Then*

(3.1.2) $$\|f(a_1)(a_2)\|_2 \leq \|f\| \|a_1\|_1 \|a_2\|_2$$

PROOF. From the equation (3.1.1), it follows that

(3.1.3) $$\frac{\|f(a_1)(a_2)\|_2}{\|a_1\|_1 \|a_2\|_2} \leq \sup \frac{\|f(a_1)(a_2)\|_2}{\|a_1\|_1 \|a_2\|_2} = \|f\|$$

The inequality (3.1.2) follows from the inequality (3.1.3). □

THEOREM 3.1.3. *Let*
$$f : A_1 \dashrightarrow A_2$$
*be representation of $\Omega_1$-group $A_1$ with norm $\|x\|_1$ in $\Omega_2$-group $A_2$ with norm $\|x\|_2$. The following inequality is true*

(3.1.4) $$\|f(c_1)(c_2) - f(a_1)(a_2)\|_2 \leq \|f\|(\|c_1 - a_1\|_1 C_2 + C_1 \|c_2 - a_2\|_2)$$

*where*

(3.1.5) $$C_1 = \max(\|a_1\|_1, \|c_1\|_1) \quad C_2 = \max(\|a_2\|_2, \|c_2\|_2)$$

PROOF. According to the definitions 2.1.3, [7]-2.1.4,

(3.1.6) $$\begin{aligned} f(c_1)(c_2) - f(a_1)(a_2) &= f(c_1)(c_2) - f(a_1)(c_2) + f(a_1)(c_2) - f(a_1)(a_2) \\ &= f(c_1 - a_1)(c_2) + f(a_1)(c_2 - a_2) \end{aligned}$$

---

[3.1]See the definition [7]-2.1.4 of the representation of universal algebra. According to the definition [9]-2.1.2, module is representation of ring in Abelian group. Since ring and Abelian group are $\Omega$-groups, then module is representation of $\Omega$-group.

From the equation (3.1.6) and the statement 2.1.9.3, it follows that

(3.1.7) $$\|f(c_1)(c_2) - f(a_1)(a_2)\|_2 \leq \|f(c_1 - a_1)(c_2)\|_2 + \|f(a_1)(c_2 - a_2)\|_2$$

From the equation (3.1.5) and the theorem 3.1.2, it follows that

(3.1.8) $$\begin{aligned}\|f(c_1 - a_1)(c_2)\|_2 &\leq \|f\| \, \|c_1 - a_1\|_1 \, C_2 \\ \|f(a_1)(c_2 - a_2)\|_2 &\leq \|f\| \, C_1 \, \|c_2 - a_2\|_2\end{aligned}$$

The statement (3.1.4) follows from inequalities (3.1.7), (3.1.8). $\square$

THEOREM 3.1.4. *Let*
$$f : A_1 \relbar\joinrel\twoheadrightarrow A_2$$
*be representation of $\Omega_1$-group $A_1$ with norm $\|x\|_1$ in $\Omega_2$-group $A_2$ with norm $\|x\|_2$. For $c_1 \in A_1$, let $c_1 \in B_c(a_1 \in A_1, R_1)$. For $c_2 \in A_2$, let $c_2 \in B_c(a_2 \in A_2, R_2)$. Then*

(3.1.9) $$f(c_1)(c_2) \in B_c(f(a_1)(a_2), \|f\|(R_1 C_2 + C_1 R_2))$$

*where*

(3.1.10) $$C_1 = \|a_1\|_1 + R_1 \quad C_2 = \|a_2\|_2 + R_2$$

PROOF. According to the definition 2.1.14

(3.1.11) $$\begin{aligned}\|c_1 - a_1\|_1 &\leq R_1 \\ \|c_2 - a_2\|_2 &\leq R_2\end{aligned}$$

The inequality

(3.1.12) $$\|c_i\| = \|a_i + c_i - a_i\| \leq \|a_i\| + \|c_i - a_i\| \leq \|a_i\| + R_i$$

follows from the inequality (3.1.9) and from the statement 2.1.9.3. The inequality

(3.1.13) $$C_i \leq \max(\|a_i\|, \|a_i\| + R_i) = \|a_i\| + R_i \quad i = 1, 2$$

follows from the inequality (3.1.12) and the equation (3.1.5). The equation (3.1.10) follows from the inequality (3.1.13). The statement (3.1.9) follows from inequalities (3.1.4), (3.1.11), the equation (3.1.10) and the definition 2.1.14. $\square$

THEOREM 3.1.5. *Let*
$$f : A_1 \relbar\joinrel\twoheadrightarrow A_2$$
*be representation of $\Omega_1$-group $A_1$ with norm $\|x\|_1$ in $\Omega_2$-group $A_2$ with norm $\|x\|_2$. Let the sequence of $A_1$-numbers $a_{1 \cdot n}$ converge and*

(3.1.14) $$\lim_{n \to \infty} a_{1 \cdot n} = a_1$$

*Let the sequence of $A_2$-numbers $a_{2 \cdot n}$ converge and*

(3.1.15) $$\lim_{n \to \infty} a_{2 \cdot n} = a_2$$

*Then the sequence of $A_2$-numbers $f(a_{1 \cdot n})(a_{2 \cdot n})$ converges and*

(3.1.16) $$\lim_{n \to \infty} f(a_{1 \cdot n})(a_{2 \cdot n}) = f(a_1)(a_2)$$

PROOF. Let

(3.1.17) $$\delta_1 \in R, \ \delta_1 > 0$$

From the equation (3.1.14) and the theorem 2.1.18, it follows that for given $\delta_1$ there exists $N_1$ such that the condition $n > N_1$ implies that

(3.1.18) $$a_{1 \cdot n} \in B_o(a_1, \delta_1)$$

From the equation (3.1.15) and the theorem 2.1.18, it follows that for given $\delta_1$ there exists $N_2$ such that the condition $n > N_2$ implies that

(3.1.19) $$a_{2 \cdot n} \in B_o(a_2, \delta_1)$$

Let
$$N = \max(N_1, N_2)$$

From equations (3.1.18), (3.1.19), the theorem 3.1.4 and condition $n > N$, it follows that

(3.1.20) $$f(a_{1 \cdot n})(a_{2 \cdot n}) \in B_o(f(a_1)(a_n), \epsilon_1)$$

(3.1.21) $$\epsilon_1 = \|f\|(\delta_1 C_2 + \delta_1 C_1)$$

where

(3.1.22) $$C_i = \|a_i\| + \delta_1$$

From the equation (3.1.22) and statements 2.1.9.1, (3.1.17) it follows that

(3.1.23) $$C_i > 0 \quad \frac{dC_i}{d\delta_1} > 0$$

From equations (3.1.21), (3.1.22) and the statement (3.1.23), it follows that $\epsilon_1$ is polynomial strictly monotone increasing function of $\delta_1$ such that

$$\delta_1 = 0 \Rightarrow \epsilon_1 = 0$$

According to the theorem 2.3.5, the map (3.1.21) maps the interval $[0, \delta_1)$ into the interval $[0, \epsilon_1)$. According to the theorem 2.3.3, for given $\epsilon > 0$ there exist $\delta > 0$ such that

$$\epsilon_1(\delta) < \epsilon$$

According to construction, a value of $N$ depends on a value of $\delta_1$. We choose the value of $N$ corresponding to $\delta_1 = \delta$. Therefore, for given $\epsilon \in R, \epsilon > 0$, there exists $N$ such that the condition $n > N$ implies that

(3.1.24) $$f(a_{1 \cdot n})(a_{2 \cdot n}) \in B_o(f(a_1)(a_n), \epsilon)$$

The equation (3.1.16) follows from the equation (3.1.24) and the theorem 2.1.18. □

THEOREM 3.1.6. *The representation*

$$f : A_1 \relbar\joinrel\twoheadrightarrow A_2$$

*of $\Omega_1$-group $A_1$ with norm $\|x\|_1$ in $\Omega_2$-group $A_2$ with norm $\|x\|_2$ can be extended to representation*

$$f^* : A_1 \relbar\joinrel\twoheadrightarrow M(X, A_2)$$

*of $\Omega_1$-group $A_1$ in $\Omega_2$-group $M(X, A_2)$ where $(g \in M(X, A_2))$*

(3.1.25) $$(f^*(a_1)(g))(x) = f(a_1)(g(x))$$

PROOF. To prove the theorem, we must show that the map $f^*$ is homomorphism of $\Omega_1$-group $A_1$ and for any $a_1 \in A_1$ a map $f^*(a_1)$ is homomorphism of $\Omega_2$-group $M(X, A_2)$.

According to the definition [7]-2.1.4, the map $f$ is homomorphism of $\Omega_1$-group $A_1$. Therefore, for $n$-ari operation $\omega \in \Omega_1$ and any $a_{1 \cdot 1}, ..., a_{1 \cdot n} \in A_1$, following equation is true

$$f(a_{1 \cdot 1}...a_{1 \cdot n}\omega) = f(a_{1 \cdot 1})...f(a_{1 \cdot n})\omega \tag{3.1.26}$$

From equations (2.6.1), (3.1.26), it follows that

$$f(a_{1 \cdot 1}...a_{1 \cdot n}\omega)(a_2) = (f(a_{1 \cdot 1})(a_2))...(f(a_{1 \cdot n})(a_2))\omega \tag{3.1.27}$$

From equations (3.1.25), (3.1.27), it follows that

$$\begin{aligned}(f^*(a_{1 \cdot 1}...a_{1 \cdot n}\omega)(g))(x) &= f(a_{1 \cdot 1}...a_{1 \cdot n}\omega)(g(x)) \\ &= f(a_{1 \cdot 1})(g(x))...f(a_{1 \cdot n})(g(x))\omega \\ &= ((f^*(a_{1 \cdot 1})(g))(x))...((f^*(a_{1 \cdot n})(g))(x))\omega\end{aligned} \tag{3.1.28}$$

From equations (2.6.1), (3.1.28), it follows that

$$\begin{aligned}(f^*(a_{1 \cdot 1}...a_{1 \cdot n}\omega)(g))(x) &= ((f^*(a_{1 \cdot 1})(g))...(f^*(a_{1 \cdot n})(g))\omega)(x) \\ &= ((f^*(a_{1 \cdot 1})...f^*(a_{1 \cdot n})\omega)(g))(x)\end{aligned} \tag{3.1.29}$$

From the equation (3.1.29), it follows that

$$f^*(a_{1 \cdot 1}...a_{1 \cdot n}\omega) = f^*(a_{1 \cdot 1})...f^*(a_{1 \cdot n})\omega \tag{3.1.30}$$

From the equation (3.1.30) it follows that the map $f^*$ is homomorphism of $\Omega_1$-group $A_1$.

According to the definition [7]-2.1.4, the map $f(a_1)$, $a_1 \in A_1$, is homomorphism of $\Omega_2$-group $A_2$. For $n$-ari operation $\omega \in \Omega_2$ and any $a_{2 \cdot 1}, ..., a_{2 \cdot n} \in A_2$, following equation is true

$$f(a_1)(a_{2 \cdot 1}...a_{2 \cdot n}\omega) = (f(a_1)(a_{2 \cdot 1}))...(f(a_1)(a_{2 \cdot n}))\omega \tag{3.1.31}$$

From equations (3.1.25), it follows that

$$(f^*(a_1)(g_1...g_n\omega))(x) = f(a_1)((g_1...g_n\omega)(x)) \tag{3.1.32}$$

From equations (2.6.1), (3.1.32), it follows that

$$(f^*(a_1)(g_1...g_n\omega))(x) = f(a_1)(g_1(x)...g_n(x)\omega) \tag{3.1.33}$$

From equations (3.1.31), (3.1.33), it follows that

$$(f^*(a_1)(g_1...g_n\omega))(x) = (f(a_1)(g_1(x)))...(f(a_1)(g_n(x)))\omega \tag{3.1.34}$$

From equations (2.6.1), (3.1.25), (3.1.34), it follows that

$$\begin{aligned}(f^*(a_1)(g_1...g_n\omega))(x) &= ((f^*(a_1)(g_1))(x))...((f^*(a_1)(g_n))(x))\omega \\ &= ((f^*(a_1)(g_1))...(f^*(a_1)(g_n))\omega)(x)\end{aligned} \tag{3.1.35}$$

From the equation (3.1.35), it follows that

$$f^*(a_1)(g_1...g_n\omega) = (f^*(a_1)(g_1))...(f^*(a_1)(g_n))\omega \tag{3.1.36}$$

From the equation (3.1.36), it follows that for any $a_1 \in A_1$ a map $f^*(a_1)$ is homomorphism of $\Omega_2$-group $M(X, A_2)$. □

## 3.2. Representation of $\Omega$-Group of Maps

THEOREM 3.2.1. *The representation*

$$f : A_1 \relbar\joinrel\twoheadrightarrow A_2$$

*of $\Omega_1$-group $A_1$ with norm $\|x\|_1$ in $\Omega_2$-group $A_2$ with norm $\|x\|_2$ generates representation*

$$f_X : M(X, A_1) \relbar\joinrel\twoheadrightarrow M(X, A_2)$$

*of $\Omega_1$-group $M(X, A_1)$ in $\Omega_2$-group $M(X, A_2)$ where ( $g_1 \in M(X, A_1)$ , $g_2 \in M(X, A_2)$ )*

(3.2.1)
$$\begin{aligned} f_X(g_1)(g_2) &: X -> A_2 \\ (f_X(g_1)(g_2))(x) &= f(g_1(x))(g_2(x)) \end{aligned}$$

PROOF. To prove the theorem, we must show that the map $f_X$ is homomorphism of $\Omega_1$-group $M(X, A_1)$ and, for any $g_1 \in M(X, A_1)$, a map $f_X(g_1)$ is homomorphism of $\Omega_2$-group $M(X, A_2)$.

According to the definition [7]-2.1.4, the map $f$ is homomorphism of $\Omega_1$-group $A_1$. Therefore, for $n$-ari operation $\omega \in \Omega_1$ and any $a_{1 \cdot 1}, ..., a_{1 \cdot n} \in A_1$, following equation is true

(3.2.2) $$f(a_{1 \cdot 1}...a_{1 \cdot n}\omega) = f(a_{1 \cdot 1})...f(a_{1 \cdot n})\omega$$

From the equation (2.6.1), it follows that

(3.2.3) $$(g_{1 \cdot 1}...g_{1 \cdot n}\omega)(x) = g_{1 \cdot 1}(x)...g_{1 \cdot n}(x)\omega$$

From equations (3.2.2), (3.2.3), it follows that

(3.2.4) $$f((g_{1 \cdot 1}...g_{1 \cdot n}\omega)(x)) = f(g_{1 \cdot 1}(x)...g_{1 \cdot n}(x)\omega) = f(g_{1 \cdot 1}(x))...f(g_{1 \cdot n}(x))\omega$$

From equations (2.6.1), (3.2.4), it follows that

(3.2.5)
$$\begin{aligned} f((g_{1 \cdot 1}...g_{1 \cdot n}\omega)(x))(a_2) &= (f(g_{1 \cdot 1}(x))...f(g_{1 \cdot n}(x))\omega)(a_2) \\ &= (f(g_{1 \cdot 1}(x))(a_2))...(f(g_{1 \cdot n}(x))(a_2))\omega \end{aligned}$$

From equations (3.2.1), (3.2.5), it follows that

(3.2.6)
$$\begin{aligned} (f_X(g_{1 \cdot 1}...g_{1 \cdot n}\omega)(g_2))(x) &= f((g_{1 \cdot 1}...g_{1 \cdot n}\omega)(x))(g_2(x)) \\ &= (f(g_{1 \cdot 1}(x))(g_2(x)))...(f(g_{1 \cdot n}(x))(g_2(x)))\omega \\ &= ((f_X(g_{1 \cdot 1})(g_2))(x))...((f_X(g_{1 \cdot n})(g_2))(x))\omega \end{aligned}$$

From equations (2.6.1), (3.2.6), it follows that

(3.2.7)
$$\begin{aligned} (f_X(g_{1 \cdot 1}...g_{1 \cdot n}\omega)(g_2))(x) &= ((f_X(g_{1 \cdot 1})(g_2))...(f_X(g_{1 \cdot n})(g_2))\omega)(x) \\ &= ((f_X(g_{1 \cdot 1})...f_X(g_{1 \cdot n})\omega)(g_2))(x) \end{aligned}$$

From the equation (3.2.7), it follows that

(3.2.8) $$f_X(g_{1 \cdot 1}...g_{1 \cdot n}\omega) = f_X(g_{1 \cdot 1})...f_X(g_{1 \cdot n})\omega$$

From the equation (3.2.8) it follows that the map $f_X$ is homomorphism of $\Omega_1$-group $M(X, A_1)$.

According to the definition [7]-2.1.4, the map $f(a_1)$, $a_1 \in A_1$, is homomorphism of $\Omega_2$-group $A_2$. For $n$-ari operation $\omega \in \Omega_2$ and any $a_{2 \cdot 1}, ..., a_{2 \cdot n} \in A_2$, following equation is true

(3.2.9) $$f(a_1)(a_{2 \cdot 1}...a_{2 \cdot n}\omega) = (f(a_1)(a_{2 \cdot 1}))...(f(a_1)(a_{2 \cdot n}))\omega$$

From the equation (2.6.1), it follows that

$$(g_{2\cdot 1}...g_{2\cdot n}\omega)(x) = g_{2\cdot 1}(x)...g_{2\cdot n}(x)\omega \tag{3.2.10}$$

From equations (3.2.9), (3.2.10), it follows that

$$\begin{aligned} f(a_1)((g_{2\cdot 1}...g_{2\cdot n}\omega)(x)) &= f(a_1)(g_{2\cdot 1}(x)...g_{2\cdot n}(x)\omega) \\ &= (f(a_1)(g_{2\cdot 1}(x)))...(f(a_1)(g_{2\cdot n}(x)))\omega \end{aligned} \tag{3.2.11}$$

From equations (3.2.1), (3.2.11), it follows that

$$\begin{aligned} (f_X(g_1)(g_{2\cdot 1}...g_{2\cdot n}\omega))(x) &= f(g_1(x))((g_{2\cdot 1}...g_{2\cdot n}\omega)(x)) \\ &= (f(g_1(x))(g_{2\cdot 1}(x)))...(f(g_1(x))(g_{2\cdot n}(x)))\omega \\ &= ((f_X(g_1)(g_{2\cdot 1}))(x))...((f(g_1)(g_{2\cdot n}))(x))\omega \end{aligned} \tag{3.2.12}$$

From equations (2.6.1), (3.2.12), it follows that

$$(f_X(g_1)(g_{2\cdot 1}...g_{2\cdot n}\omega))(x) = ((f_X(g_1)(g_{2\cdot 1}))...(f_X(g_1)(g_{2\cdot n}))\omega)(x) \tag{3.2.13}$$

From equations (3.2.9), (3.2.13), it follows that

$$f_X(g_1)(g_{2\cdot 1}...g_{2\cdot n}\omega) = (f_X(g_1)(g_{2\cdot 1}))...(f_X(g_1)(g_{2\cdot n}))\omega \tag{3.2.14}$$

From the equation (3.2.14), it follows that, for any $g_1 \in M(X, A_1)$, a map $f_X(g_1)$ is homomorphism of $\Omega_2$-group $M(X, A_2)$. □

THEOREM 3.2.2. *Let*

$$f : A_1 \relbar\joinrel\twoheadrightarrow A_2$$

*be representation of complete $\Omega_1$-group $A_1$ with norm $\|x\|_1$ in complete $\Omega_2$-group $A_2$ with norm $\|x\|_2$. Let sequence of maps $g_{1\cdot n} \in M(X, A_1)$, $n = 1, ...,$ converge uniformly to the map $g_1$. Let sequence of maps $g_{2\cdot n} \in M(X, A_2)$, $n = 1, ...,$ converge uniformly to the map $g_2$. Let range of the map $g_i$, $i = 1, 2$, be compact set. Then the sequence of maps $f_X(g_{1\cdot n})(g_{2\cdot n})$ converge uniformly to the map $f_X(g_1)(g_2)$.*

PROOF. From the equation

$$g_i(x) = \lim_{n \to \infty} g_{i\cdot n}(x) \tag{3.2.15}$$

and the theorem 2.6.4, it follows that for given

$$\delta_1 \in R, \ \delta_1 > 0 \tag{3.2.16}$$

there exists $N_i$ such that the condition $n > N_i$ implies that

$$g_{i\cdot n}(x) \in B_o(g_i(x), \delta_1) \tag{3.2.17}$$

Let

$$N = \max(N_1, N_2)$$

From the equation (3.2.17), the theorem 3.1.4 and condition $n > N$, it follows that

$$f(g_{1\cdot n}(x))(g_{2\cdot n}(x)) \in B_o(f(g_1(x))(g_2(x)), \epsilon_1) \tag{3.2.18}$$

$$\epsilon_1 = \|f\|(\delta_1 C_2 + \delta_1 C_1) \tag{3.2.19}$$

where

$$C_i = \|f_i(x)\| + \delta_1 \tag{3.2.20}$$

The value of $C_i$ depends on $x$. So $\epsilon_1$ also depends on $x$. To estimate uniform convergence, we need to choose a maximum value of $\epsilon_1$. For given $\delta_1$, value of $\epsilon_1$ is polylinear map of values of $C_i$. Therefore, $\epsilon_1$ has maximal value when each $C_i$ has maximal value. Since we consider the estimation of the value of $\epsilon_1$ from the right, whether $C_i$ acquire maximal values at the same time is of no concern to us. Since the range of the map $f_i$ is compact set, then, according to the theorem 2.3.11, the following value is defined

$$F_i = \sup \|f_i(x)\| \tag{3.2.21}$$

From the equation (3.2.21) and the statement 2.1.9.1, it follows that

$$F_i \geq 0 \tag{3.2.22}$$

From equations (3.2.20), (3.2.21), it follows that we can put

$$C_i = F_i + \delta_1 \tag{3.2.23}$$

From the equation (3.2.23) and statements (3.2.22), (3.2.16) it follows that

$$C_i > 0 \quad \frac{dC_i}{d\delta_1} > 0 \tag{3.2.24}$$

From equations (3.2.19), (3.2.23) and the statement (3.2.24), it follows that $\epsilon_1$ is polynomial strictly monotone increasing function of $\delta_1$ such that

$$\delta_1 = 0 \Rightarrow \epsilon_1 = 0$$

According to the theorem 2.3.5, the map (3.2.19) maps the interval $[0, \delta_1)$ into the interval $[0, \epsilon_1)$. According to the theorem 2.3.3, for given $\epsilon > 0$ there exist $\delta > 0$ such that

$$\epsilon_1(\delta) < \epsilon$$

According to construction, a value of $N$ depends on a value of $\delta_1$. We choose the value of $N$ corresponding to $\delta_1 = \delta$. Therefore, for given $\epsilon \in R$, $\epsilon > 0$, there exists $N$ such that the condition $n > N$ implies that

$$f(g_{1\cdot n}(x))(g_{2\cdot n}(x)) \in B_o(f(g_1(x))(g_2(x)), \epsilon) \tag{3.2.25}$$

The equation

$$f(g_1(x))(g_2(x)) = \lim_{n \to \infty} f(g_{1\cdot n}(x))(g_{2\cdot n}(x))$$

follows from the equation (3.2.25) and the theorem 2.6.4. □

### 3.3. Completion of Representation of $\Omega$-Group

DEFINITION 3.3.1. Let

$$f: A_1 \dashrightarrow A_2$$

be representation of normed $\Omega_1$-group $A_1$ with norm $\|x\|_1$ in normed $\Omega_2$-group $A_2$ with norm $\|x\|_2$. The representation

$$g: B_1 \dashrightarrow B_2$$

of complete $\Omega_1$-group $B_1$ in complete $\Omega_2$-group $B_2$ is called **completion of representation** $f$, if

3.3.1.1: $\Omega$-group $B_1$ is completion of normed $\Omega_1$-group $A_1$.
3.3.1.2: $\Omega$-group $B_2$ is completion of normed $\Omega_2$-group $A_2$.
3.3.1.3: Since $a_1 \in A_1$, $a_2 \in A_2$, then

$$g(a_1)(a_2) = f(a_1)(a_2) \tag{3.3.1}$$

□

Theorem 3.3.2. *Let*
$$f : A_1 \dashrightarrow A_2$$
*be representation of normed $\Omega_1$-group $A_1$ with norm $\|x\|_1$ in normed $\Omega_2$-group $A_2$ with norm $\|x\|_2$. Let representations*
$$g_1 : B_{11} \dashrightarrow B_{12}$$
$$g_2 : B_{21} \dashrightarrow B_{22}$$
*be completions of representation $f$.*

3.3.2.1: *For $i = 1, 2$, there exists isomorphism of $\Omega$-group*

(3.3.2) $$h_i : B_{i1} \to B_{i2}$$

*such that*

(3.3.3) $$h_i(a_i) = a_i \quad a_i \in A_i$$

(3.3.4) $$\|h_i(b_i)\|_i = \|b_i\|_i$$

3.3.2.2: *Tuple of maps*

(3.3.5) $$(h_1 : B_{11} \to B_{12}, h_2 : B_{21} \to B_{22})$$

*is morphism of representations from $g_1$ into $g_2$.*

Proof. The statement 3.3.2.1 follows from the theorem 2.5.2.

Let $b_1 \in B_{11}$, $b_2 \in B_{12}$. According to the theorem 2.2.14, there exist sequences of $A_i$-numbers $a_{i \cdot n}$ such that

(3.3.6) $$b_1 = \lim_{n \to \infty} a_{1 \cdot n} \quad b_2 = \lim_{n \to \infty} a_{2 \cdot n}$$

From the equation (3.3.6), it follows that the sequence of $A_i$-numbers $a_{i \cdot n}$ is fundamental sequence in $\Omega$-group $A_i$. For any $n$, the following equation is evident

(3.3.7) $$f(a_{1 \cdot n})(a_{2 \cdot n}) = f(a_{1 \cdot n})(a_{2 \cdot n})$$

Since $f(a_{1 \cdot n})(a_{2 \cdot n}) \in A_2$, then the equation

(3.3.8) $$f(a_{1 \cdot n})(a_{2 \cdot n}) = h_2(g_1(a_{1 \cdot n})(a_{2 \cdot n})) = g_2(h_1(a_{1 \cdot n}))(h_2(a_{2 \cdot n}))$$

follows from the equation (3.3.3). The continuity of maps $h_1$, $h_2$ follows from the equation (3.3.4). Therefore, the equation

(3.3.9) $$h_2(g_1(b_1)(b_2)) = g_2(h_1(b_1))(h_2(b_2))$$

follows from equations (3.3.6), (3.3.8) and the theorem 3.1.5. The theorem follows from equations (3.3.9), [7]-(2.2.4) and the definition [7]-2.2.2.  □

Lemma 3.3.3. *Let*
$$f : A_1 \dashrightarrow A_2$$
*be representation of normed $\Omega_1$-group $A_1$ with norm $\|x\|_1$ in normed $\Omega_2$-group $A_2$ with norm $\|x\|_2$. Let $a_{1 \cdot n}$ be a fundamental sequence of $A_1$-numbers. Let $a_{2 \cdot n}$ be a fundamental sequence of $A_2$-numbers. Then a sequence $f(a_{1 \cdot n})(a_{2 \cdot n})$ is a fundamental sequence of $A_2$-numbers.*

PROOF. Let

(3.3.10) $$\delta_1 \in R, \ \delta_1 > 0$$

Since the sequence $a_{1 \cdot n}$ is a fundamental sequence of $A_1$-numbers, then according to the theorem 2.1.20, it follows that for given $\delta_1$ there exists positive integer $N_1$ depending on $\delta_1$ and such, that

(3.3.11) $$a_{1 \cdot q} \in B_o(a_{1 \cdot p}, \delta_1)$$

for every $p, q > N_1$. Since the sequence $a_{2 \cdot n}$ is a fundamental sequence of $A_2$-numbers, then according to the theorem 2.1.20, it follows that for given $\delta_1$ there exists positive integer $N_2$ depending on $\delta_1$ and such, that

(3.3.12) $$a_{2 \cdot q} \in B_o(a_{2 \cdot p}, \delta_1)$$

for every $p, q > N_2$. Let

$$N = \max(N_1, N_2)$$

From statements (3.3.11), (3.3.12) and the theorem 3.1.4, it follows that

(3.3.13) $$f(a_{1 \cdot q})(a_{2 \cdot q}) \in B_o(f(a_{1 \cdot p})(a_{2 \cdot p}), \epsilon_1)$$

(3.3.14) $$\epsilon_1 = \|f\|(\delta_1 C_2 + \delta_1 C_1)$$

where

(3.3.15) $$C_i = \|a_{i \cdot p}\|_i + \delta_1$$

From the equation (3.3.15) and statements 2.1.9.1, (3.3.10) it follows that

(3.3.16) $$C_i > 0 \quad \frac{dC_i}{d\delta_1} > 0$$

From equations (3.3.14), (3.3.15) and the statement (3.3.16), it follows that $\epsilon_1$ is polynomial strictly monotone increasing function of $\delta_1$ such that

$$\delta_1 = 0 \Rightarrow \epsilon_1 = 0$$

According to the theorem 2.3.5, the map (3.3.14) maps the interval $[0, \delta_1)$ into the interval $[0, \epsilon_1)$. According to the theorem 2.3.3, for given $\epsilon > 0$ there exist $\delta > 0$ such that

$$\epsilon_1(\delta) < \epsilon$$

According to construction, a value of $N$ depends on a value of $\delta_1$. We choose the value of $N$ corresponding to $\delta_1 = \delta$. Therefore, for given $\epsilon \in R, \epsilon > 0$, there exists $N$ such that the condition $n > N$ implies that

(3.3.17) $$f(a_{1 \cdot q})(a_{2 \cdot q}) \in B_o(f(a_{1 \cdot p})(a_{2 \cdot p}), \epsilon)$$

From the statement (3.3.17) and the theorem 2.1.20, it follows that the sequence $f(a_{1 \cdot n})(a_{2 \cdot n})$ is a fundamental sequence. $\square$

THEOREM 3.3.4. *There exists completion of representation of normed $\Omega$-group.*

PROOF. Let $A_1$ be normed $\Omega_1$-group. Let $A_2$ be normed $\Omega_2$-group. According to the theorem 2.5.17, there exist completion of $B_1$ normed $\Omega_1$-group $A_1$ and completion of $B_2$ normed $\Omega_2$-group $A_2$.

To prove existence of completion $g$ of the representation $f$

$$g : B_1 \relbar\joinrel\twoheadrightarrow B_2$$

we will use notation introduced in the section 2.5. Consider equivalence relation on the set of fundamental sequences of $A_i$-numbers

(3.3.18) $$a_{i\cdot n} \sim b_{i\cdot n} \Leftrightarrow \lim_{n\to\infty}(a_{i\cdot n} - b_{i\cdot n}) = 0$$

Let $B_i$ be the set of classes of equivalent fundamental sequences of $A_i$-numbers. We use notation

(3.3.19) $$[a_{i\cdot n}] = \{b_{i\cdot n} : a_{i\cdot n} \sim b_{i\cdot n}\}$$

for class of sequences which are equivalent to the sequence $a_{i\cdot n}$.

LEMMA 3.3.5. *Let $a_{1\cdot n}$, $b_{1\cdot n}$ be fundamental sequences of $A_1$-numbers,*

(3.3.20) $$a_{1\cdot n} \sim b_{1\cdot n}$$

*Let $a_{2\cdot n}$, $b_{2\cdot n}$ be fundamental sequences of $A_2$-numbers,*

(3.3.21) $$a_{2\cdot n} \sim b_{2\cdot n}$$

*Then*

(3.3.22) $$f(a_{1\cdot n})(a_{2\cdot n}) \sim f(b_{1\cdot n})(b_{2\cdot n})$$

PROOF. Let

(3.3.23) $$\delta_1 \in R,\ \delta_1 > 0$$

From the statement (3.3.20) and the lemma 2.5.4, it follows that for given $\delta_1$ there exists positive integer $N_1$ depending on $\delta_1$ and such, that

(3.3.24) $$b_{1\cdot n} \in B_o(a_{1\cdot n}, \delta_1)$$

for every $n > N_1$. From the statement (3.3.21) and the theorem 2.5.4, it follows that for given $\delta_1$ there exists positive integer $N_2$ depending on $\delta_1$ and such, that

(3.3.25) $$b_{2\cdot n} \in B_o(a_{2\cdot n}, \delta_1)$$

for every $n > N_2$. Let

$$N = \max(N_1, N_2)$$

From statements (3.3.24), (3.3.25) and the theorem 3.1.4, it follows that

(3.3.26) $$f(b_{1\cdot n})(b_{2\cdot n}) \in B_o(f(a_{1\cdot n})(a_{2\cdot n}), \epsilon_1)$$

(3.3.27) $$\epsilon_1 = \|f\|(\delta_1 C_2 + \delta_1 C_1)$$

where

(3.3.28) $$C_i = \|a_{i\cdot n}\|_i + \delta_1$$

From the equation (3.3.28) and statements 2.1.9.1, (3.3.23) it follows that

(3.3.29) $$C_i > 0 \quad \frac{dC_i}{d\delta_1} > 0$$

From equations (3.3.27), (3.3.28) and the statement (3.3.29), it follows that $\epsilon_1$ is strictly monotone increasing function of $\delta_1$ such that

$$\delta_1 = 0 \Rightarrow \epsilon_1 = 0$$

According to the theorem 2.3.5, the map (3.3.27) maps the interval $[0, \delta_1)$ into the interval $[0, \epsilon_1)$. According to the theorem 2.3.3, for given $\epsilon > 0$ there exist $\delta > 0$ such that

$$\epsilon_1(\delta) < \epsilon$$

According to construction, a value of $N$ depends on a value of $\delta_1$. We choose the value of $N$ corresponding to $\delta_1 = \delta$. Therefore, for given $\epsilon \in R$, $\epsilon > 0$, there exists $N$ such that the condition $n > N$ implies that

$$(3.3.30) \qquad f(b_{1 \cdot n})(b_{2 \cdot n}) \in B_o(f(a_{1 \cdot n})(a_{2 \cdot n}), \epsilon)$$

The statement (3.3.22) follows from the statement (3.3.30) and the lemma 2.5.4. ⊙

According to lemmas 3.3.3, 3.3.5, the map

$$(3.3.31) \qquad g([a_{1 \cdot n}])([a_{2 \cdot n}]) = [f(a_{1 \cdot n})(a_{2 \cdot n})]$$

is defined properly.

Let $a_1 \in A_1$, $a_2 \in A_2$. According to the theorem 2.2.13, there exists the sequence of $A_i$-numbers $a_{i \cdot n}$, $n = 1, ...,$ such that $a_i = [a_{i \cdot n}]$. According to the theorem 3.1.5

$$(3.3.32) \qquad f(a_1)(a_2) = \lim_{n \to \infty} f(a_{1 \cdot n})(a_{2 \cdot n})$$

From equations (3.3.31), (3.3.32), it follows that the map $g$ satisfies to the statement 3.3.1.3

$$g(a_1)(a_2) = g([a_{1 \cdot n}])([a_{2 \cdot n}]) = f(a_1)(a_2)$$

LEMMA 3.3.6. *The map $g$ is homomorphism of $\Omega_1$-group $B_1$.*

PROOF. Let $a_{1 \cdot 1}, ..., a_{1 \cdot n} \in B_1$. According to the theorem 2.2.14, there exists the sequence of $A_1$-numbers $a_{1 \cdot i \cdot m}$, $m = 1, ...,$ such that

$$(3.3.33) \qquad a_{1 \cdot i} = [a_{1 \cdot i \cdot m}] = \lim_{m \to \infty} a_{1 \cdot i \cdot m}$$

Let $a_2 \in B_2$. According to the theorem 2.2.14, there exists the sequence of $A_2$-numbers $a_{2 \cdot m}$, $m = 1, ...,$ such that

$$(3.3.34) \qquad a_2 = [a_{2 \cdot m}] = \lim_{m \to \infty} a_{2 \cdot m}$$

According to the definition [7]-2.1.4, the map $f$ is homomorphism of $\Omega_1$-group $A_1$. Therefore, for $n$-ari operation $\omega \in \Omega_1$ and any $m$ following equation is true

$$(3.3.35) \qquad f(a_{1 \cdot 1 \cdot m} ... a_{1 \cdot n \cdot m} \omega) = f(a_{1 \cdot 1 \cdot m}) ... f(a_{1 \cdot n \cdot m}) \omega$$

From equations (2.6.1), (3.3.35), it follows that

$$(3.3.36) \qquad \begin{aligned} f(a_{1 \cdot 1 \cdot m} ... a_{1 \cdot n \cdot m} \omega)(a_{2 \cdot m}) &= (f(a_{1 \cdot 1 \cdot m}) ... f(a_{1 \cdot n \cdot m}) \omega)(a_{2 \cdot m}) \\ &= (f(a_{1 \cdot 1 \cdot m})(a_{2 \cdot m})) ... (f(a_{1 \cdot n \cdot m})(a_{2 \cdot m})) \omega \end{aligned}$$

The equation

$$(3.3.37) \qquad \lim_{m \to \infty} a_{1 \cdot 1 \cdot m} ... a_{1 \cdot n \cdot m} \omega = a_{1 \cdot 1} ... a_{1 \cdot n} \omega$$

follows from the equation (3.3.33) and the theorem 2.4.6. The equation

$$(3.3.38) \qquad \lim_{m \to \infty} f(a_{1 \cdot 1 \cdot m} ... a_{1 \cdot n \cdot m} \omega)(a_{2 \cdot m}) = g(a_{1 \cdot 1} ... a_{1 \cdot n} \omega)(a_2)$$

follows from equations (3.3.31), (3.3.34), (3.3.37) and the theorem 3.1.5. The equation

$$(3.3.39) \qquad \lim_{m \to \infty} f(a_{1 \cdot i \cdot m})(a_{2 \cdot m}) = g(a_{1 \cdot i})(a_2)$$

follows from equations (3.3.31), (3.3.33), (3.3.34) and the theorem 3.1.5. The equation

$$
\begin{aligned}
(3.3.40) \quad &\lim_{m\to\infty} (f(a_{1\cdot 1\cdot m})(a_{2\cdot m}))...(f(a_{1\cdot n\cdot m})(a_{2\cdot m}))\omega \\
&= (g(a_{1\cdot 1})(a_2))...(g(a_{1\cdot n})(a_2))\omega
\end{aligned}
$$

follows from the equation (3.3.39) and the theorem 2.4.6. The equation

$$(3.3.41) \quad g(a_{1\cdot 1}...a_{1\cdot n}\omega)(a_2) = (g(a_{1\cdot 1})(a_2))...(g(a_{1\cdot n})(a_2))\omega$$

follows from equations (3.3.36), (3.3.38), (3.3.40). From the equation (3.3.41) it follows that the map $g$ is homomorphism of $\Omega_1$-group $B_1$. ⊙

LEMMA 3.3.7. *For any $a_1 \in B_1$, a map $g(a_1)$ is homomorphism of $\Omega_2$-group $B_2$.*

PROOF. Let $a_{2\cdot 1}, ..., a_{2\cdot n} \in B_2$. According to the theorem 2.2.14, there exists the sequence of $A_2$-numbers $a_{2\cdot i\cdot m}$, $m = 1, ...,$ such that

$$(3.3.42) \quad a_{2\cdot i} = [a_{2\cdot i\cdot m}] = \lim_{m\to\infty} a_{2\cdot i\cdot m}$$

Let $a_1 \in B_1$. According to the theorem 2.2.14, there exists the sequence of $A_1$-numbers $a_{1\cdot m}$, $m = 1, ...,$ such that

$$(3.3.43) \quad a_1 = [a_{1\cdot m}] = \lim_{m\to\infty} a_{1\cdot m}$$

According to the definition [7]-2.1.4, the map $f(a_{1\cdot m})$, $a_{1\cdot m} \in A_1$, is homomorphism of $\Omega_2$-group $A_2$. For $n$-ari operation $\omega \in \Omega_2$ and any $m$ following equation is true

$$(3.3.44) \quad f(a_{1\cdot m})(a_{2\cdot 1\cdot m}...a_{2\cdot n\cdot m}\omega) = (f(a_{1\cdot m})(a_{2\cdot 1\cdot m}))...(f(a_{1\cdot m})(a_{2\cdot n\cdot m}))\omega$$

The equation

$$(3.3.45) \quad \lim_{m\to\infty} a_{2\cdot 1\cdot m}...a_{2\cdot n\cdot m}\omega = a_{2\cdot 1}...a_{2\cdot n}\omega$$

follows from the equation (3.3.42) and the theorem 2.4.6. The equation

$$(3.3.46) \quad \lim_{m\to\infty} f(a_{1\cdot m})(a_{2\cdot 1\cdot m}...a_{2\cdot n\cdot m}\omega) = g(a_1)(a_{2\cdot 1}...a_{2\cdot n}\omega)$$

follows from equations (3.3.31), (3.3.43), (3.3.45) and the theorem 3.1.5. The equation

$$(3.3.47) \quad \lim_{m\to\infty} f(a_{1\cdot m})(a_{2\cdot i\cdot m}) = g(a_1)(a_{2\cdot i})$$

follows from equations (3.3.31), (3.3.42), (3.3.43) and the theorem 3.1.5. The equation

$$
\begin{aligned}
(3.3.48) \quad &\lim_{m\to\infty} (f(a_{1\cdot m})(a_{2\cdot 1\cdot m}))...(f(a_{1\cdot m})(a_{2\cdot n\cdot m}))\omega \\
&= (g(a_1)(a_{2\cdot 1}))...(g(a_1)(a_{2\cdot n}))\omega
\end{aligned}
$$

follows from the equation (3.3.47) and the theorem 2.4.6. The equation

$$(3.3.49) \quad g(a_1)(a_{2\cdot 1}...a_{2\cdot n}\omega) = (g(a_1)(a_{2\cdot 1}))...(g(a_1)(a_{2\cdot n}))\omega$$

follows from equations (3.3.44), (3.3.46), (3.3.48). From the equation (3.3.49), it follows that, for any $a_1 \in B_1$, a map $g(a_1)$ is homomorphism of $\Omega_2$-group $B_2$. ⊙

From lemmas 3.3.6, 3.3.7 and the definition [7]-2.1.4, it follows that the map $g$ defined by the equation (3.3.31) is representation $g$ of complete $\Omega_1$-group $B_1$ in complete $\Omega_2$-group $B_2$. □

## 3.4. Ω-Ring

Let
$$f : A_1 \dashrightarrow A_2$$
be representation of $\Omega_1$-group $A_1$ with norm $\|x\|_1$ in $\Omega_2$-group $A_2$ with norm $\|x\|_2$. According to the definition [7]-2.1.4, the map $f$ is homomorphism of $\Omega_1$-group $A_1$. Therefore, if product is defined in $\Omega_1$-group $A_1$, then the following equation is true

(3.4.1) $$f(a_{1\cdot 1} a_{1\cdot 2}) = f(a_{1\cdot 1}) f(a_{1\cdot 2})$$

However, if product is not defined in $\Omega_1$-group $A_1$, then we can introduce the product in $\Omega_1$-group $A_1$ using the equation (3.4.1). This is a reason to believe that the product is defined in $\Omega_1$-group $A_1$.

Since the map $f$ is homomorphism of $\Omega_1$-group $A_1$, then the following equations are true

(3.4.2) $$\begin{aligned} f(a_1(b_{1\cdot 1} + b_{1\cdot 2}))(a_2) &= f(a_1)(f(b_{1\cdot 1} + b_{1\cdot 2})(a_2)) \\ &= f(a_1)(f(b_{1\cdot 1})(a_2) + f(b_{1\cdot 2})(a_2)) \\ &= f(a_1)(f(b_{1\cdot 1})(a_2)) + f(a_1)(f(b_{1\cdot 2})(a_2)) \\ &= f(a_1 b_{1\cdot 1})(a_2) + f(a_1 b_{1\cdot 2})(a_2) \\ &= f(a_1 b_{1\cdot 1} + a_1 b_{1\cdot 2})(a_2) \end{aligned}$$

(3.4.3) $$\begin{aligned} f((a_{1\cdot 1} + a_{1\cdot 2}) b_1)(a_2) &= f(a_{1\cdot 1} + a_{1\cdot 2})(f(b_1)(a_2)) \\ &= f(a_{1\cdot 1})(f(b_1)(a_2)) + f(a_{1\cdot 2})(f(b_1)(a_2)) \\ &= f(a_{1\cdot 1} b_1)(a_2) + f(a_{1\cdot 2} f(b_1)(a_2) \\ &= f(a_{1\cdot 1} b_1 + a_{1\cdot 2} b_1)(a_2) \end{aligned}$$

Since the representation $f$ is effective representation, then the equations

(3.4.4) $$a_1(b_{1\cdot 1} + b_{1\cdot 2}) = a_1 b_{1\cdot 1} + a_1 b_{1\cdot 2}$$

(3.4.5) $$(a_{1\cdot 1} + a_{1\cdot 2}) b_1 = a_{1\cdot 1} b_1 + a_{1\cdot 2} b_1$$

follow from the equations (3.4.2), (3.4.3). Therefore, the product is distributive over addition. Thus, $\Omega_1$-group $A_1$ is ring with respect to addition and product.

DEFINITION 3.4.1. $\Omega$-group in which product is defined is called **$\Omega$-ring**. □

REMARK 3.4.2. Biring (the definition [6]-2.6) is an example of $\Omega$-group in which two products are defined. This is associated with two different representations of algebra of matices with entries from noncommutative algebra. Accordingly, two structures of $\Omega$-ring are defined in the biring.

We observe different picture in octonion algebra. Since the product in octonion algebra is nonassociative, then considered representation of octonion algebra determines different product. □

DEFINITION 3.4.3. Effective representation of $\Omega$-ring $A$ in Abelian group is called $A$-**module**.[3.2] □

I consider a representation in Abelian group, not in $\Omega$-group. The reason is that arbitrary operations in universal algebra make condition of linear dependency more complicated. At the same time, in universal algebra generating representation, only addition and product are important for us.

---

[3.2] See also the definition [9]-2.1.2.

# CHAPTER 4

# References

[1] S. Burris, H.P. Sankappanavar, A Course in Universal Algebra, Springer-Verlag (March, 1982),
eprint <http://www.math.uwaterloo.ca/ snburris/htdocs/ualg.html>
(The Millennium Edition)

[2] G. E. Shilov, Calculus, Single Variable Functions, Parts 1 - 2, Moscow, Nauka, 1969

[3] A. N. Kolmogorov and S. V. Fomin. Introductory Real Analysis.
Translated and edited by Richard A. Silverman.
Dover Publication, 1975, ISBN-13: 978-0486612263

[4] W.A. Coppel, Number Theory: An Introduction to Mathematics,
Springer, 2009

[5] Michael J. Field, Differential Calculus and Its Applications,
Dover Publications, 2012; ISBN-13: 978-0486497952

[6] Aleks Kleyn, Biring of Matrices,
eprint arXiv:math.OA/0612111 (2007)

[7] Aleks Kleyn, Representation of Universal Algebra,
eprint arXiv:0912.3315 (2009)

[8] Aleks Kleyn, Linear Maps of Free Algebra,
eprint arXiv:1003.1544 (2010)

[9] Aleks Kleyn, Free Algebra with Countable Basis,
eprint arXiv:1211.6965 (2012)

[10] V. I. Smirnov, A Course of Higher Mathematics, volume I.
Translated by D. E. Brown.
Translation, edited and additions made by I. N. Sneddon.
Pergamon Press, Addison-Wesley Publishing Company, 1964

[11] Paul M. Cohn, Universal Algebra, Springer, 1981

[12] Paul M. Cohn, Algebra, Volume 3, John Wiley & Sons, 1991

[13] N. Bourbaki, General Topology, Chapters 1 - 4, Springer, 1989

[14] N. Bourbaki, General Topology, Chapters 5 - 10, Springer, 1989

[15] V. I. Arnautov, S. T. Glavatsky, A. V. Mikhalev,
Introduction to the theory of topological rings and modules, Volume 1995,
Marcel Dekker, Inc, 1996

# CHAPTER 5

# Index

# CHAPTER 6

# Special Symbols and Notations

# Нормированная Ω-группа

Александр Клейн

Aleks_Kleyn@MailAPS.org
`http://AleksKleyn.dyndns-home.com:4080/`
`http://sites.google.com/site/AleksKleyn/`
`http://arxiv.org/a/kleyn_a_1`
`http://AleksKleyn.blogspot.com/`

Аннотация. Если в $\Omega$-алгебре $A$ определена операция сложения, которая не обязательно коммутативна, то $\Omega$-алгебре $A$ называется $\Omega$-группой. Я также рассмотрел представление $\Omega$-группы. Определение нормы в $\Omega$-группе позволяет рассмотреть непрерывность операций и представления.

# Оглавление

Глава 1

# Предисловие

### 1.1. Предисловие

Впервые я обратил внимание на универсальную алгебру, которая является абелевой группой относительно сложения, в главе [8]-2. Позже я узнал, что подобная универсальная алгебра называется $\Omega$-группой.

$\Omega$-группа интересна для меня, так как я полагаю, что мы можем рассмотреть математический анализ в $\Omega$-группе. Так возникло решение изучить нормированные $\Omega$-группы.

Я уделил большое внимание топологии нормированной $\Omega$-группы.

### 1.2. Соглашения

Соглашение 1.2.1. *Пусть $A$ - $\Omega_1$-алгебра. Пусть $B$ - $\Omega_2$-алгебра. Запись*

$$A \longrightarrow_* \rightarrow B$$

*означает, что определено представление $\Omega_1$-алгебры $A$ в $\Omega_2$-алгебре $B$.* □

Без сомнения, у читателя могут быть вопросы, замечания, возражения. Я буду признателен любому отзыву.

Глава 2

# Нормированная $\Omega$-группа

## 2.1. Нормированная $\Omega$-группа

Пусть в $\Omega$-алгебре[2.1] $A$ определена операция сложения, которая не обязательно коммутативна.

Определение 2.1.1. Отображение
$$f : A \to A$$
называется **аддитивным отображением**, если
$$f(a+b) = f(a) + f(b)$$

$\square$

Определение 2.1.2. Отображение
$$f : A^n \to A$$
называется **полиаддитивным отображением**, если для любого $i$, $i = 1, ..., n$,
$$f(a_1, ..., a_i + b_i, ..., a_n) = f(a_1, ..., a_i, ..., a_n) + f(a_1, ..., b_i, ..., a_n)$$

$\square$

Определение 2.1.3. Если $\Omega$-алгебра $A$ является группой относительно операции сложения и любая операция $\omega \in \Omega$ является полиаддитивным отображением, то $\Omega$-алгебра $A$ называется $\Omega$-**группой**.[2.2] Если $\Omega$-группа $A$ является ассоциативной группой относительно операции сложения, то $\Omega$-алгебра $A$ называется **ассоциативной** $\Omega$-**группой**. Если $\Omega$-группа $A$ является абелевой группой относительно операции сложения, то $\Omega$-алгебра $A$ называется **абелевой** $\Omega$-**группой**. $\square$

Пример 2.1.4. Наиболее очевидный пример $\Omega$-группы является собственно группой. $\square$

Пример 2.1.5. Кольцо является $\Omega$-группой. $\square$

Пример 2.1.6. Бикольцо матриц над телом ([6]) является $\Omega$-группой. $\square$

Соглашение 2.1.7. *Мы будем полагать, что рассматриваемая $\Omega$-группа абелева.* $\square$

Элемент $\Omega$-группы $A$ называется $A$-**числом**.

---

[2.1]Смотри определение универсальной алгебры в [1, 11].

[2.2]Вы можете найти определение $\Omega$-группы по адресу

http://ncatlab.org/nlab/show/Omega-group

Определение 2.1.8. Пусть $A$ - $\Omega$-группа. Подгруппа $B$ аддитивной группы $A$, замкнутая относительно операций из $\Omega$, называется **подгруппой $\Omega$-группы** $A$. $\square$

Определение 2.1.9. **Норма на $\Omega$-группе** $A$[2.3] - это отображение

$$d \in A \to \|d\| \in R$$

такое, что

2.1.9.1: $\|a\| \geq 0$
2.1.9.2: $\|a\| = 0$ равносильно $a = 0$
2.1.9.3: $\|a + b\| \leq \|a\| + \|b\|$
2.1.9.4: $\| - a\| = \|a\|$

$\Omega$-группа $A$, наделённая структурой, определяемой заданием на $A$ нормы, называется **нормированной $\Omega$-группой**. $\square$

Замечание 2.1.10. Если $\Omega$-группа $B$ с нормой $\|b\|_B$ является подгруппой $\Omega$-группы $A$ с нормой $\|b\|_A$, то мы требуем $\|b\|_B = \|b\|_A$ $\square$

Теорема 2.1.11. *Пусть $A$ - нормированная $\Omega$-группа. Тогда*

$$(2.1.1) \qquad \|a - b\| \geq |\,\|a\| - \|b\|\,|$$

Доказательство. Так как $a = a - b + b$, то из утверждения 2.1.9.3 следует, что

$$(2.1.2) \qquad \|a\| \leq \|a - b\| + \|b\|$$

Неравенство

$$(2.1.3) \qquad \|a - b\| \geq \|a\| - \|b\|$$

следует из неравенства (2.1.2). Из утверждения 2.1.9.4 и неравенства (2.1.3) следует, что

$$(2.1.4) \qquad \|a - b\| = \|b - a\| \geq \|b\| - \|a\|$$

Неравенство (2.1.1) следует из неравенств (2.1.3), (2.1.4). $\square$

Определение 2.1.12. Пусть $A$ - нормированная $\Omega$-группа. Для $n$-арной операции $\omega$, величина

$$(2.1.5) \qquad \|\omega\| = \sup \frac{\|a_1...a_n\omega\|}{\|a_1\|...\|a_n\|}$$

называется **нормой операции** $\omega$. $\square$

Теорема 2.1.13. *Пусть $A$ - нормированная $\Omega$-группа. Для $n$-арной операции $\omega$,*

$$(2.1.6) \qquad \|a_1...a_n\omega\| \leq \|\omega\|\|a_1\|...\|a_n\|$$

Доказательство. Из равенства (2.1.5) следует, что

$$(2.1.7) \qquad \frac{\|a_1...a_n\omega\|}{\|a_1\|...\|a_n\|} \leq \sup \frac{\|a_1...a_n\omega\|}{\|a_1\|...\|a_n\|} = \|\omega\|$$

Неравенство (2.1.6) следует из неравенства (2.1.7). $\square$

---

[2.3]Определение дано согласно определению из [13], гл. IX, §3, п°2, а также согласно определению [14]-1.1.12, с. 23.

Определение 2.1.14. Пусть $A$ - нормированная Ω-группа. Пусть $a \in A$. Множество
$$B_o(a, R) = \{b \in A : \|b - a\| < R\}$$
называется **открытым шаром** с центром в $a$. □

Определение 2.1.15. Пусть $A$ - нормированная Ω-группа. Пусть $a \in A$. Множество
$$B_c(a, R) = \{b \in A : \|b - a\| \leq R\}$$
называется **замкнутым шаром** с центром в $a$. □

Теорема 2.1.16. *Из утверждения $b \in B_o(a, R)$ следует утверждение*[2.4] *$a \in B_o(b, R)$.*

Доказательство. Теорема следует из определения 2.1.14. □

Определение 2.1.17. Пусть $A$ - нормированная Ω-группа. Элемент $a \in A$ называется **пределом последовательности** $a_n$
$$a = \lim_{n \to \infty} a_n$$
если для любого $\epsilon \in R, \epsilon > 0$, существует, зависящее от $\epsilon$, натуральное число $n_0$ такое, что $\|a_n - a\| < \epsilon$ для любого $n > n_0$. Мы будем также говорить, что **последовательность** $a_n$ **сходится** к $a$. □

Теорема 2.1.18. *Пусть $A$ - нормированная Ω-группа. Элемент $a \in A$ является пределом последовательности $a_n$*
$$a = \lim_{n \to \infty} a_n$$
*если для любого $\epsilon \in R, \epsilon > 0$, существует, зависящее от $\epsilon$, натуральное число $n_0$ такое, что*
$$a_n \in B_o(a, \epsilon)$$
*для любого $n > n_0$.*

Доказательство. Следствие определений 2.1.14, 2.1.17. □

Определение 2.1.19. Пусть $A$ - нормированная Ω-группа. Последовательность $a_n$, $a_n \in A$ называется **фундаментальной** или **последовательностью Коши**, если для любого $\epsilon \in R, \epsilon > 0$, существует, зависящее от $\epsilon$, натуральное число $n_0$ такое, что $\|a_p - a_q\| < \epsilon$ для любых $p, q > n_0$. □

Теорема 2.1.20. *Пусть $A$ - нормированная Ω-группа. Последовательность $a_n$, $a_n \in A$ является фундаментальной последовательностью, если для любого $\epsilon \in R, \epsilon > 0$, существует, зависящее от $\epsilon$, натуральное число $n_0$ такое, что*
$$a_q \in B_o(a_p, \epsilon)$$
*для любых $p, q > n_0$.*

Доказательство. Следствие определений 2.1.14, 2.1.19. □

---

[2.4]Аналогичная теорема верна для замкнутых шаров.

Теорема 2.1.21. *Пусть $A$ - нормированная $\Omega$-группа. Пусть $a_n$, $n = 1$, ..., - фундаментальная последовательность. Пусть $b_n$, $n = 1$, ..., - последовательность. Пусть*

$$\lim_{n \to \infty} (a_n - b_n) = 0 \tag{2.1.8}$$

*Тогда $b_n$ - фундаментальная последовательность.*

Доказательство. Из равенства (2.1.8) и определения 2.1.17 следует, что для заданного $\epsilon \in R$, $\epsilon > 0$, существует, зависящее от $\epsilon$, натуральное число $N_1$ такое, что

$$\|a_n - b_n\| < \frac{\epsilon}{3} \tag{2.1.9}$$

для любого $n > N_1$. Согласно определению 2.1.19, для заданного $\epsilon \in R$, $\epsilon > 0$, существует, зависящее от $\epsilon$, натуральное число $N_2$ такое, что

$$\|a_p - a_q\| < \frac{\epsilon}{3} \tag{2.1.10}$$

для любых $p$, $q > N_2$. Пусть

$$N = \max(N_1, N_2)$$

Из неравенств (2.1.9), (2.1.10) следует, что для заданного $\epsilon \in R$, $\epsilon > 0$, существует, зависящее от $\epsilon$, натуральное число $N$ такое, что

$$\|b_p - b_q\| = \|b_p - a_p + a_p - a_q + a_q - b_q\| \leq \|b_p - a_p\| + \|a_p - a_q\| + \|a_q - b_q\| < \epsilon$$

для любых $p$, $q > N$. Согласно определению 2.1.19, последовательность $b_n$ - фундаментальная последовательность. $\square$

Теорема 2.1.22. *Пусть $A$ - нормированная $\Omega$-группа. Пусть $a_n$, $b_n$, $n = 1$, ..., - фундаментальные последовательности. Пусть*

$$\lim_{n \to \infty} (a_n - b_n) = 0 \tag{2.1.11}$$

*Если последовательность $a_n$ сходится, то последовательность $b_n$ сходится и*

$$\lim_{n \to \infty} a_n = \lim_{n \to \infty} b_n \tag{2.1.12}$$

Доказательство. Из равенства (2.1.11) и определения 2.1.17 следует, что для заданного $\epsilon \in R$, $\epsilon > 0$, существует, зависящее от $\epsilon$, натуральное число $N_1$ такое, что

$$\|a_n - b_n\| < \frac{\epsilon}{2} \tag{2.1.13}$$

для любого $n > N_1$. Согласно определению 2.1.23 и теореме 2.1.24, существует предел $a$ последовательности $a_n$. Согласно определению 2.1.17, для заданного $\epsilon \in R$, $\epsilon > 0$, существует, зависящее от $\epsilon$, натуральное число $N_2$ такое, что

$$\|a_n - a\| < \frac{\epsilon}{2} \tag{2.1.14}$$

для любого $n > N_2$. Пусть

$$N = \max(N_1, N_2)$$

Из неравенств (2.1.13), (2.1.14) и утверждения 2.1.9.3 следует, что для заданного $\epsilon \in R$, $\epsilon > 0$, существует, зависящее от $\epsilon$, натуральное число $N$ такое, что

$$\|a - b_n\| = \|a - a_n + a_n - b_n\| \leq \|a - a_n\| + \|a_n - b_n\| < \epsilon$$

для любого $n > N$. Согласно определению 2.1.17, последовательность $b_n$ сходится к $a$. □

Определение 2.1.23. Нормированная Ω-группа $A$ называется **полной** если любая фундаментальная последовательность элементов Ω-группы $A$ сходится, т. е. имеет предел в этой Ω-группе. □

Теорема 2.1.24. *Пусть $A$ - полная Ω-группа. Любая фундаментальная последовательность имеет один и только один предел.*

Доказательство. Пусть $a_n$, $n = 1, ...,$ - фундаментальная последовательность. Согласно определению 2.1.19, для заданного $\epsilon \in R$, $\epsilon > 0$, существует, зависящее от $\epsilon$, натуральное число $n_1$ такое, что

(2.1.15) $$\|a_p - a_q\| < \frac{\epsilon}{3}$$

для любых $p$, $q > n_1$. Пусть

(2.1.16) $$a = \lim_{n \to \infty} a_n$$

и

(2.1.17) $$b = \lim_{n \to \infty} a_n$$

Из равенства (2.1.16) и определения 2.1.17 следует, что для заданного $\epsilon \in R$, $\epsilon > 0$, существует, зависящее от $\epsilon$, натуральное число $n_2$ такое, что

(2.1.18) $$\|a_p - a\| < \frac{\epsilon}{3}$$

для любого $p > n_2$. Из равенства (2.1.17) и определения 2.1.17 следует, что для заданного $\epsilon \in R$, $\epsilon > 0$, существует, зависящее от $\epsilon$, натуральное число $n_3$ такое, что

(2.1.19) $$\|a_q - b\| < \frac{\epsilon}{3}$$

для любого $q > n_3$. Пусть

$$n_0 = \max(n_1, n_2, n_3)$$

Из неравенств (2.1.15), (2.1.18), (2.1.19) следует, что для заданного $\epsilon \in R$, $\epsilon > 0$, существует, зависящее от $\epsilon$, натуральное число $n_0$ такое, что

$$\|a - b\| = \|a - a_p + a_p - a_q + a_q - b\| \leq \|a - a_p\| + \|a_p - a_q\| + \|a_q - b\| < \epsilon$$

для любых $p$, $q > n_0$. Следовательно, $\|a - b\| = 0$. Согласно утверждению 2.1.9.2, $a = b$. □

## 2.2. Топология $\Omega$-группы

Замечание 2.2.1. Инвариантное расстояние на аддитивной группе $\Omega$-группы $A$

$$(2.2.1) \qquad d(a,b) = \|a-b\|$$

определяет топологию метрического пространства, согласующуюся со структурой $\Omega$-группы в $A$. Множество открытых шаров нормированной $\Omega$-группы $A$ является базой топологии,[2.5] согласованной с расстоянием (2.2.1). $\square$

Определение 2.2.2. Пусть $A$ - нормированная $\Omega$-группа. Множество $U \subset A$ называется **открытым**,[2.6] если для любого $A$-числа $a \in U$ существует $\epsilon \in R, \epsilon > 0$, такое, что $B_o(a, \epsilon) \subset U$. $\square$

Определение 2.2.3. Пусть $B$ - подмножество топологического пространства $A$. Точка $x \in A$ называется **точкой прикосновения множества** $B$,[2.7] если для любого $\epsilon \in R, \epsilon > 0$, открытый шар $B_o(x, \epsilon)$ содержит хотя бы одну точку множества $B$

$$B_o(x,\epsilon) \cap B \neq \emptyset$$

Множество всех точек прикосновения множества $B$ называется **замыканием множества** $B$. Замыкание множества $B$ обозначается $[B]$. $\square$

Замечание 2.2.4. Множество $B$ замкнуто тогда и только тогда, когда $B = [B]$.

Это утверждение может быть либо определением замкнутого множества, либо теоремой. Например, рассмотрим определение замкнутого множества как дополнение к открытому множеству. Если $x \in A \setminus [B]$. то существует окрестность точки $x$, не содержащая точек множества $B$. Следовательно, $A \setminus [B]$ открытое множество. $\square$

Определение 2.2.5. Пусть $A$ - нормированная $\Omega$-группа. Множество $B \subset A$ называется **плотным в множестве**[2.8] $C \subset A$, если $C \subset [B]$. Множество $B \subset A$ называется **всюду плотным**,[2.9] если $[B] = A$. $\square$

Определение 2.2.6 (**вторая аксиома счётности**). В топологическом пространстве со счётной базой существует по крайней мере одна база, состоящая не более чем из счётного числа множеств.[2.10] $\square$

Определение 2.2.7 (Первая аксиома отделимости). Топологическое пространство $A$ называется $T_1$-**пространством**,[2.11] если для точек $x, y \in A$,

---

[2.5] Смотри определение базы топологии в [3], страница 87.

[2.6] В топологии обычно сперва определяют открытое множество, а затем базу топологии. В случае метрического или нормированного пространства, удобнее дать определение открытого множества, опираясь на определение базы топологии. В этом случае определение основано на одном из свойств базы топологии. Непосредственная проверка позволяет убедиться, что определённое таким образом открытое множество удовлетворяет основным свойствам.

[2.7] Смотри определение точки прикосновения и замыкания множества в [3], страница 85.

[2.8] Смотри также определение в [3], страница 59.

[2.9] Мы также говорим, что множество $B$ всюду плотно в $\Omega$-группе $A$.

[2.10] Смотри также определению в [3], страница 88. Согласно теореме в [3], страницы 88, 89, в топологическом пространстве со второй аксиомой счётности существует счётное всюду плотное множество.

[2.11] Смотри аналогичное определение в [3], с. 94.

$x \neq y$, существуют окрестности $U_x$, $x \in U_x$, и $U_y$, $y \in U_y$, такие, что $y \notin U_x$, $x \notin U_y$. □

Теорема 2.2.8. *Любая точка $T_1$-пространства является замкнутым множеством.*[2.12]

Доказательство. Если $x \neq y$, то согласно определению 2.2.7, существует окрестность $O_y$ точки $y$ такая, что $x \notin O_y$. Согласно определению 2.2.3 $y \notin [\{x\}]$. Следовательно, $\{x\} = [\{x\}]$. Согласно замечанию 2.2.4, $\{x\}$ - замкнутое множество. □

Теорема 2.2.9. *Нормированная $\Omega$-группа $A$ является $T_1$-пространством тогда и только тогда, когда для точек $x, y \in A$, $x \neq y$, существуют открытые шары $B_o(x, r_x)$, $B_o(y, r_y)$ такие, что $y \notin B_o(x, r_x)$, $x \notin B_o(y, r_y)$.*

Доказательство. Согласно определению 2.2.7, существуют окрестности $U_x$, $x \in U_x$, и $U_y$, $y \in U_y$, такие, что $y \notin U_x$, $x \notin U_y$. Согласно замечанию 2.2.1 и определению на странице [3]-87, существуют открытые шары $B_o(x, r_x)$, $B_o(y, r_y)$ такие, что $B_o(x, r_x) \subset U_x$, $B_o(y, r_y) \subset U_y$. Следовательно, $y \notin B_o(x, r_x)$, $x \notin B_o(y, r_y)$. □

Соглашение 2.2.10. *Для того, чтобы топология нормированной $\Omega$-группы $A$ была нетривиальной, мы потребуем следующее.*

2.2.10.1: *$\Omega$-группа $A$ является $T_1$-пространством.*
2.2.10.2: *Нормированная $\Omega$-группа $A$ удовлетворяет второй аксиоме счётности.*
2.2.10.3: *Множество $\{x\}$ не является открытым множеством.*

□

Теорема 2.2.11. *Пусть $A$ - нормированная $\Omega$-группа. Для любых $x \in A$, $r \in R$, существует $y \in B_o(x, r)$, $y \neq x$.*

Доказательство. Если мы предположим, что не существует $y \in B_o(x, r)$, $y \neq x$, то в этом случае $B_o(x, r) = \{x\}$ и множество $\{x\}$ является открытым. Это утверждение противоречит соглашению 2.2.10.3. Следовательно, предположение не верно. □

Теорема 2.2.12. *Пусть $A$ - нормированная $\Omega$-группа. Для любых $a, b \in A$, $b \neq a$, существует $c \in A, c \neq a$, такой, что*

$$(2.2.2) \qquad \|c - a\| < \|b - a\|$$

Доказательство. Согласно теореме 2.2.9, существует открытый шар $B_o(a, r)$ такой, что $b \notin B_o(a, r)$. Согласно определению 2.1.14

$$(2.2.3) \qquad r \leq \|b - a\|$$

Согласно теореме 2.2.11, существует $c \in B_o(a, r), c \neq a$. Согласно определению 2.1.14

$$(2.2.4) \qquad \|c - a\| < r$$

Неравенство (2.2.2) следует из неравенств (2.2.3), (2.2.4). □

---
[2.12]Смотри также замечание после определения в [3], с. 94.

Теорема 2.2.13. *Пусть $A$ - нормированная $\Omega$-группа. Для любого $a \in A$ существует последовательность $A$-чисел $a_n$, $a_n \neq a$, такая, что*

$$\lim_{n \to \infty} a_n = a \tag{2.2.5}$$

Доказательство. Существуют разные способы построить последовательность $a_n$.

Согласно теореме 2.2.11 для любого $n$ существует $a_n \in B(a, \frac{1}{n})$, $a_n \neq a$. Согласно теореме 2.1.18, последовательность $a_n$ сходится к $a$.

В рассматриваемой последовательности возможно, что $a_n = a_{n+1}$. Однако очевидно, что если $m > \frac{1}{\|a_n\|}$, то $a_m \neq a_n$. Мы можем построить последовательность, в которой все элементы различны.

Пусть $a_1 \in A$. Предположим, мы выбрали $a_n \in A$, $n \geq 1$. Пусть

$$r = \frac{\|a_n - a\|}{2} \tag{2.2.6}$$

Согласно теореме 2.2.11 существует $a_{n+1} \in B(a, r)$, $a_{n+1} \neq a$. Согласно равенству (2.2.6) и определению 2.1.14, $a_{n+1} \neq a_n$. Согласно теореме 2.1.18, последовательность $a_n$ сходится к $a$. □

Теорема 2.2.14. *Пусть $A$ - нормированная $\Omega$-группа. Пусть множество $B \subset A$ плотно в множестве $C \subset A$, . Тогда для любого $A$-числа $b \in C$ существует последовательность $A$-чисел $b_n$, $b_n \in B$, сходящаяся к $b$*

$$b = \lim_{n \to \infty} b_n \tag{2.2.7}$$

Доказательство. Согласно определениям 2.2.3, 2.2.5, для любого $n > 0$ существует $b_n$ такой, что

$$b_n \in B \cap B_o(b, 1/n)$$

Согласно теореме 2.1.18, последовательность $b_n$ сходится к $b$, так как $b_n \in B_o(x, \epsilon)$ для любого $n > 1/\epsilon$. □

Определение 2.2.15. Множество $T$ топологического пространства называется **компактным**, если любое его открытое покрытие содержит конечное подпокрытие.[2.13] □

Определение 2.2.16. Множество $T$ топологического пространства называется **связным**,[2.14] если существуют открытые множества $B$ и $C$ такие, что из условий

2.2.16.1: $T \subset B \cup C$
2.2.16.2: $T \cap B \neq \emptyset$
2.2.16.3: $T \cap C \neq \emptyset$

следует $T \cap B \cap C \neq \emptyset$ . □

Пример 2.2.17. Связное открытое множество в поле действительных чисел является открытым интервалом.[2.15] □

---

[2.13] Смотри также определение в [3], страница 98.

[2.14] Смотри также определения 1 и 2 в [12], страницы 169, 170.

[2.15] Смотри замечание 5 в [3], страницы 65, 66.

### 2.3. Непрерывное отображение Ω-группы

Определение 2.3.1. Отображение

$$f : A_1 \to A_2$$

нормированной $\Omega_1$-группы $A_1$ с нормой $\|x\|_1$ в нормированную $\Omega_2$-группу $A_2$ с нормой $\|y\|_2$ называется **непрерывным**, если для любого сколь угодно малого $\epsilon > 0$ существует такое $\delta > 0$, что

$$\|x' - x\|_1 < \delta$$

влечёт

$$\|f(x') - f(x)\|_2 < \epsilon$$

□

Теорема 2.3.2. *Отображение*

$$f : A_1 \to A_2$$

*нормированной $\Omega_1$-группы $A_1$ с нормой $\|x\|_1$ в нормированную $\Omega_2$-группу $A_2$ с нормой $\|y\|_2$ непрерывно тогда и только тогда, когда для любого сколь угодно малого $\epsilon > 0$ существует такое $\delta > 0$, что*

$$b \in B_o(a, \delta)$$

*влечёт*

$$f(b) \in B_o(f(a), \epsilon)$$

Доказательство. Теорема является следствием определений 2.1.14, 2.3.1. □

Теорема 2.3.3. *Отображение*

$$f : A_1 \to A_2$$

*нормированной $\Omega_1$-группы $A_1$ с нормой $\|x\|_1$ в нормированную $\Omega_2$-группу $A_2$ с нормой $\|y\|_2$ непрерывно тогда и только тогда, когда прообраз открытого множества является открытым множеством.*

Доказательство. Пусть $U \subset f(A_1) \subset A_2$ открытое множество и

(2.3.1) $$W = f^{-1}(U)$$

- Пусть отображение $f$ непрерывно. Пусть $a \in W$. Тогда $f(a) \in U$. Согласно определению 2.2.2, существует $\epsilon \in R, \epsilon > 0,$ такое, что

  (2.3.2) $$B_o(f(a), \epsilon) \subset U$$

  Согласно теореме 2.3.2, существует такое $\delta > 0$, что

  (2.3.3) $$b \in B_o(a, \delta)$$

  влечёт

  (2.3.4) $$f(b) \in B_o(f(a), \epsilon)$$

  Из равенств (2.3.2), (2.3.4) следует, что

  (2.3.5) $$f(b) \in U$$

  Из равенств (2.3.1), (2.3.5) следует, что

  (2.3.6) $$b \in W$$

Так как из утверждения (2.3.3) следует утверждение (2.3.6), то $B_o(a, \delta) \subset W$. Согласно определению 2.2.2 множество $W$ открыто.

- Пусть прообраз открытого множества является открытым множеством. Пусть $U = B_o(f(a), \epsilon)$. Так как $W$ - открытое множество, то согласно определению 2.2.2 существует $\delta > 0$ такое, что $B_o(a, \delta) \subset W$. Следовательно, из утверждения $b \in B_o(a, \delta)$ следует $f(b) \in B_o(f(a), \epsilon)$. Согласно теореме 2.3.2, отображение $f$ непрерывно.

$\square$

Теорема 2.3.4. *Пусть*
$$f : X \to Y$$
*непрерывное отображение топологического пространства $X$ в топологическое пространство $Y$. Пусть $T$ - связное множество топологического пространства $X$. Тогда $f(T)$ - связное множество топологического пространства $Y$.*[2.16]

Доказательство. Пусть $B$, $C$ - открытые множества топологического пространства $Y$. Так как отображение $f$ непрерывно, то согласно теореме 2.3.3, множества $f^{-1}(B)$, $f^{-1}(C)$ открыты. Если утверждения 2.2.16.1, 2.2.16.2, 2.2.16.3 верны для множеств $B$, $C$, $f(T)$, то утверждения 2.2.16.1, 2.2.16.2, 2.2.16.3 верны для множеств $f^{-1}(B)$, $f^{-1}(C)$, $T$. Если мы предположим, что множество $f(T)$ не связно, то согласно определению 2.2.16
$$(2.3.7) \quad B \cap C \cap f(T) = \emptyset$$
Из равенства (2.3.7) следует, что
$$(2.3.8) \quad f^{-1}(B) \cap f^{-1}(C) \cap T = \emptyset$$
Согласно определению 2.2.16, множество $T$ не связно. Это утверждение противоречит предположению теоремы. Следовательно, множество $f(T)$ связно. $\square$

Теорема 2.3.5. *Пусть*
$$f : R \to R$$
*непрерывное отображение поля действительных чисел. Тогда образ интервала является интервалом.*

Доказательство. Теорема является следствием теоремы 2.3.4 и примера 2.2.17. $\square$

Теорема 2.3.6. *Отображение*
$$f : A_1 \to A_2$$
*нормированной $\Omega_1$-группы $A_1$ с нормой $\|x\|_1$ в нормированную $\Omega_2$-группу $A_2$ с нормой $\|y\|_2$ непрерывно тогда и только тогда, когда из условия, что последовательность $A_1$-чисел $a_n$ сходится, следует, что последовательность $A_2$-чисел $f(a_n)$ сходится и верно равенство*
$$(2.3.9) \quad f(\lim_{n \to \infty} a_n) = \lim_{n \to \infty} f(a_n)$$

---

[2.16]Смотри также предложение 4 в [12], страница 171.

Доказательство. Пусть $A_1$-число $a$ является пределом последовательности $a_n$

$$a = \lim_{n \to \infty} a_n \tag{2.3.10}$$

- Пусть равенство (2.3.9) верно для $A$-числа $a$. Равенство

$$f(a) = \lim_{n \to \infty} f(a_n) \tag{2.3.11}$$

следует из равенств (2.3.9), (2.3.10). Из равенства (2.3.11) и теоремы 2.1.18, следует, что для заданного $\epsilon \in R, \epsilon > 0$, существует, зависящее от $\epsilon$, натуральное число $n_2$ такое, что

$$f(a_n) \in B_o(f(a), \epsilon) \tag{2.3.12}$$

для любого $n > n_2$. Из равенства (2.3.10) и теоремы 2.1.18, следует, что для заданного $\delta \in R, \delta > 0$, существует, зависящее от $\delta$, натуральное число $n_1$ такое, что

$$a_n \in B_o(a, \delta) \tag{2.3.13}$$

для любого $n > n_1$. Если $n_1 \le n_2$, то оценка $\delta$ завышена. Мы положим

$$\delta = \frac{\|a - a_{n_1+1}\|}{2}$$

и повторим оценку значения $n_1$. Очевидно, что новое значение $n_1$ будет больше предыдущего. Следовательно, за конечное число итераций мы найдём $\delta$ такое, что $n_1 > n_2$ и из утверждения (2.3.13) следует утверждение (2.3.12). Согласно теореме 2.3.2, отображение $f$ непрерывно.

- Пусть отображение $f$ непрерывно. Согласно теореме 2.3.2, для любого сколь угодно малого $\epsilon > 0$ существует такое $\delta > 0$, что

$$b \in B_o(a, \delta) \tag{2.3.14}$$

влечёт

$$f(b) \in B_o(f(a), \epsilon) \tag{2.3.15}$$

Из равенства (2.3.10) и теоремы 2.1.18, следует, что для заданного $\delta \in R, \delta > 0$, существует, зависящее от $\delta$, натуральное число $n_1$ такое, что

$$a_n \in B_o(a, \delta) \tag{2.3.16}$$

для любого $n > n_1$. Из равенств (2.3.14), (2.3.15), (2.3.16) следует, что для заданного $\epsilon \in R, \epsilon > 0$, существует, зависящее от $\epsilon$, натуральное число $n_1$ такое, что

$$f(a_n) \in B_o(f(a), \epsilon) \tag{2.3.17}$$

для любого $n > n_1$. Равенство

$$f(a) = \lim_{n \to \infty} f(a_n) \tag{2.3.18}$$

следует из равенства (2.3.17) и теоремы 2.1.18. Равенство (2.3.9) следует из равенств (2.3.18), (2.3.10).

□

Определение 2.3.7. Пусть
$$f:A_1\to A_2$$
отображение нормированной $\Omega_1$-группы $A_1$ с нормой $\|x\|_1$ в нормированную $\Omega_2$-группу $A_2$ с нормой $\|y\|_2$. Величина

(2.3.19) $$\|f\|=\sup\frac{\|f(x)\|_2}{\|x\|_1}$$

называется **нормой отображения** $f$.                                       □

Теорема 2.3.8. *Пусть*
$$f:A_1\to A_2$$
*аддитивное отображение нормированной $\Omega_1$-группы $A_1$ с нормой $\|x\|_1$ в нормированную $\Omega_2$-группу $A_2$ с нормой $\|y\|_2$. Отображение $f$ непрерывно, если $\|f\|<\infty$.*

Доказательство. Поскольку отображение $f$ аддитивно, то согласно определению 2.3.7
$$\|f(x)-f(y)\|_2=\|f(x-y)\|_2\le\|f\|\,\|x-y\|_1$$
Возьмём произвольное $\epsilon>0$. Положим $\delta=\dfrac{\epsilon}{\|f\|}$. Тогда из неравенства
$$\|x-y\|_1<\delta$$
следует
$$\|f(x)-f(y)\|_2\le\|f\|\,\delta=\epsilon$$
Согласно определению 2.3.1 отображение $f$ непрерывно.                □

Теорема 2.3.9. *Пусть $A$ - нормированная $\Omega$-группа. Норма, определённая в $\Omega$-группе $A$, является непрерывным отображением $\Omega$-группы $A$ в поле действительных чисел.*

Доказательство. Пусть $\epsilon\in R$, $\epsilon>0$, и $a\in A$. Согласно теореме 2.2.11 существует $b\in B_o(a,\epsilon)$, $b\ne a$. Согласно определению 2.1.14

(2.3.20) $$\|a-b\|<\epsilon$$

Из неравенств (2.3.20), (2.1.1) следует, что

(2.3.21) $$|\,\|a\|-\|b\|\,|\le\|a-b\|<\epsilon$$

Так как неравенство (2.3.21) является следствием неравенства (2.3.20), то норма непрерывна согласно определению 2.3.1.                                       □

Теорема 2.3.10. *Пусть*
$$f:A\to B$$
*непрерывное отображение нормированной $\Omega$-алгебры $A$ в нормированную $\Omega$-алгебру $B$. Если $C$ - компактное множество $\Omega$-алгебры $A$, то образ $f(C)$ является компактным множеством $\Omega$-алгебры $B$.*

Доказательство. Пусть $D$ - открытое покрытие множества $f(C)$. Любое множество $E \in D$ является открытым множеством. Согласно теореме 2.3.3, прообраз $f^{-1}(E)$ является открытым множеством. Пусть

$$D' = \{E' : E' = f^{-1}(E), E \in D\}$$

Утверждение

$$C = \bigcup_{E' \in D'} E'$$

следует из утверждения

$$f(C) = \bigcup_{E \in D} E$$

Следовательно, $D'$ является открытым покрытием множества $C$. Согласно определению 2.2.15, существует конечное подпокрытие $E'_1, ..., E'_n$. Следовательно, покрытие $D$ имеет конечное подпокрытие $E_1 = f(E'_1), ..., E_n = f(E'_n)$. Согласно определению 2.2.15, множество $f(C)$ компактно. □

Теорема 2.3.11. *Пусть $C$ - компактное множество нормированной Ω-группы $A$. Тогда норма $\|x\|$, $x \in C$, ограничена сверху и снизу.*

Доказательство. Согласно теореме 2.3.9, норма является непрерывным отображением Ω-группы $A$ в поле действительных чисел. Согласно теореме 2.3.10, множество

$$C' = \{\|x\| : x \in C\}$$

компактно в поле действительных чисел. Согласно утверждению [2]-3.96.e, страница 118, множество $C'$ ограничено и содержит свои точные границы. □

## 2.4. Непрерывность операций в Ω-группе

Из утверждения 2.1.9.3 и теоремы 2.3.8 следует, что сумма в нормированной Ω-группе непрерывна. Из определения 2.1.12 следует, что операции в нормированной Ω-группе непрерывны.[2.17] В этом разделе мы рассмотрим утверждения, связанные с непрерывностью операций.

Теорема 2.4.1. *В нормированной Ω-группе $A$, верно следующее неравенство*

$$(2.4.1) \qquad \|(c_1 + c_2) - (a_1 + a_2)\| \leq \|(c_1 - a_1)\| + \|(c_2 - a_2)\|$$

Доказательство. Из утверждения 2.1.9.3 следует, что

$$(2.4.2) \qquad \begin{aligned} \|(c_1 + c_2) - (a_1 + a_2)\| &= \|(c_1 - a_1) + (c_2 - a_2)\| \\ &\leq \|(c_1 - a_1)\| + \|(c_2 - a_2)\| \end{aligned}$$

Неравенство (2.4.1) следует из неравенства (2.4.2). □

Теорема 2.4.2. *Пусть $A$ - нормированная Ω-группа. Для $c_1$, $c_2 \in A$, пусть $c_1 \in B_c(a_1, R_1)$, $c_2 \in B_c(a_2, R_2)$. Тогда*

$$(2.4.3) \qquad c_1 + c_2 \in B_c(a_1 + a_2, R_1 + R_2)$$

---

[2.17]Сравни определение 2.1.12 нормы операции и определение [9]-3.4.8 нормы полилинейного отображения.

Доказательство. Согласно определению 2.1.14

$$\|c_1 - a_1\| \leq R_1 \tag{2.4.4}$$
$$\|c_2 - a_2\| \leq R_2$$

Из неравенств (2.4.1), (2.4.4) следует, что

$$\|(c_1 + c_2) - (a_1 + a_2)\| \leq \|(c_1 - a_1)\| + \|(c_2 - a_2)\| \leq R_1 + R_2 \tag{2.4.5}$$

Утверждение (2.4.3) следует из неравенства (2.4.5) и определения 2.1.14. □

Теорема 2.4.3. *Пусть $A$ - нормированная $\Omega$-группа. Пусть последовательность $A$-чисел $a_n$ сходится и*

$$\lim_{n \to \infty} a_n = a \tag{2.4.6}$$

*Пусть последовательность $A$-чисел $b_n$ сходится и*

$$\lim_{n \to \infty} b_n = b \tag{2.4.7}$$

*Тогда последовательность $A$-чисел $a_n + b_n$ сходится и*

$$\lim_{n \to \infty} a_n + b_n = a + b \tag{2.4.8}$$

Доказательство. Из равенства (2.4.6) и теоремы 2.1.18 следует, что для заданного $\epsilon \in R, \epsilon > 0$, существует $N_a$ такое, что из условия $n > N_a$ следует

$$a_n \in B_o(a, \epsilon/2) \tag{2.4.9}$$

Из равенства (2.4.7) и теоремы 2.1.18 следует, что для заданного $\epsilon \in R, \epsilon > 0$, существует $N_b$ такое, что из условия $n > N_b$ следует

$$b_n \in B_o(b, \epsilon/2) \tag{2.4.10}$$

Пусть

$$N = \max(N_a, N_b)$$

Из равенств (2.4.9), (2.4.10), теоремы 2.4.2 и условия $n > N$ следует, что

$$a_n + b_n \in B_o(a + b, \epsilon/2 + \epsilon/2) = B_o(a + b, \epsilon) \tag{2.4.11}$$

Равенство (2.4.8) следует из равенства (2.4.11) и теоремы 2.1.18. □

Теорема 2.4.4. *Пусть $A$ - нормированная $\Omega$-группа. Пусть $\omega \in \Omega$ - $n$-арная операция. Верно следующее неравенство*

$$\|c_1...c_n\omega - a_1...a_n\omega\| \leq \|\omega\|(\|c_1 - a_1\|C_2...C_n + ... + C_1...C_{n-1}\|c_n - a_n\|) \tag{2.4.12}$$

*где*

$$C_i = \max(\|a_i\|, \|c_i\|) \quad i = 1, ..., n \tag{2.4.13}$$

Доказательство. Согласно определениям 2.1.2, 2.1.3,

$$\begin{aligned}
c_1...c_n\omega - a_1...a_n\omega &= c_1c_2...c_n\omega - a_1c_2...c_n\omega \\
&\quad + a_1c_2...c_n\omega - a_1a_2...c_n\omega \\
&\quad ...... \\
&\quad + a_1...a_{n-1}c_n\omega - a_1...a_{n-1}a_n\omega \\
&= (c_1 - a_1)c_2...c_n\omega + a_1(c_2 - a_2)...c_n\omega \\
&\quad + ... + a_1...a_{n-1}(c_n - a_n)\omega
\end{aligned} \tag{2.4.14}$$

Из равенства (2.4.14) и утверждения 2.1.9.3 следует, что

$$(2.4.15) \quad \|c_1...c_n\omega - a_1...a_n\omega\| \leq \|(c_1 - a_1)c_2...c_n\omega\| + \|a_1(c_2 - a_2)...c_n\omega\|$$
$$+ ... + \|a_1...a_{n-1}(c_n - a_n)\omega\|$$

Из равенства (2.4.13) и определения 2.1.12 следует, что

$$\|(c_1 - a_1)c_2...c_n\omega\| \leq \|c_1 - a_1\|C_2...C_n$$
$$(2.4.16) \quad ......$$
$$\|a_1...a_{n-1}(c_n - a_n)\omega\| \leq C_1...C_{n-1}\|c_n - a_n\|$$

Неравенство (2.4.12) следует из неравенств (2.4.15), (2.4.16). □

ТЕОРЕМА 2.4.5. *Пусть $A$ - нормированная $\Omega$-группа. Пусть $\omega \in \Omega$ - $n$-арная операция. Для $c_1, ..., c_n \in A$, пусть $c_1 \in B_c(a_1, R_1), ..., c_n \in B_c(a_n, R_n)$. Тогда*

$$(2.4.17) \quad c_1...c_n\omega \in B_c(a_1...a_n\omega, R)$$

*где*

$$(2.4.18) \quad R = \|\omega\|(R_1 C_2...C_n + ... + C_1...C_{n-1}R_n)$$

$$(2.4.19) \quad C_i = \|a_i\| + R_i \quad i = 1, ..., n$$

ДОКАЗАТЕЛЬСТВО. Согласно определению 2.1.14

$$\|c_1 - a_1\| \leq R_1$$
$$(2.4.20) \quad ...$$
$$\|c_n - a_n\| \leq R_n$$

Неравенство

$$(2.4.21) \quad \|c_i\| = \|a_i + c_i - a_i\| \leq \|a_i\| + \|c_i - a_i\| \leq \|a_i\| + R_i$$

следует из неравенства (2.4.20) и утверждения 2.1.9.3. Неравенство

$$(2.4.22) \quad C_i \leq \max(\|a_i\|, \|a_i\| + R_i) = \|a_i\| + R_i \quad i = 1, ..., n$$

следует из неравенства (2.4.21) и равенства (2.4.13). Равенство (2.4.19) следует из неравенства (2.4.22). Утверждение (2.4.17) следует из неравенств (2.4.12), (2.4.20) и равенства (2.4.19). □

ТЕОРЕМА 2.4.6. *Пусть $A$ - нормированная $\Omega$-группа. Пусть $\omega \in \Omega$ - $n$-арная операция. Пусть последовательность $A$-чисел $a_{i \cdot m}$ сходится и*

$$(2.4.23) \quad \lim_{m \to \infty} a_{i \cdot m} = a_i$$

*Тогда последовательность $A$-чисел $a_{1 \cdot m}...a_{n \cdot m}\omega$ сходится и*

$$(2.4.24) \quad \lim_{m \to \infty} a_{1 \cdot m}...a_{n \cdot m}\omega = a_1...a_n\omega$$

ДОКАЗАТЕЛЬСТВО. Из равенства (2.4.23) и теоремы 2.1.18 следует, что для заданного

$$(2.4.25) \quad \delta_1 \in R, \ \delta_1 > 0$$

существует $M_i$ такое, что из условия $m > M_i$ следует

$$(2.4.26) \quad a_{i \cdot m} \in B_o(a_i, \delta_1)$$

Пусть
$$M = \max(M_1, ..., M_n)$$
Из равенств (2.4.26), теоремы 2.4.5 и условия $m > M$ следует, что

(2.4.27) $$a_{1 \cdot m}...a_{n \cdot m}\omega \in B_o(a_1...a_n\omega, \epsilon_1)$$

(2.4.28) $$\epsilon_1 = \delta_1 C_2...C_n + ... + \delta_1 C_1...C_{n-1}$$

где

(2.4.29) $$C_i = \|a_i\| + \delta_1$$

Из равенства (2.4.29) и утверждений 2.1.9.1, (2.4.25) следует, что

(2.4.30) $$C_i > 0 \quad \frac{dC_i}{d\delta_1} > 0$$

Из равенств (2.4.28), (2.4.29) и утверждения (2.4.30) следует, что $\epsilon_1$ является полиномиальной строго монотонно возрастающей функцией $\delta_1 > 0$. При этом
$$\delta_1 = 0 \Rightarrow \epsilon_1 = 0$$
Согласно теореме 2.3.5, отображение (2.4.28) отображает интервал $[0, \delta_1)$ в интервал $[0, \epsilon_1)$. Согласно теореме 2.3.3, для заданного $\epsilon > 0$ существует такое $\delta > 0$, что
$$\epsilon_1(\delta) < \epsilon$$
Согласно построению, значение $M$ зависит от значения $\delta_1$. Мы выберем значение $M$, соответствующее $\delta_1 = \delta$. Следовательно, для заданного $\epsilon \in R$, $\epsilon > 0$, существует $M$ такое, что из условия $m > M$ следует

(2.4.31) $$a_{1 \cdot m}...a_{n \cdot m}\omega \in B_o(a_1...a_n\omega, \epsilon)$$

Равенство (2.4.24) следует из равенства (2.4.31) и теоремы 2.1.18.　　□

## 2.5. Пополнение нормированной Ω-группы

Определение 2.5.1. Пусть $A$ - нормированная Ω-группа. Полная Ω-группа $B$ называется **пополнением нормированной Ω-группы** $A$,[2.18] если

2.5.1.1: Ω-группа $A$ является подгруппой Ω-группы $B$.[2.19]
2.5.1.2: Ω-группа $A$ всюду плотна в Ω-группе $B$.

□

---

[2.18]Существование пополнения топологического пространства очень важно, так как позволяет использовать непрерывность как инструмент для изучения топологического пространства. Если топологическое пространство имеет дополнительную структуру, то мы ожидаем, что пополнение имеет такую же структуру. Смотри определение пополнение метрического пространства в определении [3]-2, страница 71. Смотри определение пополнение нормированного поля в утверждении [4]-6 на странице 270 и в последующем определении на странице 271. Смотри определение пополнение нормированного векторного пространства в теореме [5]-1.11.1 на странице 55.

[2.19]Точнее, мы изоморфно отождествляем Ω-группу $A$ с некоторой подгруппой Ω-группы $B$.

Теорема 2.5.2. *Пусть $A$ - нормированная $\Omega$-группа. Пусть $B_1$, $B_2$ являются пополнением $\Omega$-группы $A$. Существует изоморфизм $\Omega$-группы*[2.20]

$$(2.5.1) \qquad f : B_1 \to B_2$$

*такой, что*

$$(2.5.2) \qquad f(a) = a \quad a \in A$$

$$(2.5.3) \qquad \|f(b)\| = \|b\|$$

Доказательство. Пусть $b_1 \in B_1$. Согласно утверждению 2.5.1.2 и теореме 2.2.14 существует последовательность $A$-чисел $a_n$ такая, что

$$(2.5.4) \qquad b_1 = \lim_{n \to \infty} a_n$$

Из равенства (2.5.4) следует, что последовательность $A$-чисел $a_n$ является фундаментальной последовательностью в $\Omega$-группе $A$. Следовательно, последовательность $a_n$ сходится к $b_2$ в полной $\Omega$-группе $B_2$. Согласно теореме 2.1.22, $B_2$-число $b_2$ не зависит от выбора последовательности $a_n$, сходящейся к $b_1$. Следовательно, $f(b_1) = b_2$.

Равенство (2.5.2) является следствием теоремы 2.2.13.

Отображение $f$ является биективным, так как выполненное построение обратимо и мы можем найти $B_1$-число $b_1$, соответствующее $B_2$-числу $b_2$.

Согласно теореме 2.3.9, в $\Omega$-группе $B_1$ верно равенство

$$(2.5.5) \qquad \|b_1\| = \lim_{n \to \infty} \|a_n\|$$

и в $\Omega$-группе $B_2$ верно равенство

$$(2.5.6) \qquad \|b_2\| = \lim_{n \to \infty} \|a_n\|$$

Равенства (2.5.3) следует из равенств (2.5.5), (2.5.6).

Из теорем 2.4.3, 2.4.6 следует, что отображение $f$ является изоморфизмом $\Omega$-группы. $\square$

Доказательство существования пополнения нормированной $\Omega$-группы $A$ (теорема 2.5.17) основано на доказательстве теоремы 2.5.2. Если $B$ - полная $\Omega$-группа, являющаяся пополнением нормированной $\Omega$-группы $A$, то каждое $B$-число является пределом фундаментальной последовательности $A$-чисел. Поэтому, чтобы построить $\Omega$-группу $B$, рассмотрим множество $B'$ фундаментальных последовательностей $A$-чисел.

Лемма 2.5.3. *Отношение $\sim$ на множестве $B'$*

$$(2.5.7) \qquad a_n \sim b_n \Leftrightarrow \lim_{n \to \infty}(a_n - b_n) = 0$$

*является отношением эквивалентности.*[2.21]

Доказательство.

---

[2.20]Мы также говорим, что пополнение $\Omega$-группы $A$ единственно с точностью до рассматриваемого изоморфизма.

[2.21]Смотри определение отношения эквивалентности в [11], страница 27.

2.5.3.1: Так как $a_n = a_n$, то отношение $\sim$ рефлексивно
$$a_n \sim a_n$$

2.5.3.2: Из утверждения 2.1.9.4 следует, что
$$\|a_n - b_n\| = \|b_n - a_n\|$$
Следовательно, отношение $\sim$ симметрично
$$a_n \sim b_n \Leftrightarrow b_n \sim a_n$$

2.5.3.3: Пусть $a_n$, $n = 1, ...,$ - фундаментальная последовательность $A$-чисел. Пусть $b_n$, $c_n$, $n = 1, ...,$ - последовательности $A$-чисел. Равенство
$$\lim_{n \to \infty}(a_n - b_n) = 0 \tag{2.5.8}$$
следует из утверждения $a_n \sim b_n$. Равенство
$$\lim_{n \to \infty}(a_n - c_n) = 0 \tag{2.5.9}$$
следует из утверждения $a_n \sim c_n$.

Из теоремы 2.1.21 и равенства (2.5.8) следует, что последовательность $b_n$ фундаментальна. Согласно определению 2.1.17, из равенства (2.5.8) следует, что для любого $\epsilon \in R$, $\epsilon > 0$, существует, зависящее от $\epsilon$, натуральное число $N_1$ такое, что
$$\|a_n - b_n\| < \frac{\epsilon}{2} \tag{2.5.10}$$
для любого $n > N_1$.

Из теоремы 2.1.21 и равенства (2.5.9) следует, что последовательность $c_n$ фундаментальна. Согласно определению 2.1.17, из равенства (2.5.9) следует, что для любого $\epsilon \in R$, $\epsilon > 0$, существует, зависящее от $\epsilon$, натуральное число $N_2$ такое, что
$$\|a_n - c_n\| < \frac{\epsilon}{2} \tag{2.5.11}$$
для любого $n > N_2$.

Пусть
$$N = \max(N_1, N_2)$$
Из неравенств (2.5.10), (2.5.11) и утверждения 2.1.9.3 следует, что для заданного $\epsilon \in R$, $\epsilon > 0$, существует, зависящее от $\epsilon$, натуральное число $N$ такое, что
$$\|b_n - c_n\| = \|b_n - a_n + a_n - c_n\| \leq \|b_n - a_n\| + \|a_n - c_n\| < \epsilon$$
для любого $n > N$. Равенство
$$\lim_{n \to \infty}(b_n - c_n) = 0 \tag{2.5.12}$$
следует из определения 2.1.17. Утверждение $b_n \sim c_n$ следует из равенства (2.5.12).

2.5.3.4: Из рассуждения 2.5.3.3 следует, что отношение $\sim$ транзитивно
$$a_n \sim b_n, a_n \sim c_n \Rightarrow b_n \sim c_n$$

Лемма следует из утверждений 2.5.3.1, 2.5.3.1, 2.5.3.4.    □

Лемма 2.5.4. $a_n \sim b_n$ тогда и только тогда, когда для любого $\epsilon \in R$, $\epsilon > 0$, существует, зависящее от $\epsilon$, натуральное число $n_0$ такое, что

(2.5.13) $$b_n \in B_o(a_n, \epsilon)$$

для любого $n > n_0$.

Доказательство. Согласно определению 2.1.14, утверждение (2.5.13) верно тогда и только тогда, когда

$$\|a_n - b_n\| < \epsilon$$

Теорема верна согласно определению 2.1.17. $\square$

Пусть $B$ является множеством классов эквивалентных фундаментальных последовательностей $A$-чисел. Мы будем пользоваться записью

(2.5.14) $$[a_n] = \{b_n : a_n \sim b_n\}$$

для класса последовательностей, эквивалентных последовательности $a_n$.

Лемма 2.5.5. Множество $B$ является абелевой группой относительно операции

(2.5.15) $$[a_n] + [b_n] = [a_n + b_n]$$

Доказательство. Пусть $A$ - нормированная Ω-группа. Пусть последовательности $A$-чисел $a_n, b_n$ являются фундаментальными последовательностями.

Лемма 2.5.6. Последовательность $A$-чисел $a_n + b_n$ является фундаментальной.

Доказательство. Так как последовательность $A$-чисел $a_n$ является фундаментальной, то согласно определению 2.1.19 для любого $\epsilon \in R$, $\epsilon > 0$, существует, зависящее от $\epsilon$, натуральное число $N_a$ такое, что

(2.5.16) $$\|a_p - a_q\| < \frac{\epsilon}{2}$$

для любых $p, q > N_a$. Так как последовательность $A$-чисел $b_n$ является фундаментальной, то согласно определению 2.1.19 для любого $\epsilon \in R$, $\epsilon > 0$, существует, зависящее от $\epsilon$, натуральное число $N_b$ такое, что

(2.5.17) $$\|b_p - b_q\| < \frac{\epsilon}{2}$$

для любых $p, q > N_b$. Пусть

$$N = \max(N_a, N_b)$$

Из неравенств (2.4.1), (2.5.16), (2.5.17) и условия $n > N$ следует, что

$$\|(a_p + b_p) - (a_q + b_q)\| \leq \|(a_p - a_q)\| + \|(b_p - b_q)\| < \epsilon$$

Следовательно, согласно определению 2.1.19, последовательность $A$-чисел $a_n + b_n$ является фундаментальной.

Лемма 2.5.7. Пусть $a_n \sim c_n$, $b_n \sim d_n$. Тогда

(2.5.18) $$a_n + b_n \sim c_n + d_n$$

Доказательство. Равенство

(2.5.19) $$\lim_{n\to\infty}(a_n - c_n) = 0$$

следует из утверждения $a_n \sim c_n$. Равенство

(2.5.20) $$\lim_{n\to\infty}(b_n - d_n) = 0$$

следует из утверждения $b_n \sim d_n$.

Из теоремы 2.1.21 и равенства (2.5.19) следует, что последовательность $c_n$ фундаментальна. Согласно определению 2.1.17, из равенства (2.5.19) следует, что для любого $\epsilon \in R$, $\epsilon > 0$, существует, зависящее от $\epsilon$, натуральное число $N_1$ такое, что

(2.5.21) $$\|a_n - c_n\| < \frac{\epsilon}{2}$$

для любого $n > N_1$.

Из теоремы 2.1.21 и равенства (2.5.20) следует, что последовательность $d_n$ фундаментальна. Согласно определению 2.1.17, из равенства (2.5.20) следует, что для любого $\epsilon \in R$, $\epsilon > 0$, существует, зависящее от $\epsilon$, натуральное число $N_2$ такое, что

(2.5.22) $$\|b_n - d_n\| < \frac{\epsilon}{2}$$

для любого $n > N_2$.

Пусть
$$N = \max(N_1, N_2)$$

Из неравенств (2.4.1), (2.5.21), (2.5.22) следует, что для заданного $\epsilon \in R$, $\epsilon > 0$, существует, зависящее от $\epsilon$, натуральное число $N$ такое, что

(2.5.23) $$\|(a_n + b_n) - (c_n + d_n)\| \leq \|a_n - c_n\| + \|b_n - d_n\| < \epsilon$$

для любого $n > N$. Равенство

(2.5.24) $$\lim_{n\to\infty}((a_n + b_n) - (c_n + d_n)) = 0$$

следует из определения 2.1.17. Следовательно, утверждение (2.5.18) следует из равенства (2.5.24).

Из леммы 2.5.6 следует, что последовательность $a_n + b_n$ является фундаментальной последовательностью. Из леммы 2.5.7 следует, что операция (2.5.15) определена корректно.

2.5.7.1: Из равенства
$$a_n + (b_n + c_n) = (a_n + b_n) + c_n$$
следует, что операция (2.5.15) ассоциативна
$$[a_n] + ([b_n] + [c_n]) = ([a_n] + [b_n]) + [c_n]$$

2.5.7.2: Из равенства
$$a_n + b_n = b_n + a_n$$
следует, что операция (2.5.15) коммутативна
$$[a_n] + [b_n] = [b_n] + [a_n]$$

2.5.7.3: Если положить $b_n = 0$, то из равенства
$$a_n + 0 = a_n$$
следует
$$[a_n] + [0] = [a_n]$$

2.5.7.4: Из равенства
$$a_n + (-a_n) = 0$$
следует
$$[a_n] + [-a_n] = [0]$$
Следовательно, мы можем положить
$$-[a_n] = [-a_n]$$

Лемма является следствием утверждений 2.5.7.1, 2.5.7.2, 2.5.7.3, 2.5.7.4. □

Лемма 2.5.8. Мы определим операцию $\omega \in \Omega$ на множестве $B$, пользуясь равенством

(2.5.25) $$[a_{1 \cdot m}]...[a_{n \cdot m}]\omega = [a_{1 \cdot m}...a_{n \cdot m}\omega]$$

Доказательство. Пусть $A$ - нормированная Ω-группа.

Лемма 2.5.9. Пусть последовательности $A$-чисел $a_{1 \cdot m}, ..., a_{n \cdot m}$, $m = 1, ...,$ являются фундаментальными последовательностями. Последовательность $A$-чисел $a_{1 \cdot m}...a_{n \cdot m}\omega$ является фундаментальной.

Доказательство. Так как последовательность $a_{i \cdot m}$ является фундаментальной последовательностью $A$-чисел, то согласно теореме 2.1.20 следует, что для заданного

(2.5.26) $$\delta_1 \in R, \ \delta_1 > 0$$

существует, зависящее от $\delta_1$, натуральное число $M_i$ такое, что

(2.5.27) $$a_{i \cdot q} \in B_o(a_{i \cdot p}, \delta_1)$$

для любых $p, q > M_i$. Пусть
$$M = \max(M_1, ..., M_n)$$

Из утверждений (2.5.27) и теоремы 2.4.5 следует, что
$$a_{1 \cdot q}...a_{n \cdot q}\omega \in B_o(a_{1 \cdot p}...a_{n \cdot p}\omega, \epsilon_1)$$

(2.5.28) $$\epsilon_1 = \delta_1 C_2...C_n + ... + \delta_1 C_1...C_{n-1}$$

где

(2.5.29) $$C_i = \|a_{i \cdot p}\|_i + \delta_1$$

Из равенства (2.5.29) и утверждений 2.1.9.1, (2.5.26) следует, что

(2.5.30) $$C_i > 0 \quad \frac{dC_i}{d\delta_1} > 0$$

Из равенств (2.5.28), (2.5.29) и утверждения (2.5.30) следует, что $\epsilon_1$ является полиномиальной строго монотонно возрастающей функцией $\delta_1 > 0$. При этом
$$\delta_1 = 0 \Rightarrow \epsilon_1 = 0$$
Согласно теореме 2.3.5, отображение (2.5.28) отображает интервал $[0, \delta_1)$ в интервал $[0, \epsilon_1)$. Согласно теореме 2.3.3, для заданного $\epsilon > 0$ существует такое $\delta > 0$, что
$$\epsilon_1(\delta) < \epsilon$$
Согласно построению, значение $M$ зависит от значения $\delta_1$. Мы выберем значение $M$, соответствующее $\delta_1 = \delta$. Следовательно, для заданного $\epsilon \in R$, $\epsilon > 0$, существует $M$ такое, что из условия $p, q > M$ следует

(2.5.31) $$a_{1 \cdot q}...a_{n \cdot q}\omega \in B_o(a_{1 \cdot p}...a_{n \cdot p}\omega, \epsilon)$$

Из утверждения (2.5.31) и теоремы 2.1.20 следует, что последовательность $A$-чисел $a_{1 \cdot m}...a_{n \cdot m}\omega$ является фундаментальной последовательностью. □

ЛЕММА 2.5.10. *Пусть последовательности $A$-чисел $a_{1 \cdot m}, ..., a_{n \cdot m}$, $m = 1, ...,$ $c_{1 \cdot m}, ..., c_{n \cdot m}$, $m = 1, ...,$ являются фундаментальными последовательностями. Пусть*

(2.5.32) $$a_{i \cdot m} \sim c_{i \cdot m}, \quad i = 1, ..., n$$

*Тогда*

(2.5.33) $$a_{1 \cdot m}...a_{n \cdot m}\omega \sim c_{1 \cdot m}...c_{n \cdot m}\omega$$

ДОКАЗАТЕЛЬСТВО. Из утверждения (2.5.32) и леммы 2.5.4 следует, что для заданного

(2.5.34) $$\delta_1 \in R, \ \delta_1 > 0$$

существует, зависящее от $\delta_1$, натуральное число $M_1$ такое, что

(2.5.35) $$c_{i \cdot m} \in B_o(a_{i \cdot m}, \delta_1)$$

для любого $m > M_1$. Пусть
$$M = \max(M_1, ..., M_n)$$

Из утверждений (2.5.35) и теоремы 2.4.5 следует, что
$$c_{1 \cdot m}...c_{n \cdot m}\omega \in B_o(a_{1 \cdot m}...a_{n \cdot m}\omega, \epsilon_1)$$

(2.5.36) $$\epsilon_1 = \delta_1 C_2...C_n + ... + \delta_1 C_1...C_{n-1}$$

где

(2.5.37) $$C_i = \|a_{i \cdot n}\|_i + \delta_1$$

Из равенства (2.5.37) и утверждений 2.1.9.1, (3.3.23) следует, что

(2.5.38) $$C_i > 0 \quad \frac{dC_i}{d\delta_1} > 0$$

Из равенств (2.5.36), (2.5.37) и утверждения (2.5.38) следует, что $\epsilon_1$ является строго монотонно возрастающей функцией $\delta_1 > 0$. При этом
$$\delta_1 = 0 \Rightarrow \epsilon_1 = 0$$

Согласно теореме 2.3.5, отображение (2.5.36) отображает интервал $[0, \delta_1)$ в интервал $[0, \epsilon_1)$. Согласно теореме 2.3.3, для заданного $\epsilon > 0$ существует такое $\delta > 0$, что
$$\epsilon_1(\delta) < \epsilon$$
Согласно построению, значение $M$ зависит от значения $\delta_1$. Мы выберем значение $M$, соответствующее $\delta_1 = \delta$. Следовательно, для заданного $\epsilon \in R, \epsilon > 0$, существует $M$ такое, что из условия $m > M$ следует

(2.5.39) $$c_{1 \cdot m}...c_{n \cdot m}\omega \in B_o(a_{1 \cdot m}...a_{n \cdot m}\omega, \epsilon)$$

Утверждение (2.5.33) следует из утверждения (2.5.39) и леммы 2.5.4.

Из леммы 2.5.9 следует, что последовательность $a_{1 \cdot m}...a_{n \cdot m}\omega$ является фундаментальной последовательностью. Из леммы 2.5.10 следует, что операция (2.5.25) определена корректно. □

Лемма 2.5.11. *Абелева группа $B$ является абелевой Ω-группой.*

Доказательство. Согласно определению 2.1.3, операция $\omega \in \Omega$ полиаддитивна. Из равенства
$$a_{1 \cdot m}...(a_{i \cdot m} + b_{i \cdot m})...a_{n \cdot m}\omega = a_{1 \cdot m}...a_{i \cdot m}...a_{n \cdot m}\omega + a_{1 \cdot m}...b_{i \cdot m}...a_{n \cdot m}\omega$$
следует, что операция (2.5.25) полиаддитивна
$$[a_{1 \cdot m}]...[a_{i \cdot m} + b_{i \cdot m}]...[a_{n \cdot m}]\omega = [a_{1 \cdot m}]...[a_{i \cdot m}]...[a_{n \cdot m}]\omega + [a_{1 \cdot m}]...[b_{i \cdot m}]...[a_{n \cdot m}]\omega$$
Лемма является следствием определения 2.1.3. □

Лемма 2.5.12. *Пусть последовательность $a_n$ является фундаментальной. Тогда последовательность $\|a_n\|$ является фундаментальной.*

Доказательство. Согласно определению 2.1.19, для заданного $\epsilon \in R$, $\epsilon > 0$, существует, зависящее от $\epsilon$, натуральное число $N$ такое, что

(2.5.40) $$\|a_p - a_q\| < \epsilon$$

для любых $p, q > N$. Из неравенств (2.1.1), (2.5.40) следует, что для любого $\epsilon \in R, \epsilon > 0$, существует, зависящее от $\epsilon$, натуральное число $N$ такое, что

(2.5.41) $$|\|a_p\| - \|a_q\|| \leq \|a_p - a_q\| < \epsilon$$

для любого $n > N$. Согласно определению 2.1.19, последовательность $\|a_n\|$ является фундаментальной. □

Лемма 2.5.13. *Пусть $a_n \sim b_n$. Тогда*

(2.5.42) $$\lim_{n \to \infty} \|a_n\| = \lim_{n \to \infty} \|b_n\|$$

Доказательство. Равенство

(2.5.43) $$\lim_{n \to \infty}(a_n - b_n) = 0$$

следует из утверждения $a_n \sim b_n$. Согласно определению 2.1.17, из равенства (2.5.43) следует, что для любого $\epsilon \in R, \epsilon > 0$, существует, зависящее от $\epsilon$, натуральное число $N$ такое, что
$$\|a_n - b_n\| < \epsilon$$

для любого $n > N$. Из неравенства (2.1.1) следует, что для любого $\epsilon \in R, \epsilon > 0$, существует, зависящее от $\epsilon$, натуральное число $N$ такое, что

$$(2.5.44) \qquad |\,\|a_n\| - \|b_n\|\,| \le \|a_n - b_n\| < \epsilon$$

для любого $n > N$. Согласно определению 2.1.17, следует, что

$$(2.5.45) \qquad \lim_{n \to \infty}(\|a_n\| - \|b_n\|) = 0$$

Из леммы 2.5.12 следует, что последовательности $\|a_n\|$ и $\|b_n\|$ фундаментальны. Равенство (2.5.42) следует из теоремы 2.1.22. $\square$

Лемма 2.5.14. Мы определим норму на $\Omega$-группе $B$ равенством

$$(2.5.46) \qquad \|\,[a_n]\,\| = \lim_{n \to \infty} \|a_n\|$$

Доказательство. Из лемм 2.5.12, 2.5.13 следует, что равенство (2.5.46) является корректным определением отображения

$$B \to R$$

2.5.14.1: Согласно утверждению 2.1.9.1, $\|a_n\| \ge 0$. Следовательно,[2.22]

$$(2.5.47) \qquad \lim_{n \to \infty} \|a_n\| \ge 0$$

Из равенства (2.5.46) и неравенства (2.5.47) следует, что $\|\,[a_n]\,\| \ge 0$. Следовательно, утверждение 2.1.9.1 верно для отображения (2.5.46).

2.5.14.2: Пусть $\|\,[a]\,\| = 0$. Из равенства (2.5.46) следует, что

$$(2.5.48) \qquad \lim_{n \to \infty} \|a_n\| = 0$$

Согласно определению 2.1.17, из равенства (2.5.48) следует, что для любого $\epsilon \in R, \epsilon > 0$, существует, зависящее от $\epsilon$, натуральное число $N$ такое, что $\|a_n\| < \epsilon$ для любого $n > N$. Согласно определению 2.1.17,

$$(2.5.49) \qquad \lim_{n \to \infty} a_n = \lim_{n \to \infty}(a_n - 0) = 0$$

Из равенств (2.5.49), (2.5.7) следует, что $a_n \sim 0$. Согласно определению (2.5.14), $[a_n] = [0]$. Следовательно, утверждение 2.1.9.2 верно для отображения (2.5.46).

2.5.14.3: Согласно утверждению 2.1.9.3,

$$\|a_n + b_n\| \le \|a_n\| + \|b_n\|$$

Следовательно,[2.23]

$$(2.5.50) \qquad \lim_{n \to \infty} \|a_n + b_n\| \le \lim_{n \to \infty} \|a_n\| + \lim_{n \to \infty} \|b_n\|$$

Из равенства (2.5.46) и неравенства (2.5.50) следует, что

$$(2.5.51) \qquad \|\,[a_n + b_n]\,\| \le \|\,[a_n]\,\| + \|\,[b_n]\,\|$$

Из неравенства (2.5.51) и равенства (2.5.15) следует, что

$$\|\,[a_n] + [b_n]\,\| \le \|\,[a_n]\,\| + \|\,[b_n]\,\|$$

Следовательно, утверждение 2.1.9.3 верно для отображения (2.5.46).

---

[2.22]Смотри, например, [10], страница 58, утверждение 3.

[2.23]Смотри, например, [10], страница 58, утверждение 3.

2.5.14.4: Согласно утверждению 2.1.9.4, $\|-a_n\| = \|a_n\|$. Согласно теореме 2.1.22

$$\lim_{n \to \infty} \|-a_n\| = \lim_{n \to \infty} \|a_n\| \tag{2.5.52}$$

Из равенств (2.5.46), (2.5.52) следует, что

$$\|[-a_n]\| = \|[a_n]\| \tag{2.5.53}$$

Из равенства (2.5.53) и утверждения 2.5.7.4 следует, что

$$\|-[a_n]\| = \|[a_n]\|$$

Следовательно, утверждение 2.1.9.4 верно для отображения (2.5.46).

Согласно утверждениям 2.5.14.1, 2.5.14.2, 2.5.14.3, 2.5.7.4 и определению 2.1.9, отображение (2.5.46) является нормой на Ω-группе $B$. $\square$

Лемма 2.5.15. Ω-группа $A$ является подгруппой Ω-группы $B$. Из условий $a \in A$ и

$$a = \lim_{n \to \infty} a_n$$

следует

$$\|a\| = \|[a_n]\| \tag{2.5.54}$$

Доказательство. Пусть $a \in A$. Согласно определениям 2.1.17, 2.1.19, последовательность $a_n = a$, $n = 1, ...,$ является фундаментальной и сходится к $a$. Согласно определению (2.5.7) и теореме 2.1.22, из условия $a_n \sim b_n$ следует

$$a = \lim_{n \to \infty} b_n$$

Поэтому мы можем отождествить $a$ и $[a_n]$. Согласно теоремам 2.4.3, 2.4.6, это отождествление является гомоморфизмом Ω-группы $A$ в Ω-группу $B$. Из равенства (2.5.46) и теоремы 2.3.9 следует, что

$$\|[a_n]\| = \lim_{n \to \infty} \|a_n\| = \|a\| \tag{2.5.55}$$

Равенство (2.5.54) следует из равенства (2.5.55). $\square$

Лемма 2.5.16. Для любого $B$-числа $a$ существует последовательность $A$-чисел такая, что

$$a = \lim_{n \to \infty} a_n \tag{2.5.56}$$

Доказательство. Существует фундаментальная последовательность $A$-чисел $a_n$ такая, что $a = [a_n]$. Согласно определению 2.1.19, для заданного $\epsilon \in R$, $\epsilon > 0$, существует, зависящее от $\epsilon$, натуральное число $N_1$ такое, что

$$\|a_p - a_q\| < \frac{\epsilon}{2} \tag{2.5.57}$$

для любых $p, q > N_1$. Согласно леммам 2.5.5, 2.5.15, для любого $q > N_1$

$$a - a_q = [a_p] - [a_q] = [a_p - a_q] \tag{2.5.58}$$

Согласно лемме 2.5.14

$$\|a - a_q\| = \lim_{p \to \infty} \|a_p - a_q\| \tag{2.5.59}$$

Согласно определению 2.1.17, из равенства (2.5.59) следует, что для любого $\epsilon \in R$, $\epsilon > 0$, существует, зависящее от $\epsilon$, натуральное число $N_2$ такое, что

(2.5.60) $$|\,\|a - a_q\| - \|a_p - a_q\|\,| < \frac{\epsilon}{2}$$

для любого $n > N_2$. Неравенство[2.24]

(2.5.63) $$\|a - a_q\| - \|a_p - a_q\| < \frac{\epsilon}{2}$$

следует из неравенства (2.5.60). Пусть

$$N = \max(N_1, N_2)$$

Из неравенств (2.5.57), (2.5.63) следует, что для заданного $\epsilon \in R$, $\epsilon > 0$, существует, зависящее от $\epsilon$, натуральное число $N$ такое, что

(2.5.64) $$\|a - a_q\| < \|a_p - a_q\| + \frac{\epsilon}{2} < \epsilon$$

для любых $q > N$. Равенство (2.5.56) следует из определения 2.1.17. □

Теорема 2.5.17. *Пополнение нормированной Ω-группы A существует.*

Доказательство. Согласно лемме 2.5.15, Ω-группа $A$ является подгруппой Ω-группы $B$. Из леммы 2.5.16 и теоремы 2.2.14 следует, что Ω-группа $A$ всюду плотна в Ω-группе $B$. Согласно определению 2.5.1, чтобы доказать теорему, необходимо доказать, что Ω-группа $B$ является полной.

Пусть $b_n$ - фундаментальная последовательность $B$-чисел. Согласно определениям 2.1.14, 2.2.3, 2.2.5, для любого $n > 0$ существует $a_n \in A$ такой, что

(2.5.65) $$\|a_n - b_n\| < \frac{1}{n}$$

Равенство

(2.5.66) $$\lim_{n \to \infty}(a_n - b_n) = 0$$

следует из определения 2.1.17. Из равенства (2.5.66) и теоремы 2.1.21 следует, что последовательность $A$-чисел $a_n$ является фундаментальной. Согласно определению (2.5.14) существует $B$-число $[a_n]$. Согласно теореме 2.1.22 и лемме 2.5.16

(2.5.67) $$\lim_{n \to \infty} b_n = \lim_{n \to \infty} a_n = [a_n]$$

Из равенства (2.5.67) следует, что Ω-группа $B$ является полной. □

---

[2.24]Неравенство

(2.5.61) $$\|a - a_q\| - \|a_p - a_q\| > -\frac{\epsilon}{2}$$

также следует из неравенства (2.5.60). Неравенство

(2.5.62) $$\|a - a_q\| > \|a_p - a_q\| - \frac{\epsilon}{2}$$

следует из неравенства (2.5.61). Из неравенств (2.5.57), (2.5.62) и утверждения 2.1.9.1 следует, что

$$\|a - a_q\| \geq 0$$

Однако это утверждение очевидно и не интересно для нас.

### 2.6. Ω-группа отображений

Теорема 2.6.1. *Пусть $M(X, A)$ - множество отображений множества $X$ в Ω-группу $A$. Мы можем определить структуру Ω-группы на множестве $M(X, A)$.*

Доказательство. Пусть $f, g \in M(X, A)$. Тогда мы положим
$$(f + g)(x) = f(x) + g(x)$$
Пусть $\omega \in \Omega$ - $n$-арная операция. Для отображений $f_i \in M(X, A), i = 1, ..., n,$ мы положим
$$(2.6.1) \qquad (f_1...f_n\omega)(x) = f_1(x)...f_n(x)\omega$$
□

Так как $X$ - произвольное множество, мы не можем определить норму в Ω-группе $M(X, A)$. Однако мы можем определить сходимость последовательности в $M(X, A)$; следовательно, мы можем определить топологию в $M(X, A)$

Определение 2.6.2. Пусть $f_n \in M(X, A), n = 1, ...,$ - последовательность отображений в нормированную Ω-группу $A$. Отображение $f \in M(X, A)$ называется **пределом последовательности** $f_n$, если для любого $x \in X$
$$f(x) = \lim_{n\to\infty} f_n(x)$$
Мы будем также говорить, что **последовательность** $f_n$ **сходится** к отображению $f$. □

Из определений 2.1.17, 2.6.2 следует, что если последовательность $f_n$ сходится к $f$, то для любого $\epsilon \in R, \epsilon > 0,$ существует $N(x)$ такое, что
$$\|f_n(x) - f(x)\| < \epsilon$$
для любого $n > N(x)$.

Определение 2.6.3. Пусть $f_n \in M(X, A), n = 1, ...,$ - последовательность отображений в нормированную Ω-группу $A$. **Последовательность** $f_n$ **сходится равномерно** к отображению $f$, если для любого $\epsilon \in R, \epsilon > 0,$ существует $N$ такое, что
$$\|f_n(x) - f(x)\| < \epsilon$$
для любого $n > N$. □

Теорема 2.6.4. *Последовательность отображений $f_n \in M(X, A), n = 1, ...,$ в нормированную Ω-группу $A$ сходится равномерно к отображению $f$, если для любого $\epsilon \in R, \epsilon > 0,$ существует $N$ такое, что $f_n(x) \in B_o(f(x), \epsilon)$ для любого $n > N$.*

Доказательство. Следствие определений 2.1.14, 2.6.3. □

Теорема 2.6.5. *Последовательность отображений $f_n \in M(X, A), n = 1, ...,$ в нормированную Ω-группу $A$ сходится равномерно к отображению $f$, если для любого $\epsilon \in R, \epsilon > 0,$ существует $N$ такое, что*
$$(2.6.2) \qquad \|f_n(x) - f_m(x)\| < \epsilon$$
*для любых $n, m > N$.*

Доказательство. Согласно определению 2.6.3, для любого $\epsilon \in R$, $\epsilon > 0$, существует $N$ такое, что

$$\|f_n(x) - f(x)\| < \frac{\epsilon}{2} \tag{2.6.3}$$

для любого $n > N$. Для любых $n, m > N$, из утверждения 2.1.9.3 и неравенства (2.6.3) следует, что

$$\begin{aligned}\|f_n(x) - f_m(x)\| &= \|f_n(x) - f(x) + f(x) - f_m(x)\| \\ &\leq \|f_n(x) - f(x)\| + \|f_m(x) - f(x)\| < \epsilon\end{aligned} \tag{2.6.4}$$

Неравенство (2.6.2) следует из неравенства (2.6.4). $\square$

Теорема 2.6.6. *Последовательность отображений $f_n \in M(X, A)$, $n = 1, ...$, в нормированную $\Omega$-группу $A$ сходится равномерно к отображению $f$, если для любого $\epsilon \in R$, $\epsilon > 0$, существует $N$ такое, что $f_n(x) \in B_o(f_m(x), \epsilon)$ для любых $n, m > N$.*

Доказательство. Следствие определения 2.1.14 и теоремы 2.6.5. $\square$

Теорема 2.6.7. *Пусть последовательность отображений $f_n \in M(X, A)$, $n = 1, ...$, в полную $\Omega$-группу $A$ сходится равномерно к отображению $f$. Пусть последовательность отображений $g_n \in M(X, A)$, $n = 1, ...$, в полную $\Omega$-группу $A$ сходится равномерно к отображению $g$. Тогда последовательность отображений*

$$h_n = f_n + g_n$$

*в полную $\Omega$-группу $A$ сходится равномерно к отображению*

$$h = f + g \tag{2.6.5}$$

Доказательство. Из равенства

$$f(x) = \lim_{n \to \infty} f_n(x)$$

и теоремы 2.6.4 следует, что для заданного $\epsilon \in R$, $\epsilon > 0$, существует $N_a$ такое, что из условия $n > N_a$ следует

$$f_n(x) \in B_o(f(x), \epsilon/2) \tag{2.6.6}$$

Из равенства

$$g(x) = \lim_{n \to \infty} g_n(x)$$

и теоремы 2.6.4 следует, что для заданного $\epsilon \in R$, $\epsilon > 0$, существует $N_b$ такое, что из условия $n > N_b$ следует

$$g_n(x) \in B_o(g(x), \epsilon/2) \tag{2.6.7}$$

Пусть

$$N = \max(N_a, N_b)$$

Из равенств (2.6.6), (2.6.7), теоремы 2.4.2 и условия $n > N$ следует, что

$$h_n(x) = f_n(x) + g_n(x) \in B_o(f(x) + g(x), \epsilon/2 + \epsilon/2) = B_o(h(x), \epsilon) \tag{2.6.8}$$

Равенство (2.6.5) следует из равенства (2.6.8) и теоремы 2.6.4. $\square$

Теорема 2.6.8. *Пусть $A$ - полная $\Omega$-группа. Пусть $\omega \in \Omega$ - $n$-арная операция. Пусть последовательность отображений $f_{i \cdot m} \in M(X, A)$, $i = 1$, ..., $n$, $m = 1$, ..., в полную $\Omega$-группу $A$ сходится равномерно к отображению $f_i$. Пусть множество значений отображения $f_i$ компактно. Тогда последовательность отображений*

$$h_m = f_{1 \cdot m} ... f_{n \cdot m} \omega$$

*в полную $\Omega$-группу $A$ сходится равномерно к отображению*

(2.6.9) $$h = f_1 ... f_n \omega$$

Доказательство. Из равенства

$$f_i(x) = \lim_{m \to \infty} f_{i \cdot m}(x)$$

и теоремы 2.6.4 следует, что для заданного

(2.6.10) $$\delta_1 \in R, \ \delta_1 > 0$$

существует $M_i$ такое, что из условия $m > M_i$ следует

(2.6.11) $$f_{i \cdot m}(x) \in B_o(f_i(x), \delta_1)$$

Пусть

$$M = \max(M_1, ..., M_n)$$

Из равенств (2.6.11), теоремы 2.4.5 и условия $m > M$ следует, что

(2.6.12) $$h_m(x) = f_{1 \cdot m}(x) ... f_{n \cdot m}(x) \omega \in B_o(f_1(x) ... f_n(x) \omega, \epsilon_1) = B_o(h(x), \epsilon_1)$$

(2.6.13) $$\epsilon_1 = \delta_1 C_2 ... C_n + ... + \delta_1 C_1 ... C_{n-1}$$

где

(2.6.14) $$C_i = \|f_i(x)\| + \delta_1$$

Значение $C_i$ зависит от $x$. Поэтому $\epsilon_1$ также зависит от $x$. Для оценки равномерной сходимости нам необходимо выбрать максимальное значение $\epsilon_1$. Для заданного $\delta_1$, значение $\epsilon_1$ является полилинейным отображением значений $C_i$. Следовательно, $\epsilon_1$ имеет максимальное значение, когда каждое $C_i$ имеет максимальное значение. Так как мы рассматриваем оценку значения $\epsilon_1$ справа, для нас неважно одновременно ли $C_i$ приобретают максимальные значения. Так как множество значений отображения $f_i$ компактно, то, согласно теореме 2.3.11, определена следующая величина

(2.6.15) $$F_i = \sup \|f_i(x)\|$$

Из равенства (2.6.15) и утверждения 2.1.9.1 следует, что

(2.6.16) $$F_i \geq 0$$

Из равенств (2.6.14), (2.6.15) следует, что мы можем положить

(2.6.17) $$C_i = F_i + \delta_1$$

Из равенства (2.6.17) и утверждений (2.6.16), (2.6.10) следует, что

(2.6.18) $$C_i > 0 \quad \frac{dC_i}{d\delta_1} > 0$$

Из равенств (2.6.13), (2.6.17) и утверждения (2.6.18) следует, что $\epsilon_1$ является полиномиальной строго монотонно возрастающей функцией $\delta_1 > 0$. При этом

$$\delta_1 = 0 \Rightarrow \epsilon_1 = 0$$

Согласно теореме 2.3.5, отображение (2.6.13) отображает интервал $[0, \delta_1)$ в интервал $[0, \epsilon_1)$. Согласно теореме 2.3.3, для заданного $\epsilon > 0$ существует такое $\delta > 0$, что

$$\epsilon_1(\delta) < \epsilon$$

Согласно построению, значение $M$ зависит от значения $\delta_1$. Мы выберем значение $M$, соответствующее $\delta_1 = \delta$. Следовательно, для заданного $\epsilon \in R$, $\epsilon > 0$, существует $M$ такое, что из условия $m > M$ следует

(2.6.19) $\quad h_m(x) = f_{1 \cdot m}(x)...f_{n \cdot m}(x)\omega \in B_o(f_1(x)...f_n(x)\omega, \epsilon) = B_o(h(x), \epsilon)$

Равенство (2.6.9) следует из равенства (2.6.19) и теоремы 2.6.4. $\square$

Глава 3

# Представление $\Omega$-группы

### 3.1. Представление $\Omega$-группы

Определение 3.1.1. Пусть

$$f: A_1 \dashrightarrow A_2$$

представление[3.1] $\Omega_1$-группы $A_1$ с нормой $\|x\|_1$ в $\Omega_2$-группе $A_2$ с нормой $\|x\|_2$. Величина

$$(3.1.1) \qquad \|f\| = \sup \frac{\|f(a_1)(a_2)\|_2}{\|a_1\|_1 \|a_2\|_2}$$

называется **нормой представления** $f$. $\square$

Теорема 3.1.2. *Пусть*

$$f: A_1 \dashrightarrow A_2$$

*представление $\Omega_1$-группы $A_1$ с нормой $\|x\|_1$ в $\Omega_2$-группе $A_2$ с нормой $\|x\|_2$. Тогда*

$$(3.1.2) \qquad \|f(a_1)(a_2)\|_2 \le \|f\| \|a_1\|_1 \|a_2\|_2$$

Доказательство. Из равенства (3.1.1) следует, что

$$(3.1.3) \qquad \frac{\|f(a_1)(a_2)\|_2}{\|a_1\|_1 \|a_2\|_2} \le \sup \frac{\|f(a_1)(a_2)\|_2}{\|a_1\|_1 \|a_2\|_2} = \|f\|$$

Неравенство (3.1.2) следует из неравенства (3.1.3). $\square$

Теорема 3.1.3. *Пусть*

$$f: A_1 \dashrightarrow A_2$$

*представление $\Omega_1$-группы $A_1$ с нормой $\|x\|_1$ в $\Omega_2$-группе $A_2$ с нормой $\|x\|_2$. Верно следующее неравенство*

$$(3.1.4) \qquad \|f(c_1)(c_2) - f(a_1)(a_2)\|_2 \le \|f\|(\|c_1 - a_1\|_1 C_2 + C_1 \|c_2 - a_2\|_2)$$

*где*

$$(3.1.5) \qquad C_1 = \max(\|a_1\|_1, \|c_1\|_1) \quad C_2 = \max(\|a_2\|_2, \|c_2\|_2)$$

---

[3.1]Смотри определение [7]-2.1.4, представления универсальной алгебры. Согласно определению [9]-2.1.2, модуль - это представление кольца в абелевой группе. Так как кольцо и абелева группа являются $\Omega$-группами, то модуль - это представление $\Omega$-группы.

Доказательство. Согласно определениям 2.1.3, [7]-2.1.4,

$$(3.1.6) \quad \begin{aligned} f(c_1)(c_2) - f(a_1)(a_2) &= f(c_1)(c_2) - f(a_1)(c_2) + f(a_1)(c_2) - f(a_1)(a_2) \\ &= f(c_1 - a_1)(c_2) + f(a_1)(c_2 - a_2) \end{aligned}$$

Из равенства (3.1.6) и утверждения 2.1.9.3 следует, что

$$(3.1.7) \quad \|f(c_1)(c_2) - f(a_1)(a_2)\|_2 \leq \|f(c_1 - a_1)(c_2)\|_2 + \|f(a_1)(c_2 - a_2)\|_2$$

Из равенства (3.1.5) и теоремы 3.1.2 следует, что

$$(3.1.8) \quad \begin{aligned} \|f(c_1 - a_1)(c_2)\|_2 &\leq \|f\| \, \|c_1 - a_1\|_1 \, C_2 \\ \|f(a_1)(c_2 - a_2)\|_2 &\leq \|f\| \, C_1 \, \|c_2 - a_2\|_2 \end{aligned}$$

Утверждение (3.1.4) следует из неравенств (3.1.7), (3.1.8). $\square$

Теорема 3.1.4. *Пусть*

$$f : A_1 \dashrightarrow A_2$$

*представление $\Omega_1$-группы $A_1$ с нормой $\|x\|_1$ в $\Omega_2$-группе $A_2$ с нормой $\|x\|_2$. Для $c_1 \in A_1$, пусть $c_1 \in B_c(a_1 \in A_1, R_1)$. Для $c_2 \in A_2$, пусть $c_2 \in B_c(a_2 \in A_2, R_2)$. Тогда*

$$(3.1.9) \quad f(c_1)(c_2) \in B_c(f(a_1)(a_2), \|f\|(R_1 C_2 + C_1 R_2))$$

*где*

$$(3.1.10) \quad C_1 = \|a_1\|_1 + R_1 \quad C_2 = \|a_2\|_2 + R_2$$

Доказательство. Согласно определению 2.1.14

$$(3.1.11) \quad \begin{aligned} \|c_1 - a_1\|_1 &\leq R_1 \\ \|c_2 - a_2\|_2 &\leq R_2 \end{aligned}$$

Неравенство

$$(3.1.12) \quad \|c_i\| = \|a_i + c_i - a_i\| \leq \|a_i\| + \|c_i - a_i\| \leq \|a_i\| + R_i$$

следует из неравенства (3.1.11) и утверждения 2.1.9.3. Неравенство

$$(3.1.13) \quad C_i \leq \max(\|a_i\|, \|a_i\| + R_i) = \|a_i\| + R_i \quad i = 1, 2$$

следует из неравенства (3.1.12) и равенства (3.1.5). Равенство (3.1.10) следует из неравенства (3.1.13). Утверждение (3.1.9) следует из неравенств (3.1.4), (3.1.11), равенства (3.1.10) и определения 2.1.14. $\square$

Теорема 3.1.5. *Пусть*

$$f : A_1 \dashrightarrow A_2$$

*представление $\Omega_1$-группы $A_1$ с нормой $\|x\|_1$ в $\Omega_2$-группе $A_2$ с нормой $\|x\|_2$. Пусть последовательность $A_1$-чисел $a_{1 \cdot n}$ сходится и*

$$(3.1.14) \quad \lim_{n \to \infty} a_{1 \cdot n} = a_1$$

*Пусть последовательность $A_2$-чисел $a_{2 \cdot n}$ сходится и*

$$(3.1.15) \quad \lim_{n \to \infty} a_{2 \cdot n} = a_2$$

*Тогда последовательность $A_2$-чисел $f(a_{1 \cdot n})(a_{2 \cdot n})$ сходится и*

$$(3.1.16) \quad \lim_{n \to \infty} f(a_{1 \cdot n})(a_{2 \cdot n}) = f(a_1)(a_2)$$

Доказательство. Пусть

$$\delta_1 \in R, \ \delta_1 > 0 \tag{3.1.17}$$

Из равенства (3.1.14) и теоремы 2.1.18 следует, что для заданного $\delta_1$ существует $N_1$ такое, что из условия $n > N_1$ следует

$$a_{1 \cdot n} \in B_o(a_1, \delta_1) \tag{3.1.18}$$

Из равенства (3.1.15) и теоремы 2.1.18 следует, что для заданного $\delta_1$ существует $N_2$ такое, что из условия $n > N_2$ следует

$$a_{2 \cdot n} \in B_o(a_2, \delta_1) \tag{3.1.19}$$

Пусть
$$N = \max(N_1, N_2)$$

Из равенств (3.1.18), (3.1.19), теоремы 3.1.4 и условия $n > N$ следует, что

$$f(a_{1 \cdot n})(a_{2 \cdot n}) \in B_o(f(a_1)(a_n), \epsilon_1) \tag{3.1.20}$$

$$\epsilon_1 = \|f\|(\delta_1 C_2 + \delta_1 C_1) \tag{3.1.21}$$

где

$$C_i = \|a_i\| + \delta_1 \tag{3.1.22}$$

Из равенства (3.1.22) и утверждений 2.1.9.1, (3.1.17) следует, что

$$C_i > 0 \quad \frac{dC_i}{d\delta_1} > 0 \tag{3.1.23}$$

Из равенств (3.1.21), (3.1.22) и утверждения (3.1.23) следует, что $\epsilon_1$ является полиномиальной строго монотонно возрастающей функцией $\delta_1 > 0$. При этом

$$\delta_1 = 0 \Rightarrow \epsilon_1 = 0$$

Согласно теореме 2.3.5, отображение (3.1.21) отображает интервал $[0, \delta_1)$ в интервал $[0, \epsilon_1)$. Согласно теореме 2.3.3, для заданного $\epsilon > 0$ существует такое $\delta > 0$, что

$$\epsilon_1(\delta) < \epsilon$$

Согласно построению, значение $N$ зависит от значения $\delta_1$. Мы выберем значение $N$, соответствующее $\delta_1 = \delta$. Следовательно, для заданного $\epsilon \in R, \epsilon > 0$, существует $N$ такое, что из условия $N > N$ следует

$$f(a_{1 \cdot n})(a_{2 \cdot n}) \in B_o(f(a_1)(a_n), \epsilon) \tag{3.1.24}$$

Равенство (3.1.16) следует из равенства (3.1.24) и теоремы 2.1.18. $\square$

Теорема 3.1.6. *Представление*

$$f : A_1 \overset{*}{\longrightarrow} A_2$$

$\Omega_1$-*группы* $A_1$ *с нормой* $\|x\|_1$ *в* $\Omega_2$-*группе* $A_2$ *с нормой* $\|x\|_2$ *может быть продолжено до представления*

$$f^* : A_1 \overset{*}{\longrightarrow} M(X, A_2)$$

$\Omega_1$-*группы* $A_1$ *в* $\Omega_2$-*группе* $M(X, A_2)$ *где* $(g \in M(X, A_2))$

$$(f^*(a_1)(g))(x) = f(a_1)(g(x)) \tag{3.1.25}$$

Доказательство. Чтобы доказать теорему, мы должны показать, что отображение $f^*$ является гомоморфизмом $\Omega_1$-группы $A_1$ и для любого $a_1 \in A_1$ отображение $f^*(a_1)$ является гомоморфизмом $\Omega_2$-группы $M(X, A_2)$.

Согласно определению [7]-2.1.4, отображение $f$ является гомоморфизмом $\Omega_1$-группы $A_1$. Следовательно, для $n$-арной операции $\omega \in \Omega_1$ и произвольных $a_{1\cdot 1}, ..., a_{1\cdot n} \in A_1$ следующее равенство верно

$$(3.1.26) \qquad f(a_{1\cdot 1}...a_{1\cdot n}\omega) = f(a_{1\cdot 1})...f(a_{1\cdot n})\omega$$

Из равенств (2.6.1), (3.1.26) следует, что

$$(3.1.27) \qquad f(a_{1\cdot 1}...a_{1\cdot n}\omega)(a_2) = (f(a_{1\cdot 1})(a_2))...(f(a_{1\cdot n})(a_2))\omega$$

Из равенств (3.1.25), (3.1.27) следует, что

$$(3.1.28) \qquad \begin{aligned}(f^*(a_{1\cdot 1}...a_{1\cdot n}\omega)(g))(x) &= f(a_{1\cdot 1}...a_{1\cdot n}\omega)(g(x)) \\ &= f(a_{1\cdot 1})(g(x))...f(a_{1\cdot n})(g(x))\omega \\ &= ((f^*(a_{1\cdot 1})(g))(x))...((f^*(a_{1\cdot n})(g))(x))\omega\end{aligned}$$

Из равенств (2.6.1), (3.1.28) следует, что

$$(3.1.29) \qquad \begin{aligned}(f^*(a_{1\cdot 1}...a_{1\cdot n}\omega)(g))(x) &= ((f^*(a_{1\cdot 1})(g))...(f^*(a_{1\cdot n})(g))\omega)(x) \\ &= ((f^*(a_{1\cdot 1})...f^*(a_{1\cdot n})\omega)(g))(x)\end{aligned}$$

Из равенства (3.1.29) следует, что

$$(3.1.30) \qquad f^*(a_{1\cdot 1}...a_{1\cdot n}\omega) = f^*(a_{1\cdot 1})...f^*(a_{1\cdot n})\omega$$

Из равенства (3.1.30) следует, что отображение $f^*$ является гомоморфизмом $\Omega_1$-группы $A_1$.

Согласно определению [7]-2.1.4, отображение $f(a_1)$, $a_1 \in A_1$, является гомоморфизмом $\Omega_2$-группы $A_2$. Следовательно, для $n$-арной операции $\omega \in \Omega_2$ и произвольных $a_{2\cdot 1}, ..., a_{2\cdot n} \in A_2$ следующее равенство верно

$$(3.1.31) \qquad f(a_1)(a_{2\cdot 1}...a_{2\cdot n}\omega) = (f(a_1)(a_{2\cdot 1}))...(f(a_1)(a_{2\cdot n}))\omega$$

Из равенств (3.1.25), следует, что

$$(3.1.32) \qquad (f^*(a_1)(g_1...g_n\omega))(x) = f(a_1)((g_1...g_n\omega)(x))$$

Из равенств (2.6.1), (3.1.32) следует, что

$$(3.1.33) \qquad (f^*(a_1)(g_1...g_n\omega))(x) = f(a_1)(g_1(x)...g_n(x)\omega)$$

Из равенств (3.1.31), (3.1.33) следует, что

$$(3.1.34) \qquad (f^*(a_1)(g_1...g_n\omega))(x) = (f(a_1)(g_1(x)))...(f(a_1)(g_n(x)))\omega$$

Из равенств (2.6.1), (3.1.25), (3.1.34) следует, что

$$(3.1.35) \qquad \begin{aligned}(f^*(a_1)(g_1...g_n\omega))(x) &= ((f^*(a_1)(g_1))(x))...((f^*(a_1)(g_n))(x))\omega \\ &= ((f^*(a_1)(g_1))...(f^*(a_1)(g_n))\omega)(x)\end{aligned}$$

Из равенства (3.1.35) следует, что

$$(3.1.36) \qquad f^*(a_1)(g_1...g_n\omega) = (f^*(a_1)(g_1))...(f^*(a_1)(g_n))\omega$$

Из равенства (3.1.36) следует, что для любого $a_1 \in A_1$ отображение $f^*(a_1)$ является гомоморфизмом $\Omega_2$-группы $M(X, A_2)$. $\square$

### 3.2. Представление Ω-группы отображений

Теорема 3.2.1. *Представление*

$$f : A_1 \longrightarrow\!\!\!* A_2$$

$\Omega_1$-*группы* $A_1$ *с нормой* $\|x\|_1$ *в* $\Omega_2$-*группе* $A_2$ *с нормой* $\|x\|_2$ *порождает представление*

$$f_X : M(X, A_1) \longrightarrow\!\!\!* M(X, A_2)$$

$\Omega_1$-*группы* $M(X, A_1)$ *в* $\Omega_2$-*группе* $M(X, A_2)$ *где ( $g_1 \in M(X, A_1)$, $g_2 \in M(X, A_2)$ )*

$$\begin{aligned}f_X(g_1)(g_2) &: X -> A_2 \\ (f_X(g_1)(g_2))(x) &= f(g_1(x))(g_2(x))\end{aligned} \tag{3.2.1}$$

Доказательство. Чтобы доказать теорему, мы должны показать, что отображение $f_X$ является гомоморфизмом $\Omega_1$-группы $M(X, A_1)$ и для любого $g_1 \in M(X, A_1)$ отображение $f_X(g_1)$ является гомоморфизмом $\Omega_2$-группы $M(X, A_2)$.

Согласно определению [7]-2.1.4, отображение $f$ является гомоморфизмом $\Omega_1$-группы $A_1$. Следовательно, для $n$-арной операции $\omega \in \Omega_1$ и произвольных $a_{1\cdot 1}, ..., a_{1\cdot n} \in A_1$ следующее равенство верно

$$f(a_{1\cdot 1}...a_{1\cdot n}\omega) = f(a_{1\cdot 1})...f(a_{1\cdot n})\omega \tag{3.2.2}$$

Из равенства (2.6.1) следует, что

$$(g_{1\cdot 1}...g_{1\cdot n}\omega)(x) = g_{1\cdot 1}(x)...g_{1\cdot n}(x)\omega \tag{3.2.3}$$

Из равенств (3.2.2), (3.2.3) следует, что

$$f((g_{1\cdot 1}...g_{1\cdot n}\omega)(x)) = f(g_{1\cdot 1}(x)...g_{1\cdot n}(x)\omega) = f(g_{1\cdot 1}(x))...f(g_{1\cdot n}(x))\omega \tag{3.2.4}$$

Из равенств (2.6.1), (3.2.4), следует, что

$$\begin{aligned}f((g_{1\cdot 1}...g_{1\cdot n}\omega)(x))(a_2) &= (f(g_{1\cdot 1}(x))...f(g_{1\cdot n}(x))\omega)(a_2) \\ &= (f(g_{1\cdot 1}(x))(a_2))...(f(g_{1\cdot n}(x))(a_2))\omega\end{aligned} \tag{3.2.5}$$

Из равенств (3.2.1), (3.2.5) следует, что

$$\begin{aligned}(f_X(g_{1\cdot 1}...g_{1\cdot n}\omega)(g_2))(x) &= f((g_{1\cdot 1}...g_{1\cdot n}\omega)(x))(g_2(x)) \\ &= (f(g_{1\cdot 1}(x))(g_2(x)))...(f(g_{1\cdot n}(x))(g_2(x)))\omega \\ &= ((f_X(g_{1\cdot 1})(g_2))(x))...((f_X(g_{1\cdot n})(g_2))(x))\omega\end{aligned} \tag{3.2.6}$$

Из равенств (2.6.1), (3.2.6) следует, что

$$\begin{aligned}(f_X(g_{1\cdot 1}...g_{1\cdot n}\omega)(g_2))(x) &= ((f_X(g_{1\cdot 1})(g_2))...(f_X(g_{1\cdot n})(g_2))\omega)(x) \\ &= ((f_X(g_{1\cdot 1})...f_X(g_{1\cdot n})\omega)(g_2))(x)\end{aligned} \tag{3.2.7}$$

Из равенства (3.2.7) следует, что

$$f_X(g_{1\cdot 1}...g_{1\cdot n}\omega) = f_X(g_{1\cdot 1})...f_X(g_{1\cdot n})\omega \tag{3.2.8}$$

Из равенства (3.2.8) следует, что отображение $f_X$ является гомоморфизмом $\Omega_1$-группы $M(X, A_1)$.

Согласно определению [7]-2.1.4, отображение $f(a_1)$, $a_1 \in A_1$, является гомоморфизмом $\Omega_2$-группы $A_2$. Следовательно, для $n$-арной операции $\omega \in \Omega_2$ и произвольных $a_{2\cdot 1}, ..., a_{2\cdot n} \in A_2$ следующее равенство верно

$$(3.2.9) \qquad f(a_1)(a_{2\cdot 1}...a_{2\cdot n}\omega) = (f(a_1)(a_{2\cdot 1}))...(f(a_1)(a_{2\cdot n}))\omega$$

Из равенства (2.6.1) следует, что

$$(3.2.10) \qquad (g_{2\cdot 1}...g_{2\cdot n}\omega)(x) = g_{2\cdot 1}(x)...g_{2\cdot n}(x)\omega$$

Из равенств (3.2.9), (3.2.10) следует, что

$$(3.2.11) \qquad \begin{aligned} f(a_1)((g_{2\cdot 1}...g_{2\cdot n}\omega)(x)) &= f(a_1)(g_{2\cdot 1}(x)...g_{2\cdot n}(x)\omega) \\ &= (f(a_1)(g_{2\cdot 1}(x)))...(f(a_1)(g_{2\cdot n}(x)))\omega \end{aligned}$$

Из равенств (3.2.1), (3.2.11) следует, что

$$(3.2.12) \qquad \begin{aligned} (f_X(g_1)(g_{2\cdot 1}...g_{2\cdot n}\omega))(x) &= f(g_1(x))((g_{2\cdot 1}...g_{2\cdot n}\omega)(x)) \\ &= (f(g_1(x))(g_{2\cdot 1}(x)))...(f(g_1(x))(g_{2\cdot n}(x)))\omega \\ &= ((f_X(g_1)(g_{2\cdot 1}))(x))...((f(g_1)(g_{2\cdot n}))(x))\omega \end{aligned}$$

Из равенств (2.6.1), (3.2.12) следует, что

$$(3.2.13) \qquad (f_X(g_1)(g_{2\cdot 1}...g_{2\cdot n}\omega))(x) = ((f_X(g_1)(g_{2\cdot 1}))...(f_X(g_1)(g_{2\cdot n}))\omega)(x)$$

Из равенств (3.2.9), (3.2.13) следует, что

$$(3.2.14) \qquad f_X(g_1)(g_{2\cdot 1}...g_{2\cdot n}\omega) = (f_X(g_1)(g_{2\cdot 1}))...(f_X(g_1)(g_{2\cdot n}))\omega$$

Из равенства (3.2.14) следует, что, для любого $g_1 \in M(X, A_1)$, отображение $f_X(g_1)$ является гомоморфизмом $\Omega_2$-группы $M(X, A_2)$. $\square$

Теорема 3.2.2. *Пусть*

$$f : A_1 \relbar\joinrel\ast\joinrel\rightarrow A_2$$

*представление полной $\Omega_1$-группы $A_1$ с нормой $\|x\|_1$ в полной $\Omega_2$-группе $A_2$ с нормой $\|x\|_2$. Пусть последовательность отображений $g_{1\cdot n} \in M(X, A_1)$, $n = 1, ...,$ сходится равномерно к отображению $g_1$. Пусть последовательность отображений $g_{2\cdot n} \in M(X, A_2)$, $n = 1, ...,$ сходится равномерно к отображению $g_2$. Пусть множество значений отображения $g_i$, $i = 1, 2$, компактно. Тогда последовательность отображений $f_X(g_{1\cdot n})(g_{2\cdot n})$ сходится равномерно к отображению $f_X(g_1)(g_2)$.*

Доказательство. Из равенства

$$(3.2.15) \qquad g_i(x) = \lim_{n \to \infty} g_{i\cdot n}(x)$$

и теоремы 2.6.4 следует, что для заданного

$$(3.2.16) \qquad \delta_1 \in R, \ \delta_1 > 0$$

существует $N_i$ такое, что из условия $n > N_i$ следует

$$(3.2.17) \qquad g_{i\cdot n}(x) \in B_o(g_i(x), \delta_1)$$

Пусть

$$N = \max(N_1, N_2)$$

Из равенства (3.2.17), теоремы 3.1.4 и условия $n > N$ следует, что

$$(3.2.18) \qquad f(g_{1\cdot n}(x))(g_{2\cdot n}(x)) \in B_o(f(g_1(x))(g_2(x)), \epsilon_1)$$

$$(3.2.19) \quad \epsilon_1 = \|f\|(\delta_1 C_2 + \delta_1 C_1)$$

где

$$(3.2.20) \quad C_i = \|f_i(x)\| + \delta_1$$

Значение $C_i$ зависит от $x$. Поэтому $\epsilon_1$ также зависит от $x$. Для оценки равномерной сходимости нам необходимо выбрать максимальное значение $\epsilon_1$. Для заданного $\delta_1$, значение $\epsilon_1$ является полилинейным отображением значений $C_i$. Следовательно, $\epsilon_1$ имеет максимальное значение, когда каждое $C_i$ имеет максимальное значение. Так как мы рассматриваем оценку значения $\epsilon_1$ справа, для нас неважно одновременно ли $C_i$ приобретают максимальные значения. Так как множество значений отображения $f_i$ компактно, то, согласно теореме 2.3.11, определена следующая величина

$$(3.2.21) \quad F_i = \sup \|f_i(x)\|$$

Из равенства (3.2.21) и утверждения 2.1.9.1 следует, что

$$(3.2.22) \quad F_i \geq 0$$

Из равенств (3.2.20), (3.2.21) следует, что мы можем положить

$$(3.2.23) \quad C_i = F_i + \delta_1$$

Из равенства (3.2.23) и утверждений (3.2.22), (3.2.16) следует, что

$$(3.2.24) \quad C_i > 0 \quad \frac{dC_i}{d\delta_1} > 0$$

Из равенств (3.2.19), (3.2.23) и утверждения (3.2.24) следует, что $\epsilon_1$ является полиномиальной строго монотонно возрастающей функцией $\delta_1 > 0$. При этом

$$\delta_1 = 0 \Rightarrow \epsilon_1 = 0$$

Согласно теореме 2.3.5, отображение (3.2.19) отображает интервал $[0, \delta_1)$ в интервал $[0, \epsilon_1)$. Согласно теореме 2.3.3, для заданного $\epsilon > 0$ существует такое $\delta > 0$, что

$$\epsilon_1(\delta) < \epsilon$$

Согласно построению, значение $N$ зависит от значения $\delta_1$. Мы выберем значение $N$, соответствующее $\delta_1 = \delta$. Следовательно, для заданного $\epsilon \in R$, $\epsilon > 0$, существует $N$ такое, что из условия $n > N$ следует

$$(3.2.25) \quad f(g_{1 \cdot n}(x))(g_{2 \cdot n}(x)) \in B_o(f(g_1(x))(g_2(x)), \epsilon)$$

Равенство

$$f(g_1(x))(g_2(x)) = \lim_{n \to \infty} f(g_{1 \cdot n}(x))(g_{2 \cdot n}(x))$$

следует из равенства (3.2.25) и теоремы 2.6.4. $\square$

### 3.3. Пополнение представления $\Omega$-группы

Определение 3.3.1. Пусть

$$f : A_1 \relbar\joinrel\twoheadrightarrow A_2$$

представление нормированной $\Omega_1$-группы $A_1$ с нормой $\|x\|_1$ в нормированную $\Omega_2$-группе $A_2$ с нормой $\|x\|_2$. Представление

$$g : B_1 \relbar\joinrel\twoheadrightarrow B_2$$

полной $\Omega_1$-группы $B_1$ в полной $\Omega_2$-группе $B_2$ называется **пополнением представления** $f$, если

3.3.1.1: $\Omega$-группа $B_1$ является пополнением нормированной $\Omega_1$-группы $A_1$.
3.3.1.2: $\Omega$-группа $B_2$ является пополнением нормированной $\Omega_2$-группы $A_2$.
3.3.1.3: Если $a_1 \in A_1$, $a_2 \in A_2$, то

$$(3.3.1) \quad g(a_1)(a_2) = f(a_1)(a_2)$$

□

Теорема 3.3.2. *Пусть*

$$f : A_1 \relbar\joinrel\ast\joinrel\rightarrow A_2$$

*представление нормированной $\Omega_1$-группы $A_1$ с нормой $\|x\|_1$ в нормированную $\Omega_2$-группе $A_2$ с нормой $\|x\|_2$. Пусть представления*

$$g_1 : B_{11} \relbar\joinrel\ast\joinrel\rightarrow B_{12}$$

$$g_2 : B_{21} \relbar\joinrel\ast\joinrel\rightarrow B_{22}$$

*являются пополнениями представления $f$.*

3.3.2.1: *Для $i = 1$, $2$, существует изоморфизм $\Omega$-группы*

$$(3.3.2) \quad h_i : B_{i1} \to B_{i2}$$

*такой, что*

$$(3.3.3) \quad h_i(a_i) = a_i \quad a_i \in A_i$$

$$(3.3.4) \quad \|h_i(b_i)\|_i = \|b_i\|_i$$

3.3.2.2: *Пара отображений*

$$(3.3.5) \quad (h_1 : B_{11} \to B_{12}, h_2 : B_{21} \to B_{22})$$

*является морфизмом представлений из $g_1$ в $g_2$.*

Доказательство. Утверждение 3.3.2.1 является следствием теоремы 2.5.2.

Пусть $b_1 \in B_{11}$, $b_2 \in B_{12}$. Согласно теореме 2.2.14 существуют последовательности $A_i$-чисел $a_{i \cdot n}$ такая, что

$$(3.3.6) \quad b_1 = \lim_{n \to \infty} a_{1 \cdot n} \quad b_2 = \lim_{n \to \infty} a_{2 \cdot n}$$

Из равенства (3.3.6) следует, что последовательность $A_I$-чисел $a_{i \cdot n}$ является фундаментальной последовательностью в $\Omega$-группе $A_i$. Для любого $n$ следующее равенство очевидно

$$(3.3.7) \quad f(a_{1 \cdot n})(a_{2 \cdot n}) = f(a_{1 \cdot n})(a_{2 \cdot n})$$

Так как $f(a_{1 \cdot n})(a_{2 \cdot n}) \in A_2$, то равенство

$$(3.3.8) \quad f(a_{1 \cdot n})(a_{2 \cdot n}) = h_2(g_1(a_{1 \cdot n})(a_{2 \cdot n})) = g_2(h_1(a_{1 \cdot n}))(h_2(a_{2 \cdot n}))$$

является следствием равенства (3.3.3). Непрерывность отображений $h_1$, $h_2$ следует из равенства (3.3.4). Следовательно, равенство

$$(3.3.9) \quad h_2(g_1(b_1)(b_2)) = g_2(h_1(b_1))(h_2(b_2))$$

следует из равенств (3.3.6), (3.3.8) и теоремы 3.1.5. Теорема следует из равенств (3.3.9), [7]-(2.2.4) и определения [7]-2.2.2. □

Лемма 3.3.3. Пусть

$$f : A_1 \dashrightarrow A_2$$

представление нормированной $\Omega_1$-группы $A_1$ с нормой $\|x\|_1$ в нормированную $\Omega_2$-группе $A_2$ с нормой $\|x\|_2$. Пусть $a_{1\cdot n}$ является фундаментальной последовательностью $A_1$-чисел. Пусть $a_{2\cdot n}$ является фундаментальной последовательностью $A_2$-чисел. Тогда последовательность $f(a_{1\cdot n})(a_{2\cdot n})$ является фундаментальной последовательностью $A_2$-чисел.

Доказательство. Пусть

(3.3.10) $$\delta_1 \in R, \ \delta_1 > 0$$

Так как последовательность $a_{1\cdot n}$ является фундаментальной последовательностью $A_1$-чисел, то согласно теореме 2.1.20 следует, что для заданного $\delta_1$ существует, зависящее от $\delta_1$, натуральное число $N_1$ такое, что

(3.3.11) $$a_{1\cdot q} \in B_o(a_{1\cdot p}, \delta_1)$$

для любых $p, q > N_1$. Так как последовательность $a_{2\cdot n}$ является фундаментальной последовательностью $A_2$-чисел, то согласно теореме 2.1.20, следует, что для заданного $\delta_1$ существует, зависящее от $\delta_1$, натуральное число $N_2$ такое, что

(3.3.12) $$a_{2\cdot q} \in B_o(a_{2\cdot p}, \delta_1)$$

для любых $p, q > N_2$. Пусть

$$N = \max(N_1, N_2)$$

Из утверждений (3.3.11), (3.3.12) и теоремы 3.1.4 следует, что

(3.3.13) $$f(a_{1\cdot q})(a_{2\cdot q}) \in B_o(f(a_{1\cdot p})(a_{2\cdot p}), \epsilon_1)$$

(3.3.14) $$\epsilon_1 = \|f\|(\delta_1 C_2 + \delta_1 C_1)$$

где

(3.3.15) $$C_i = \|a_{i\cdot p}\|_i + \delta_1$$

Из равенства (3.3.15) и утверждений 2.1.9.1, (3.3.10) следует, что

(3.3.16) $$C_i > 0 \quad \frac{dC_i}{d\delta_1} > 0$$

Из равенств (3.3.14), (3.3.15) и утверждения (3.3.16) следует, что $\epsilon_1$ является полиномиальной строго монотонно возрастающей функцией $\delta_1 > 0$. При этом

$$\delta_1 = 0 \Rightarrow \epsilon_1 = 0$$

Согласно теореме 2.3.5, отображение (3.3.14) отображает интервал $[0, \delta_1)$ в интервал $[0, \epsilon_1)$. Согласно теореме 2.3.3, для заданного $\epsilon > 0$ существует такое $\delta > 0$, что

$$\epsilon_1(\delta) < \epsilon$$

Согласно построению, значение $N$ зависит от значения $\delta_1$. Мы выберем значение $N$, соответствующее $\delta_1 = \delta$. Следовательно, для заданного $\epsilon \in R, \epsilon > 0$, существует $N$ такое, что из условия $n > N$ следует

(3.3.17) $$f(a_{1\cdot q})(a_{2\cdot q}) \in B_o(f(a_{1\cdot p})(a_{2\cdot p}), \epsilon)$$

Из утверждения (3.3.17) и теоремы 2.1.20 следует, что последовательность $f(a_{1\cdot n})(a_{2\cdot n})$ является фундаментальной последовательностью. □

Теорема 3.3.4. *Пополнение представления нормированной Ω-группы существует.*

Доказательство. Пусть $A_1$ - нормированная $\Omega_1$-группа. Пусть $A_2$ - нормированная $\Omega_2$-группа. Согласно теореме 2.5.17, существует пополнение $B_1$ нормированной $\Omega_1$-группы $A_1$ и пополнение $B_2$ нормированной $\Omega_2$-группы $A_2$.

Чтобы доказать существование пополнения $g$ представления $f$

$$g : B_1 \longrightarrow_* B_2$$

мы будем пользоваться записью, предложенной в разделе 2.5. Рассмотрим отношением эквивалентности на множестве фундаментальных последовательностей $A_i$-чисел

$$(3.3.18) \qquad a_{i\cdot n} \sim b_{i\cdot n} \Leftrightarrow \lim_{n \to \infty}(a_{i\cdot n} - b_{i\cdot n}) = 0$$

Пусть $B_i$ является множеством классов эквивалентных фундаментальных последовательностей $A_i$-чисел. Мы будем пользоваться записью

$$(3.3.19) \qquad [a_{i\cdot n}] = \{b_{i\cdot n} : a_{i\cdot n} \sim b_{i\cdot n}\}$$

для класса последовательностей, эквивалентных последовательности $a_{i\cdot n}$.

Лемма 3.3.5. *Пусть $a_{1\cdot n}$, $b_{1\cdot n}$ являются фундаментальными последовательностями $A_1$-чисел,*

$$(3.3.20) \qquad a_{1\cdot n} \sim b_{1\cdot n}$$

*Пусть $a_{2\cdot n}$, $b_{2\cdot n}$ являются фундаментальными последовательностями $A_2$-чисел,*

$$(3.3.21) \qquad a_{2\cdot n} \sim b_{2\cdot n}$$

*Тогда*

$$(3.3.22) \qquad f(a_{1\cdot n})(a_{2\cdot n}) \sim f(b_{1\cdot n})(b_{2\cdot n})$$

Доказательство. Пусть

$$(3.3.23) \qquad \delta_1 \in R, \ \delta_1 > 0$$

Из утверждения (3.3.20) и леммы 2.5.4 следует, что для заданного $\delta_1$ существует, зависящее от $\delta_1$, натуральное число $N_1$ такое, что

$$(3.3.24) \qquad b_{1\cdot n} \in B_o(a_{1\cdot n}, \delta_1)$$

для любого $n > N_1$. Из утверждения (3.3.21) и теоремы 2.5.4 следует, что для заданного $\delta_1$ существует, зависящее от $\delta_1$, натуральное число $N_2$ такое, что

$$(3.3.25) \qquad b_{2\cdot n} \in B_o(a_{2\cdot n}, \delta_1)$$

для любого $n > N_2$. Пусть

$$N = \max(N_1, N_2)$$

Из утверждений (3.3.24), (3.3.25) и теоремы 3.1.4 следует, что

$$(3.3.26) \qquad f(b_{1\cdot n})(b_{2\cdot n}) \in B_o(f(a_{1\cdot n})(a_{2\cdot n}), \epsilon_1)$$

$$(3.3.27) \qquad \epsilon_1 = \|f\|(\delta_1 C_2 + \delta_1 C_1)$$

где

$$(3.3.28) \qquad C_i = \|a_{i\cdot n}\|_i + \delta_1$$

Из равенства (3.3.28) и утверждений 2.1.9.1, (3.3.23) следует, что

$$(3.3.29) \qquad C_i > 0 \quad \frac{dC_i}{d\delta_1} > 0$$

Из равенств (3.3.27), (3.3.28) и утверждения (3.3.29) следует, что $\epsilon_1$ является строго монотонно возрастающей функцией $\delta_1 > 0$. При этом

$$\delta_1 = 0 \Rightarrow \epsilon_1 = 0$$

Согласно теореме 2.3.5, отображение (3.3.27) отображает интервал $[0, \delta_1)$ в интервал $[0, \epsilon_1)$. Согласно теореме 2.3.3, для заданного $\epsilon > 0$ существует такое $\delta > 0$, что

$$\epsilon_1(\delta) < \epsilon$$

Согласно построению, значение $N$ зависит от значения $\delta_1$. Мы выберем значение $N$, соответствующее $\delta_1 = \delta$. Следовательно, для заданного $\epsilon \in R$, $\epsilon > 0$, существует $N$ такое, что из условия $n > N$ следует

$$(3.3.30) \qquad f(b_{1 \cdot n})(b_{2 \cdot n}) \in B_o(f(a_{1 \cdot n})(a_{2 \cdot n}), \epsilon)$$

Утверждение (3.3.22) следует из утверждения (3.3.30) и леммы 2.5.4. $\odot$

Согласно леммам 3.3.3, 3.3.5 отображение

$$(3.3.31) \qquad g([a_{1 \cdot n}])([a_{2 \cdot n}]) = [f(a_{1 \cdot n})(a_{2 \cdot n})]$$

определено корректно.

Пусть $a_1 \in A_1$, $a_2 \in A_2$. Согласно теореме 2.2.13, существует последовательность $A_i$-чисел $a_{i \cdot n}$, $n = 1, ...,$ такая, что $a_i = [a_{i \cdot n}]$. Согласно теореме 3.1.5

$$(3.3.32) \qquad f(a_1)(a_2) = \lim_{n \to \infty} f(a_{1 \cdot n})(a_{2 \cdot n})$$

Из равенств (3.3.31), (3.3.32) следует, что отображение $g$ удовлетворяет утверждению 3.3.1.3

$$g(a_1)(a_2) = g([a_{1 \cdot n}])([a_{2 \cdot n}]) = f(a_1)(a_2)$$

Лемма 3.3.6. *Отображение $g$ является гомоморфизмом $\Omega_1$-группы $B_1$.*

Доказательство. Пусть $a_{1.1}, ..., a_{1.n} \in B_1$. Согласно теореме 2.2.14, существует последовательность $A_1$-чисел $a_{1 \cdot i \cdot m}$, $m = 1, ...,$ такая, что

$$(3.3.33) \qquad a_{1 \cdot i} = [a_{1 \cdot i \cdot m}] = \lim_{m \to \infty} a_{1 \cdot i \cdot m}$$

Пусть $a_2 \in B_2$. Согласно теореме 2.2.14, существует последовательность $A_2$-чисел $a_{2 \cdot m}$, $m = 1, ...,$ такая, что

$$(3.3.34) \qquad a_2 = [a_{2 \cdot m}] = \lim_{m \to \infty} a_{2 \cdot m}$$

Согласно определению [7]-2.1.4, отображение $f$ является гомоморфизмом $\Omega_1$-группы $A_1$. Следовательно, для $n$-арной операции $\omega \in \Omega_1$ и любого $m$ следующее равенство верно

$$(3.3.35) \qquad f(a_{1 \cdot 1 \cdot m} ... a_{1 \cdot n \cdot m} \omega) = f(a_{1 \cdot 1 \cdot m}) ... f(a_{1 \cdot n \cdot m}) \omega$$

Из равенств (2.6.1), (3.3.35), следует, что

$$(3.3.36) \qquad \begin{aligned} f(a_{1 \cdot 1 \cdot m} ... a_{1 \cdot n \cdot m} \omega)(a_{2 \cdot m}) &= (f(a_{1 \cdot 1 \cdot m}) ... f(a_{1 \cdot n \cdot m}) \omega)(a_{2 \cdot m}) \\ &= (f(a_{1 \cdot 1 \cdot m})(a_{2 \cdot m})) ... (f(a_{1 \cdot n \cdot m})(a_{2 \cdot m})) \omega \end{aligned}$$

Равенство

$$(3.3.37) \qquad \lim_{m\to\infty} a_{1\cdot 1\cdot m}...a_{1\cdot n\cdot m}\omega = a_{1\cdot 1}...a_{1\cdot n}\omega$$

следует из равенства (3.3.33) и теоремы 2.4.6. Равенство

$$(3.3.38) \qquad \lim_{m\to\infty} f(a_{1\cdot 1\cdot m}...a_{1\cdot n\cdot m}\omega)(a_{2\cdot m}) = g(a_{1\cdot 1}...a_{1\cdot n}\omega)(a_2)$$

следует из равенств (3.3.31), (3.3.34), (3.3.37) и теоремы 3.1.5. Равенство

$$(3.3.39) \qquad \lim_{m\to\infty} f(a_{1\cdot i\cdot m})(a_{2\cdot m}) = g(a_{1\cdot i})(a_2)$$

следует из равенств (3.3.31), (3.3.33), (3.3.34) и теоремы 3.1.5. Равенство

$$(3.3.40) \begin{aligned} \lim_{m\to\infty} (f(a_{1\cdot 1\cdot m})(a_{2\cdot m}))...(f(a_{1\cdot n\cdot m})(a_{2\cdot m}))\omega \\ = (g(a_{1\cdot 1})(a_2))...(g(a_{1\cdot n})(a_2))\omega \end{aligned}$$

следует из равенства (3.3.39) и теоремы 2.4.6. Равенство

$$(3.3.41) \qquad g(a_{1\cdot 1}...a_{1\cdot n}\omega)(a_2) = (g(a_{1\cdot 1})(a_2))...(g(a_{1\cdot n})(a_2))\omega$$

следует из равенств (3.3.36), (3.3.38), (3.3.40). Из равенства (3.3.41) следует, что отображение $g$ является гомоморфизмом $\Omega_1$-группы $B_1$. $\odot$

Лемма 3.3.7. Для любого $a_1 \in B_1$ отображение $g(a_1)$ является гомоморфизмом $\Omega_2$-группы $B_2$.

Доказательство. Пусть $a_{2\cdot 1}, ..., a_{2\cdot n} \in B_2$. Согласно теореме 2.2.14, существует последовательность $A_2$-чисел $a_{2\cdot i\cdot m}$, $m = 1, ...,$ такая, что

$$(3.3.42) \qquad a_{2\cdot i} = [a_{2\cdot i\cdot m}] = \lim_{m\to\infty} a_{2\cdot i\cdot m}$$

Пусть $a_1 \in B_1$. Согласно теореме 2.2.14, существует последовательность $A_1$-чисел $a_{1\cdot m}$, $m = 1, ...,$ такая, что

$$(3.3.43) \qquad a_1 = [a_{1\cdot m}] = \lim_{m\to\infty} a_{1\cdot m}$$

Согласно определению [7]-2.1.4, отображение $f(a_{1\cdot m})$, $a_{1\cdot m} \in A_1$, является гомоморфизмом $\Omega_2$-группы $A_2$. Следовательно, для $n$-арной операции $\omega \in \Omega_2$ и произвольного $m$ следующее равенство верно

$$(3.3.44) \qquad f(a_{1\cdot m})(a_{2\cdot 1\cdot m}...a_{2\cdot n\cdot m}\omega) = (f(a_{1\cdot m})(a_{2\cdot 1\cdot m}))...(f(a_{1\cdot m})(a_{2\cdot n\cdot m}))\omega$$

Равенство

$$(3.3.45) \qquad \lim_{m\to\infty} a_{2\cdot 1\cdot m}...a_{2\cdot n\cdot m}\omega = a_{2\cdot 1}...a_{2\cdot n}\omega$$

следует из равенства (3.3.42) и теоремы 2.4.6. Равенство

$$(3.3.46) \qquad \lim_{m\to\infty} f(a_{1\cdot m})(a_{2\cdot 1\cdot m}...a_{2\cdot n\cdot m}\omega) = g(a_1)(a_{2\cdot 1}...a_{2\cdot n}\omega)$$

следует из равенств (3.3.31), (3.3.43), (3.3.45) и теоремы 3.1.5. Равенство

$$(3.3.47) \qquad \lim_{m\to\infty} f(a_{1\cdot m})(a_{2\cdot i\cdot m}) = g(a_1)(a_{2\cdot i})$$

следует из равенств (3.3.31), (3.3.42), (3.3.43) и теоремы 3.1.5. Равенство

$$(3.3.48) \begin{aligned} \lim_{m\to\infty} (f(a_{1\cdot m})(a_{2\cdot 1\cdot m}))...(f(a_{1\cdot m})(a_{2\cdot n\cdot m}))\omega \\ = (g(a_1)(a_{2\cdot 1}))...(g(a_1)(a_{2\cdot n}))\omega \end{aligned}$$

следует из равенства (3.3.47) и теоремы 2.4.6. Равенство

$$(3.3.49) \qquad g(a_1)(a_{2\cdot 1}...a_{2\cdot n}\omega) = (g(a_1)(a_{2\cdot 1}))...(g(a_1)(a_{2\cdot n}))\omega$$

следует из равенств (3.3.44), (3.3.46), (3.3.48). Из равенства (3.3.49) следует, что, для любого $a_1 \in B_1$, отображение $g(a_1)$ является гомоморфизмом $\Omega_2$-группы $B_2$. ⊙

Из лемм 3.3.6, 3.3.7 и определения [7]-2.1.4 следует, что отображение $g$, определённое равенством (3.3.31), является представлением полной $\Omega_1$-группы $B_1$ в полной $\Omega_2$-группе $B_2$. □

## 3.4. Ω-кольцо

Пусть

$$f : A_1 \relbar\joinrel\twoheadrightarrow A_2$$

представление $\Omega_1$-группы $A_1$ с нормой $\|x\|_1$ в $\Omega_2$-группе $A_2$ с нормой $\|x\|_2$. Согласно определение [7]-2.1.4, отображение $f$ является гомоморфизмом $\Omega_1$-группы $A_1$. Следовательно, если произведение определено в $\Omega_1$-группе $A_1$, то следующее равенство верно

$$(3.4.1) \qquad f(a_{1\cdot 1}a_{1\cdot 2}) = f(a_{1\cdot 1})f(a_{1\cdot 2})$$

Однако, если произведение не определено в $\Omega_1$-группе $A_1$, то мы можем определить произведение в $\Omega_1$-группе $A_1$, пользуясь равенством (3.4.1). Это даёт основание полагать, что в $\Omega_1$-группе $A_1$ определено произведение.

Так как отображение $f$ является гомоморфизмом $\Omega_1$-группы $A_1$, то следующие равенства верны

$$(3.4.2) \quad \begin{aligned} f(a_1(b_{1\cdot 1} + b_{1\cdot 2}))(a_2) &= f(a_1)(f(b_{1\cdot 1} + b_{1\cdot 2})(a_2)) \\ &= f(a_1)(f(b_{1\cdot 1})(a_2) + f(b_{1\cdot 2})(a_2)) \\ &= f(a_1)(f(b_{1\cdot 1})(a_2)) + f(a_1)(f(b_{1\cdot 2})(a_2)) \\ &= f(a_1 b_{1\cdot 1})(a_2) + f(a_1 b_{1\cdot 2})(a_2) \\ &= f(a_1 b_{1\cdot 1} + a_1 b_{1\cdot 2})(a_2) \end{aligned}$$

$$(3.4.3) \quad \begin{aligned} f((a_{1\cdot 1} + a_{1\cdot 2})b_1)(a_2) &= f(a_{1\cdot 1} + a_{1\cdot 2})(f(b_1)(a_2)) \\ &= f(a_{1\cdot 1})(f(b_1)(a_2)) + f(a_{1\cdot 2})(f(b_1)(a_2)) \\ &= f(a_{1\cdot 1} b_1)(a_2) + f(a_{1\cdot 2}f(b_1)(a_2) \\ &= f(a_{1\cdot 1} b_1 + a_{1\cdot 2} b_1)(a_2) \end{aligned}$$

Если представление $f$ эффективно, то равенства

$$(3.4.4) \qquad a_1(b_{1\cdot 1} + b_{1\cdot 2}) = a_1 b_{1\cdot 1} + a_1 b_{1\cdot 2}$$

$$(3.4.5) \qquad (a_{1\cdot 1} + a_{1\cdot 2})b_1 = a_{1\cdot 1} b_1 + a_{1\cdot 2} b_1$$

являются следствием равенств (3.4.2), (3.4.3). Следовательно, умножение дистрибутивно по отношению к сложению. Таким образом, $\Omega_1$-группа $A_1$ является кольцом относительно сложения и умножения.

Определение 3.4.1. Ω-группа, в которой определено произведение, называется Ω-**кольцом**. □

Замечание 3.4.2. Бикольцо (определение [6]-2.6) является примером Ω-группы, в которой определены две операции произведение. Это связано с двумя различными представлениями алгебры матриц с элементами из некоммутативной алгебры. Соответственно, в бикольце определены две структуры Ω-кольца.

В алгебре октонионов мы наблюдаем другую картину. Так как произведение в алгебре октонионов неассоциативно, то рассматриваемое представление алгебры октонионов определяет иную операцию произведения. □

Определение 3.4.3. Эффективное представление Ω-кольца $A$ в абелевой группе называется $A$**-модулем**.[3.2] □

Я рассматриваю представление в абелевой группе, но не в Ω-группе. Это связано с тем, что произвольные операции в универсальной алгебре усложняют условие линейной зависимости. В тоже время, в универсальной алгебре, порождающей представление, для нас важны только операции сложения и умножения.

---

[3.2]Смерти также определение [9]-2.1.2.

Глава 4

# Список литературы

[1] S. Burris, H.P. Sankappanavar, A Course in Universal Algebra, Springer-Verlag (March, 1982),
eprint [http://www.math.uwaterloo.ca/ snburris/htdocs/ualg.html](http://www.math.uwaterloo.ca/ snburris/htdocs/ualg.html)
(The Millennium Edition)

[2] Г. Е. Шилов. Математический анализ, Функции одного переменного, части 1 - 2, М., Наука, 1969

[3] А. Н. Колмогоров, С. В. Фомин. Элементы теории функций и функционального анализа. М., Наука, 1976

[4] W.A. Coppel, Number Theory: An Introduction to Mathematics,
Springer, 2009

[5] Michael J. Field, Differential Calculus and Its Applications,
Dover Publications, 2012; ISBN-13: 978-0486497952

[6] Александр Клейн, Бикольцо матриц,
eprint [arXiv:math.OA/0612111](arXiv:math.OA/0612111) (2007)

[7] Александр Клейн, Представление универсальной алгебры,
eprint [arXiv:0912.3315](arXiv:0912.3315) (2010)

[8] Александр Клейн, Линейные отображения свободной алгебры,
eprint [arXiv:1003.1544](arXiv:1003.1544) (2010)

[9] Александр Клейн, Свободная алгебра со счётным базисом,
eprint [arXiv:1211.6965](arXiv:1211.6965) (2012)

[10] В. И. Смирнов, Курс высшей математики, том первый.
М., Наука, 1974

[11] П. Кон, Универсальная алгебра, М., Мир, 1968

[12] Н. Бурбаки, Общая топология, основные структуры, перевод с французского Д. А. Райкова, М. Наука, 1968

[13] Н. Бурбаки, Общая топология. Использование вещественных чисел в общей топологии.
перевод с французского С. Н. Крачковского под редакцией Д. А. Райкова,
М. Наука, 1975

[14] V. I. Arnautov, S. T. Glavatsky, A. V. Mikhalev,
Introduction to the theory of topological rings and modules, Volume 1995,
Marcel Dekker, Inc, 1996

# Глава 5

# Предметный указатель

# Глава 6

# Специальные символы и обозначения